\documentclass[10pt,final]{article}

\usepackage{arxiv}
\usepackage{mathpazo}

\usepackage{a4}
\usepackage{amsmath}%
\usepackage{amstext}%
\usepackage{amssymb}%
\usepackage{amsthm}
\usepackage{comment}
\usepackage{listings}

\usepackage[final]{graphicx}
\usepackage{epstopdf}
\usepackage[justification=centering]{caption}
\usepackage{subcaption}

\usepackage{diagbox}
\usepackage{multirow}
\usepackage{tabulary, capt-of}

\usepackage{hyperref}

\hypersetup{
   colorlinks=true,
   urlbordercolor={1 1 1},
   linkcolor={blue},
   citecolor={green}
}

\usepackage{mathtools}
\newtagform{blue}{\color{red}(}{)}
\newtagform{redandblue}[\textcolor{red}]{\color{red}(}{)}

\usetagform{blue}

\newtheorem{theorem}{Theorem}

\newtheorem{corollary}[theorem]{Corollary}

\newtheorem{example}[theorem]{Example}

\newtheorem{lemma}[theorem]{Lemma}

\newtheorem{proposition}[theorem]{Proposition}
\newtheorem{rem}[theorem]{Remark}
\newenvironment{remark}{\rem\rm}{\endrem}

\newcounter{unnumber}

\newtheorem{prf}[unnumber]{Proof}
\newenvironment{Proof}{\prf\rm}{\hfill{$\blacksquare$}\endprf}%

\DeclareMathOperator*\argmin{argmin}

\author{Cristian Daniel Alecsa \thanks{\textbf{Corresponding author. This work was supported by a grant of Ministry of Research and Innovation, CNCS - UEFISCDI, project number PN-III-P1-1.1-TE-2016-0266, within PNCDI III.}} \\
Tiberiu Popoviciu Institute of Numerical Analysis \\
Romanian Academy \\
Cluj-Napoca, Romania \\
\texttt{cristian.alecsa@ictp.acad.ro} \\ \\
Romanian Institute of Science and Technology \\
Cluj-Napoca, Romania \\
\texttt{alecsa@rist.ro}}

\title{The long time behavior and the rate of convergence of symplectic convex algorithms obtained via splitting discretizations of inertial damping systems.}

\begin{document}
\maketitle

\noindent \textbf{Abstract.} 
In this paper we propose new numerical algorithms in the setting of unconstrained optimization problems and we study the rate of convergence in the iterates of the objective function. Furthermore, our algorithms are based upon splitting and symplectic methods and they preserve the energy properties of the inherent continuous dynamical system that contains a Hessian perturbation. In particular, we show that Nesterov gradient method is equivalent to a Lie-Trotter splitting applied to a Hessian driven damping system.  Finally, some numerical experiments are presented in order to validate the theoretical results.
\vspace{1ex}

 \vspace{0.5ex}

\section{Preliminaries}\label{S1}
The theory of accelerated optimization algorithms initiated from the work of Polyak (see \cite{Polyak}) and it peaked with the convex gradient method developed by Nesterov in \cite{Nesterov}. The theory of dynamical systems associated to Nesterov accelerated method was introduced first in \cite{Boyd} by Su, Boyd and Cand\` es, where they pointed out that the method they have developed was informal. Also for the case of convex objective functions, the first true model for Nesterov discretization was considered very recently by Muehlebach and Jordan in \cite{Jordan} and in the generalization given in \cite{ALP} by Alecsa, L\' aszl\' o and Pin\c ta, where under the gradient was considered a perturbation of the exact solution of the dynamical system. On the other hand, also recently, in \cite{Defazio} Defazio showed that the accelerated gradient method can be seen as a proximal point method in a Riemannian manifold, but the theoretical analysis was done under the assumption of strong convexity. The first contribution of the present paper is to give a rigorous alternative to these approaches and show that Nesterov method is in fact a Lie-Trotter discretization of a Hessian-driven damping dynamical system (see Theorem \ref{T2} and Subsection \ref{S21}). Our second important contribution is to show that many well-known optimization algorithms are in fact a combination between splitting and symplectic methods and we further show that the algorithms based upon Hessian-driven damping discretizations (similar to those in \cite{Attouch_Chbani}) have a continuous dynamical system that is also constructed with respect to the aforementioned methods (see Subsection \ref{S22} and Theorem \ref{T4}). On the other hand, we introduce some generalized optimization algorithms associated to Hessian-driven damping dynamical systems that have many advantages from a practical point of view(see \ref{LT-S-IGAHD} and Section \ref{S4}), and this represents our third and last contribution of the present article. \\ \\
We consider the following unconstrained optimization problem 
\begin{align}\label{Opt_Pb}
\min\limits_{x \in \mathbb{R}^N} f(x)
\end{align}
where $f : \mathbb{R}^N \to \mathbb{R}$ is a convex and Fr\' echet differentiable function, for which its gradient $\nabla f : \mathbb{R}^N \to \mathbb{R}^N$ is $L$-Lipschitz, i.e. there exists $L > 0$ such that for each $x,y \in \mathbb{R}^N$
\begin{align}\label{Lipschitz}
\| \nabla f(x) - \nabla f(y) \| \leq L \| x-y \| \ .
\end{align}
For simplicity, we consider that all of the second order equations that appear in this paper are coupled with some initial conditions $x(t_0) = x_0 \in \mathbb{R}^N$ and $\dot{x}(t_0) = v_0 \in \mathbb{R}^N$, respectively (where $t_0 > 0$). This assumption is imposed only for the clarity of the paper, since (as in \cite{Attouch_Chbani} and other similar papers), we are interested in the asymptotic behavior of a trajectory of the underlying dynamical system. On the other hand, this is also true for the algorithmic counterparts, where, given $x_0 \in \mathbb{R}^N$ and taking $x_1 \in \mathbb{R}^N$ arbitrary, the accelerated inertial methods start from $n=1$ in order to compute $x_2$. \\
Now, we give a brief presentation regarding the discrete and the continuous analysis of optimization problems. In \cite{Polyak}, Polyak has introduced the heavy ball with friction dynamical system, defined as
\begin{align}\label{HBF}\tag{HBF$_\gamma$}
\ddot{x}(t) + \gamma \dot{x}(t) = - \nabla f(x(t)) \ .
\end{align}
Here $\gamma > 0$ is a fixed viscous damping parameter that represents the coefficient of friction of the dynamical system (\ref{HBF}). Furthermore, the combination between the second derivative of the solution and $\gamma$ is used in order to accelerate the gradient methods using inertial coefficients. For the (\ref{HBF}) the objective function $f$ does not decrease along the trajectories. On the other hand, the global energy (kinetic + potential) $\mathcal{E}(t) = \frac{1}{2} \| \dot{x}(t) \|^2 + f(x(t))$ is monotone decreasing, by being a continuous Lyapunov-type functional. Furthermore, we remind that Attouch. et. al. (in \cite{Attouch_Goudou}) have studied the asymptotic properties of this dynamical system.  \\
The natural discretization of (\ref{HBF}) is Polyak's inertial method (see \cite{Polyak}), i.e.
\begin{align}\label{PIM}\tag{PIM$_\gamma$}
\begin{cases}
y_n = x_n + (1-h \gamma)(x_n - x_{n-1}) \\
x_{n+1} = y_n - h^2 \nabla f(x_n)
\end{cases}
\end{align}
where, in contrast with descent methods, it has an additional inertial term $x_n - x_{n-1}$. If we divide this difference between two consecutive iteration values by the stepsize $h$, then we obtain a basic approximation to the first order derivative of the solution $x(t)$ of (\ref{HBF}). As a side note, we mention that in \cite{Sun}, the authors made the first non-ergodic analysis of the rate of convergence for a variant of Polyak's algorithm (\ref{PIM}), that contains a non-constant stepsize, related to the inertial term, under the assumption that the objective function is coercive. Now, we shift our attention to the most popular inertial algorithm, developed by Y. Nesterov in \cite{Nesterov}. The following update rule is often called  Nesterov's second accelerated method and it takes the following form :
\begin{align}\label{AGM2}\tag{AGM2$_\alpha$}
\begin{cases}
y_n = x_n + \alpha_n (x_n - x_{n-1}) \\
x_{n+1} = y_n - h^2 \nabla f(y_n)
\end{cases}
\end{align}
where, in practical applications, $\alpha_n = \dfrac{n}{n+\alpha}$ (see the work of Chambolle and Dossal \cite{Dossal}). A different set of update rules is given by the so-called Nesterov's first accelerated method (briefly AGM1), which gives the best convergence rate in the objective function, and it is defined by a sequence $(\lambda_n)_{n \in \mathbb{N}}$ that has the general term defined recursively as $\lambda_n = \left( 1 + \sqrt{1 + 4 \lambda^{2}_{n-1}} \right)/ 2$. The equivalence between (AGM1) and (\ref{AGM2}) for $\alpha = 2$ has been used in \cite{Bansal} regarding a Lyapunov analysis of Nesterov's algorithm. Also, in \cite{Tseng} some unified generalizations were made using Bregman distances. In contrast with Polyak's method (\ref{PIM}), in \cite{Nesterov} it was shown that Nesterov's algorithm has an improved rate of convergence, namely $O \left( 1/n^2 \right)$. Further, we point out that, in a nonconvex setting, a two-step inertial method similar to (\ref{AGM2}) has been recently studied in \cite{Laszlo}. In \cite{Boyd}, Su et. al. introduced the continuous counterpart of Nesterov's dynamical system of the form
\begin{align}\label{AVD}\tag{AVD$_\alpha$}
\ddot{x}(t) + \dfrac{\alpha}{t} \dot{x}(t) = - \nabla f(x(t))
\end{align}
The dynamical system with asymptotic vanishing damping was introduced in the framework of convex optimization problems. For $\alpha \geq 3$, it was shown that each trajectory of (\ref{AVD}) satisfies the rate of convergence :
$ f(x(t)) - \min f = O \left( 1/t^2 \right) \text{ as } t \to \infty $.
The relationship between Nesterov's algorithm (\ref{AGM2}) and the dynamical system (\ref{AVD}) is based upon vanishing stepsize arguments, and so the Nesterov's optimization scheme can not be considered a proper discretization of the associated differential equation (even the authors of \cite{Boyd} pointed out that their derivation of the dynamical system is informal). A general version of \ref{AVD} can be found in the seminal papers \cite{Cabot}, \cite{Gadat} and in their extensions given in \cite{AttouchCabot_AS}, where $\dfrac{\alpha}{t}$ is replaced by a general function $a = a(t)$. \\ 
In \cite{Bot}, a perturbed variant of (\ref{AVD}) was analyzed in the setting of nonconvex minimization problems, using the Kurdyka-Lojasiewicz property imposed on a regularization of the objective function.
Very recently, in \cite{Jordan}, it was shown that Nesterov convex gradient method is the natural discretization of a modified second order dynamical system. Further, a generalized version of the one from \cite{Jordan} was studied in \cite{ALP}, in which the authors proved an $O \left( 1/t^2 \right)$ rate of convergence as $t \to \infty$ for the trajectories of the dynamical system.
that contains a perturbation under the gradient (see \cite{ALP}):
\begin{align}\label{Extended-AVD}\tag{Extended-AVD}
\ddot{x}(t) + \dfrac{\alpha}{t} \dot{x}(t) = - \nabla f \left( x \left( t \right) + \left( \gamma - \dfrac{\alpha \gamma^2}{t} \right) \dot{x}(t) \right) 
\end{align}
On the other hand, for this continuous version of Nesterov's algorithm some theoretical results were derived related to the regularity of the solution of the associated initial value problem. \\
Now, we shift our attention to generalizations of the aforementioned dynamical systems, in which the derivative of the solution is perturbed with the value of the Hessian along the trajectories. In \cite{Attouch_Chbani}, Attouch et. al. introduced the Hessian-driven damping dynamical system, i.e.
\begin{align}\label{DIN-AVD}\tag{DIN-AVD$_{\tfrac{\alpha}{t},\beta(t),b(t)}$}
\ddot{x}(t) + \dfrac{\alpha}{t} \dot{x}(t) + \beta(t) \nabla^2 f(x(t)) \dot{x}(t) = -b(t) \nabla f(x(t)) \ .
\end{align}
Furthermore, the inertial system with constant damping that represents a generalization of the heavy ball with friction was first studied in \cite{Alvarez}. Also see \cite{Attouch_Redont} for the case of $\left( \text{DIN-AVD}_{\tfrac{\alpha}{t},\beta,1} \right)$. The idea behind the Hessian dynamics is the fact that the oscillations that are present both in the continuous and in the discrete case are neutralized in the presence of convex objective functions. Now, we present some particular cases that play a key role in the present paper. The Hessian system $\left( \text{DIN-AVD}_{\tfrac{3}{t}, 0,1} \right)$ was considered in \cite{Oberman}, where M. Laborde and A.M. Oberman pointed out that Nesterov's discretization is in fact a forward Euler discretization for the Hessian dynamical system. In fact, they showed the method is a forward Euler discretization with a predictor value under the gradient, but not a true discretization. So, in this paper, we claim that Nesterov accelerated gradient method is indeed a proper discretization scheme of $\left( \text{DIN-AVD}_{\tfrac{\alpha}{t},\beta,1 + \tfrac{\beta}{t}} \right)$, namely a full Lie-Trotter splitting discretization composed of forward Euler methods. 
Attouch et. al. proposed a Nesterov-like algorithm but with two consecutive gradients that represent a discretization of the Hessian. This is linked with the term $\nabla^2 f(x(t)) \dot{x}(t)$ that acts like an additional dissipation term. We shall show a link between the true discretization of a modified \ref{DIN-AVD} type system that is a Polyak-type method and the algorithm introduced by Attouch et. al. 
On the other hand, another particular case is $\left( \text{DIN-AVD}_{\alpha,\beta,1} \right)$, where in \cite{Castera}, the forward Euler method was constructed as an optimizer to deep learning neural networks using this dynamical system. At the same time, this discretization scheme is a forward Euler applied to a first order dynamical system that is equivalent to the second order differential equation, taking account the regularity of $f$, i.e.
\begin{align}
\begin{cases}
\dot{x}(t) + \beta \nabla f(x(t)) + \left( \alpha - \dfrac{1}{\beta} \right) x(t) + \dfrac{1}{\beta} v(t) = 0 \\
\dot{v}(t) + \left( \alpha - \dfrac{1}{\beta} \right) x(t) + \dfrac{1}{\beta} v(t) = 0
\end{cases}
\end{align}
We mention that the particular case of the vanishing damping Hessian system $\left( \text{DIN-AVD}_{\tfrac{3}{t},\beta, 1 + \tfrac{\beta}{t}} \right)$ can also be written as a first order dynamical system. The following form is much simpler than the former first order system. This will be used frequently in our section of main results.
\begin{align}\label{DIN-AVD-NESTEROV}
\begin{cases}
\dot{x}(t) = \dfrac{2}{t} \left( v \left ( t \right) -x \left( t \right) \right) - \beta \nabla f(x(t)) \\
\dot{v}(t) = - \dfrac{t}{2} \nabla f(x(t))
\end{cases}
\end{align}
Finally, we end this section by pointing out that, for $\beta = 0$, the dynamical system (\ref{DIN-AVD-NESTEROV}) reduces to system that is equivalent to the second order evolution equation (\ref{AVD}) with $\alpha = 3$. For other details, we refer to \cite{Oberman}.

\subsection{Some notes on splitting operator techniques}\label{S11}

In this subsection we present the most used splitting techniques for a general autonomous evolution equation of the form
\begin{align}\label{Ev_Eq}
\dfrac{d}{dt} \varphi(t) = F(\varphi(t))
\end{align}
where $\varphi : \mathbb{R}^N \to \mathbb{R}$, $\varphi = \varphi(t)$ is the flow of the dynamical system. If the vector field from the right hand side can be decomposed into two vector fields, i.e. $F = F^{[1]} + F^{[2]}$, then we attach the associated sub-problems of (\ref{Ev_Eq}), namely
\begin{align}
(P_{1}) :
\dfrac{d}{dt} \varphi(t) = F^{[1]}(\varphi(t))
\quad \text{ and } \quad (P_{2}) :
\dfrac{d}{dt} \varphi(t) = F^{[2]}(\varphi(t))
\end{align}
Let $\varphi^{[1]} = \varphi^{[1]}(t)$ be the solution of $(P_{1})$ and $\varphi^{[2]} = \varphi^{[2]}(t)$ be the solution of $(P_{2})$. Then, the Lie-Trotter splitting or sequential splitting (see \cite{Trotter}) after one stepsize $h$ is :
\begin{align}\label{Lie-Trotter-1}\tag{LTS$_1$}
\psi_{h} = \varphi^{[1]}_{h} \circ \varphi^{[2]}_{h}
\end{align}
Another Lie-Trotter splitting operator is to first apply the exact solution of $(P_1)$ as an initial condition to $(P_2)$, i.e.
\begin{align}\label{Lie-Trotter-2}\tag{LTS$_2$}
\psi^{\ast}_{h} = \varphi^{[2]}_{h} \circ \varphi^{[1]}_{h}
\end{align}
Here we have employed the usual notation: $\varphi^{[i]}_{t} := \varphi^{[i]}(t)$ for $i \in \lbrace 1, 2 \rbrace$.
For a discrete time $t_n = nh$, (\ref{Lie-Trotter-1}), becomes
$ \psi(t_n) = (\varphi^{[1]} \circ \varphi^{[2]})(t_n) = \left( \varphi^{[1]}_{h} \circ \varphi^{[2]}_{h} \right)^n (h) $
and (\ref{Lie-Trotter-2}) is equivalent with
$ \psi^{\ast} (t_n) = \left( \varphi^{[2]}_{h} \circ \varphi^{[1]}_{h} \right)^n (h)$.
A more advanced type of splitting technique is the so-called Strang or Marchuk splitting (see \cite{Marchuk} and \cite{Strang}). From a numerical point of view, this has an advandage over the Lie-Trotter splitting by the fact that is second-order accurate. Again, giving a stepsize $h$, the Strang splitting at $t = h$, is defined as
\begin{align}\label{Strang-1}\tag{SS$_1$}
\psi_{h} = \varphi^{[1]}_{h/2} \circ \varphi^{[2]}_{h} \circ \varphi^{[1]}_{h/2}
\end{align}
Again, the Strang splitting operator can appear into a different symmetric form, namely
\begin{align}\label{Strang-2}\tag{SS$_2$}
\psi^{\ast}_{h} = \varphi^{[2]}_{h/2} \circ \varphi^{[1]}_{h} \circ \varphi^{[2]}_{h/2}
\end{align}
Applying recursively, we obtain the splitting operator formula for the discrete time $t_n = nh$. Quite interestingly, if, for example, the exact flows $\varphi^{[1]}(t)$ and $\varphi^{[2]}(t)$ are not available, then one can use consistent approximations to these solutions. In this manner, from now on, $\tilde{\varphi}^{[1]}_{h}$ and $\tilde{\varphi}^{[2]}_{h}$ represents numerical discretization of the flows on the two sub-problems $(P_1)$ and $(P_2)$ at time $t_n = nh$, respectively. Then, for example the Lie-Trotter splitting (\ref{Lie-Trotter-1}) can be defined as
\begin{align}
\psi_{h} = \tilde{\varphi}^{[1]}_{h} \circ \tilde{\varphi}^{[2]}_{h}
\end{align}
This is valid also for the Strang splitting operator formula (see \cite{Sanz-Serna}). On the other hand, for the case of multiple vector fields, when $F = F^{[1]} + \ldots + F^{[m]}$, we consider $m$ sub-problems, i.e.
\begin{align}
(P_{1}) :
\dfrac{d}{dt} \varphi(t) = F^{[1]}(\varphi(t))
\quad \ldots\ldots \quad (P_{m}) :
\dfrac{d}{dt} \varphi(t) = F^{[m]}(\varphi(t))
\end{align}
In this case, the Lie-Trotter splitting becomes (see \cite{Hairer_Book})
\begin{align}\label{Multiple-LT}
\psi_{h} = \varphi^{[1]}_{h} \circ \ldots \circ \varphi^{[m]}_{h}
\end{align}
For example, since in the next sections we will use the case of $m = 3$, we observe that the formula (\ref{Multiple-LT}) is in fact 
\begin{align}
\psi_{h} = \left( \varphi^{[1]}_{h} \circ \varphi^{[2]}_{h} \right) \circ \varphi^{[3]}_{h}
\end{align}
This is equivalent to the fact that we divide the evolution equation into two-subproblems and we use a numerical approximation to the solution of the second sub-problem. Furthermore, the first sub-problem will also be divided into another two-subproblems that can also be decomposed according to the Lie-Trotter splitting. \\
The case of unbounded linear operator has been studied in \cite{HO}, when $F(\varphi(t)) = A \varphi(t) + B \varphi(t)$. Furthermore, in \cite{HKO}, Hansen et. al. considered the splitting for the case of semilinear evolution equations, when one of the sub-problems is not linear, i.e. $F(\varphi(t)) = A \varphi(t) + F^{[2]}(\varphi(t))$. For both this cases, Lie-Trotter splitting and Strang-Marchuk splitting ca be defined with the above formulas, using the exact flows or some approximations to the exact solution on the sub-problems. \\
Now, we turn our attention to the case of non-autonomous evolution equations of the form
\begin{align}\label{NonAut_Ev_Eq}
\dfrac{d}{dt} \varphi(t) = F(t, \varphi(t))
\end{align}
There are some approaches for this kind of problem. In \cite{Batkai}, Batkai et. al. have used a similar approach as Hansen, Kramer and Ostermann. They have considered the case when $F = F^{[1]} + F^{[2]}$. In their paper, they applied the Lie-Trotter splitting, when $F^{[1]}(\varphi(t)) = A(t) \varphi(t)$ and $F^{[2]}(t) = B(t) \varphi(t)$. For brevity, we recall the sequential splitting (or Lie-Trotter splitting) on the the interval of discretization $(t_{n-1}, t_{n}]$ :
\begin{align}\label{LT-NonAut}
\begin{cases}
\dfrac{d}{dt} \varphi^{[1]}(t) = B(t) \varphi^{[1]}(t) \quad , \quad t \in (t_{n-1},t_n] \quad \text{ with } \varphi^{[1]}(t_{n-1}) = \varphi^{LT}(t_{n-1}) \\
\dfrac{d}{dt} \varphi^{[2]}(t) = A(t) \varphi^{[2]}(t) \quad , \quad t \in (t_{n-1},t_n] \quad \text{ with } \varphi^{[2]}(t_{n-1}) = \varphi^{[1]}(t_{n-1}) 
\end{cases}
\end{align}
where $\varphi^{LT}(t_{n-1})$ represent the previous discretization solution that was applied to $(t_{n-2}, t_{n-1}]$. Also, from the above numerical scheme we get that $\varphi^{LT}(t_n) = \varphi^{[2]}(t_n)$. On the other hand, taking $t = t_{n-1}$, we obtain the full discretization of Lie-Trotter splitting, given by the formula (\ref{Lie-Trotter-2}). In the following sections, we will employ this method, even though, Batkai et. al. have provided another formula when the evolution equation $\dfrac{d}{dt} \varphi(t) = A(r) \varphi(t) + B(r) \varphi(t)$ can be solved for each $r \in \mathbb{R}$. The only difference is that in (\ref{LT-NonAut}), we can evaluate $B$ exactly in the current discretization time element $t_n$. Since we consider the updated value of $t_n$ much more appropriately for the discretization schemes, we will use the previous formula. Also, we provide a reference for some splitting techniques that can  be applied to optimization problems. For this, see the work of Blanes et. al. \cite{Blanes}. We emphasize the fact that, recently, in \cite{APB} the authors have used theses splitting techniques in order to devise some optimizers used in neural network training. The analysis is in contrast with the present article by the fact that in \cite{APB}, the optimization algorithms were developed in an heuristic manner based upon algebraic manipulations. Nevertheless, this shows the applicability of the splitting-type methods to machine learning. Finally, we point out that splitting techniques have been used recently also in \cite{FrancaSplitting}, where the authors have employed balanced and rebalanced splitting methods in connection with the equilibrium states of some ODE's associated to proximal algorithms. This differs from our analysis by the fact that we use only the Lie-Trotter discretization scheme.

\subsection{A brief review of symplectic methods}\label{S12}

We consider the setting of autonomous evolution equation given by the formula (\ref{Ev_Eq}). Let $\varphi = (x, v) : \mathbb{R}^N \to \mathbb{R}^N$ and consider $H : \mathbb{R}^N \to \mathbb{R}$ to be the Hamiltonian. We write the first order differential equation as a Hamiltonian dynamical system, using $F$ from (\ref{Ev_Eq}) with $F(\varphi(t)) = ( \nabla_v H(x(t), v(t)), -  \nabla_x H(x(t), v(t)))^T$ in the following form (see \cite{Franca}):
\begin{align}\label{Hamiltonian}\tag{HS}
\begin{pmatrix}
\dfrac{d}{dt} x(t) \\ \\
\dfrac{d}{dt} v(t)
\end{pmatrix}
= 
\begin{pmatrix}
 \nabla_v H(x(t), v(t)) \\ \\
- \nabla_x H(x(t), v(t))
\end{pmatrix}
\end{align}
For the pioneering works on differential geometry and basic Hamiltonian systems, we refer to \cite{Feng} and \cite{Vogelaere}. Furthermore, for brevity, we recall (see \cite{Hairer_Book} and \cite{Hairer}) the symplectic Euler methods related to (\ref{Hamiltonian})
\begin{align}\label{SE_1}\tag{SE$_1$}
\begin{cases}
x_{n+1} = x_n + h \nabla_v H(x_{n+1}, v_n) \\
v_{n+1} = v_n - h \nabla_x H(x_{n+1}, v_n)
\end{cases}
\end{align}
respectively
\begin{align}\label{SE_2}\tag{SE$_2$}
\begin{cases}
x_{n+1} = x_n + h \nabla_v H(x_{n}, v_{n+1}) \\
v_{n+1} = v_n - h \nabla_x H(x_{n}, v_{n+1})
\end{cases}
\end{align}
Other symplectic-type algorithms can be constructed with combinations of the form $(x_n, v_{n+1})$ and $(x_{n+1}, v_n)$. Also from the work of Hairer, Lubich and Wanner \cite{Hairer_Book}, we recall a more advanced type of implicit-explicit symplectic integrators, namely the Stormer-Verlet schemes :
\begin{align}\label{SV_1}\tag{SV$_1$}
\begin{cases}
x_{n+1/2} = x_n + \dfrac{h}{2} \nabla_v H(x_{n+1/2}, v_n) \\
v_{n+1} = v_n - \dfrac{h}{2} \left( \nabla_x H(x_{n+1/2}, v_n) + \nabla_x H(x_{n+1/2}, v_{n+1}) \right) \\
x_{n+1} = x_{n+1/2} + \dfrac{h}{2} \nabla_v H(x_{n+1/2}, v_{n+1}) \\
\end{cases}
\end{align}
and a similar symmetric method
\begin{align}\label{SV_2}\tag{SV$_2$}
\begin{cases}
v_{n+1/2} = v_n - \dfrac{h}{2} \nabla_x H(x_n, v_{n+1/2}) \\
x_{n+1} = x_n + \dfrac{h}{2} \left( \nabla_v H(x_n, v_{n+1/2}) + \nabla_v H(x_{n+1}, v_{n+1/2}) \right) \\
v_{n+1} = v_{n+1/2} - \dfrac{h}{2} \nabla_x H(x_{n+1}, v_{n+1/2}) \\
\end{cases}
\end{align}
Fro brevity, we point out that $x_{n+1/2}$ from (\ref{SV_1}) and $v_{n+1/2}$ from (\ref{SV_2}) represent just a notation for intermediate steps that are used in the construction of the algorithms. Interestingly, under general assumptions, for a symplectic integrator of convergence order $\mu$, one has that $H(x_n,v_n) = H(x_0,v_0) + O(h^{\mu})$, where the constant does not depend on the iteration $n$ over exponentially long time intervals, namely $nh \leq e^{h_0/2h}$. This leads to a very good energy behavior of symplectic methods and so these methods represent a good choice for the discretization of a conservative Hamiltonian system. \\
Also, it is known that if the Hamiltonian is separable, i.e. $H(x,v) = T(v) + U(x)$, then (\ref{Hamiltonian}) becomes
\begin{align}\label{Part_Hamiltonian}
\begin{pmatrix}
\dfrac{d}{dt} x(t) \\ \\
\dfrac{d}{dt} v(t)
\end{pmatrix}
= 
\begin{pmatrix}
\nabla_v T(v(t)) \\ \\
- \nabla_x U(x(t))
\end{pmatrix}
\end{align}
Furthermore, dividing (\ref{Part_Hamiltonian}) into two sub-problems
\begin{align}
(P_{1}) :
\begin{cases}
\dfrac{d}{dt} x(t) = 0 \\ \\
\dfrac{d}{dt} v(t) = - \nabla_x U(x(t))
\end{cases}
\text{ and } \quad (P_{2}) :
\begin{cases}
\dfrac{d}{dt} x(t) = \nabla_v T(v(t)) \\ \\
\dfrac{d}{dt} v(t) = 0
\end{cases}
\end{align}
and applying the explicit forward Euler discretizations to both the sub-problems, then the symplectic Euler methods (\ref{SE_1}) and (\ref{SE_2}) lead to the Lie-Trotter full discretizations (\ref{Lie-Trotter-1}) and (\ref{Lie-Trotter-2}). The same analysis can be depicted for the Stormer Verlet methods (\ref{SV_1}) and (\ref{SV_2}), where these methods are equivalent to the Strang-Marchuk discretization schemes. \\
Recently, in \cite{Franca}, Fran\c ca et. al. devised some optimization schemes based on some Hamiltonian systems perturbed by a dissipation term $\gamma > 0$, of the form
\begin{align}\label{Hamiltonian_Dissip}\tag{HSD$_\gamma$}
\begin{pmatrix}
\dfrac{d}{dt} x(t) \\ \\
\dfrac{d}{dt} v(t)
\end{pmatrix}
= 
\begin{pmatrix}
\nabla_v H(x(t), v(t)) \\ \\
- \nabla_x H(x(t), v(t)) - \gamma v(t)
\end{pmatrix}
\end{align}
where the separable Hamiltonian $H : \mathbb{R}^N \to \mathbb{R}$ is defined such as $H(x,v) = \dfrac{\| v \|^2}{2m} + f(x)$, with $m > 0$ related to the inertial positive coefficient $\gamma$. Here, the first term $T(v) = \dfrac{\| v \|^2}{2m}$ is the kinetic energy and $U(x) = f(x)$ is the potential energy of the system. Furthermore, the heavy ball with friction system (\ref{HBF}) is a particular case of (\ref{Hamiltonian_Dissip}) and can be written as
\begin{align}\label{Hamiltonian_Poylak}
\begin{pmatrix}
\dfrac{d}{dt} x(t) \\ \\
\dfrac{d}{dt} v(t) 
\end{pmatrix}
= 
\begin{pmatrix}
v(t) \\ \\
- \nabla f(x(t)) - \gamma v(t)
\end{pmatrix}
\end{align}
and the Hamiltonian of the problem is a Lyapunov function and satisfies $\dfrac{d}{dt} H(x(t),v(t)) = - \gamma \| v(t) \|^2 \leq 0$ along the trajectories. Furthermore, it was shown that Polyak's method (\ref{PIM}) is a conformal symplectic Euler method for (\ref{Hamiltonian_Poylak}), but Nesterov's method is not. Also, Polyak's algorithm can be shown to be equivalent to a Lie-Trotter splitting (in connection with the decomposition of \ref{Hamiltonian_Dissip}), where for the conservative subproblem one applies a symplectic Euler method and for the dissipative subproblem an explicit forward Euler for the exact solution is used (for details see \cite{Franca}). Finally, we remind that for heavy ball system (\ref{HBF}) without friction (when $\gamma = 0$), the symplectic Euler and the symplectic Stormer-Verlet algorithms can be considered as the natural discretizations of the underlying ODE (see \cite{Hairer_Book} and \cite{Shi}).

\subsection{The organization of the paper}\label{S13}
The outline of the present paper is the following. In Subsection \ref{S21}, we show that Nesterov's algorithm can be viewed as a Lie-Trotter or sequential splitting for a Hessian driven damping system. Furthermore, we also present an alternative version of writing the momentum based methods with asymptotic vanishing damping that resembles the continuous dynamical system. Also, we present some well known results for the heavy ball system and the extended dynamical system with asymptotic vanishing damping in connection to symplectic Euler integrator.
In Subsection \ref{S23} and Subsection \ref{S24}, we consider some results regarding various algorithms linked to symplectic integrators for the extended AVD system. In Section \ref{S3}, we prove the rate of convergence of the Hessian-based algorithms in the framework of convex , differentiable objective functions, such that the gradient satisfies the Lipschitz continuity property. Finally, in Section \ref{S4}, we validate our theoretical results with some numerical experiments involving various choices concerning the parameters of our newly introduced algorithms in the framework of classical convex test functions. For conciseness, we point out a brief analysis that will be done in the following sequel. \\ \\
In Subsection \ref{S22}, we show some connections between the continuous version and the discrete version of some dynamical systems. Then dynamical system (\ref{Extended-AVD}) can be seen as a generalization of (\ref{AVD}) with an extra perturbation under the gradient that is similar to the updates in Nesterov's algorithm iterations. On the other hand, the (\ref{AVD}) dynamical system can be generalized through the Hessian driven damping evolution equation (\ref{DIN-AVD}), by using a perturbation term consisting on the product of the Hessian matrix and the gradient of the underlying solution. For conciseness, we remind that the difference between Polyak's method (\ref{PIM}) and (\ref{AGM2}) lies on the fact that the momentum method is a splitting method for a system similar to (\ref{AVD}) but with constant damping, while Nesterov's method is also a splitting based-discretization combined with a symplectic Euler method, but for the extended continuous system. Quite suprinsingly, our novelty was given in the second section concerning the fact that Nesterov's method can be seen also a discretization (by using a triple splitting technique) but for the Hessian driven damping system (\ref{DIN-AVD}) with $\beta(t) = \beta$ and $b(t) = 1 + \beta / t$. At the same time, regarding the Hessian-based dynamical system, as in the Nesterov's case, IGAHD algorithm (see \cite{Attouch_Chbani}), is also a splitting-based method where on the potential forces we apply the symplectic Euler method. The key role in our analysis is the fact that the symplectic methods are playing a special role in the discretizations of the continuous dynamical systems. We emphasize the fact that we have four types of symplectic methods, namely (\ref{SE_1}), (\ref{SE_2}), (\ref{SV_1}) and (\ref{SV_2}). By a trial and error approach, these methods, in combination with splitting discretizations, applied on the Hessian driven damping system (\ref{DIN-AVD}) lead to unfeasible algorithms due to the fact that the continuous evolution equation is of rather complicated nature. Also, using these techniques, one obtains algorithms containing at least with three gradient evaluations or with a product on the form $\nabla^2 f(x(t)) \cdot \nabla f(x(t))$ that is a computational burden (even after a forward-type discretization). So, our aim Subsection \ref{S23} and Subsection \ref{S24}, is to focus only on the (\ref{Extended-AVD}) dynamical system. We also have observed that the alternative variant of the symplectic Euler in contrast with the one used in \cite{Jordan} leads to the same empirical convergence as Nesterov's algorithm but with an additional gradient evaluation. On the other hand, Stormer-Verlet algorithms leads to three gradient evaluations, which is unfeasible. In this way, we will analyze in an intuitive manner only the Stormer-Verlet integrator (\ref{SV_2}) and the symplectic Euler integrator (\ref{SE_1}) that induces two gradient evaluations. This is similar to the two gradient scheme IGAHD. Moreover, we will introduce another dynamical system that for a Lie-Trotter discretization we will obtain the perturbed Nesterov-type algorithm introduced in \cite{Nguyen}.

\section{The discretizations of the dynamical systems}\label{S2}

\subsection{Nesterov's accelerated method is a Lie-Trotter splitting method applied to a Hessian-based dynamical system}\label{S21}

In \cite{Oberman}, it was shown that the system $\left( \text{DIN-AVD}_{\tfrac{3}{t},\beta,1 + \tfrac{\beta}{t}} \right)$ with the vanishing damping $3/t$ can also be written as a first order dynamical system, namely (\ref{DIN-AVD-NESTEROV}). For brevity, in our first result, we will show that the general vanishing damping system with Hessian perturbation $\left( \text{DIN-AVD}_{\tfrac{3}{t},\beta,1 + \tfrac{\beta}{t}} \right)$ can be written as a first order system. This result will be used also in this section when we will prove that Nesterov's method is actually a Lie-Trotter splitting applied to three sub-problems of this dynamical system.
\begin{proposition}\label{T1}
Let $f \in C^2(\mathbb{R}^n, \mathbb{R})$ and let $\beta \in \mathbb{R}$ and $\alpha > 1$. Further, consider the $\left( \text{DIN-AVD}_{\tfrac{\alpha}{t},\beta,1 + \tfrac{\beta}{t}} \right)$ second order evolution equation, i.e.
\begin{align}\label{Hessian_VD}
\ddot{x}(t) + \dfrac{\alpha}{t} \dot{x}(t) + \beta \nabla^2 f(x(t)) \dot{x}(t) = - \left( 1 + \dfrac{\beta}{t} \right) \nabla f(x(t))
\end{align}
On the other hand, for $t > 0$ consider the non-autonomous first order dynamical system :
\begin{align}\label{FirstOrder_VD}
\begin{cases}
\dot{x}(t) = \dfrac{\alpha - 1}{t} (v(t)-x(t)) - \beta \nabla f(x(t)) \\
\dot{v}(t) = - \dfrac{t}{\alpha - 1} \nabla f(x(t))
\end{cases}
\end{align}
Then, the dynamical systems (\ref{Hessian_VD}) and (\ref{FirstOrder_VD}) are equivalent, in the sense that each solution of the second order dynamical system is a solution of the first order system and viceversa.
\end{proposition}

\begin{Proof}
Let $(x(t),v(t)) \in \mathbb{R}^n \times \mathbb{R}^N$ be a solution for (\ref{FirstOrder_VD}). We will show that $x = x(t)$ also satisfies the Hessian driven damping system (\ref{Hessian_VD}). Firstly, using the fact that 
$$v(t) = x(t) + \dfrac{t}{\alpha - 1} \left( \dot{x}(t) + \beta \nabla f(x(t)) \right),$$
by differentiation we obtain that 
$$\dot{v}(t) = \dot{x}(t) + \dfrac{1}{\alpha - 1} \left( \dot{x}(t) + \beta \nabla f(x(t)) \right) + \dfrac{t}{\alpha - 1} \left( \ddot{x}(t) + \beta \nabla^2 f(x(t)) \dot{x}(t) \right) .$$
Since $\dot{v}(t) = - \dfrac{t}{\alpha-1} \nabla f(x(t))$, one obtains that 
$$ -\dfrac{t}{\alpha - 1} \nabla f(x(t)) = \dot{x}(t) + \dfrac{1}{\alpha - 1} \dot{x}(t) + \dfrac{\beta}{\alpha - 1} \nabla f(x(t)) + \dfrac{t}{\alpha-1} \ddot{x}(t) + \dfrac{\beta t}{\alpha-1} \nabla^2 f(x(t)) \dot{x}(t) \ , $$
which is equivalent to (\ref{Hessian_VD}). \\
On the other hand, let $x = x(t)$ be a solution to the Hessian damping system (\ref{Hessian_VD}). We will show that, if we define $v = v(t)$ as in the first equation of (\ref{FirstOrder_VD}), then the derivative of the function $v$ satisfies the second equation of (\ref{FirstOrder_VD}). In this sense, by the definition of $v = v(t)$, we have that 
$$ v(t) = x(t) + \dfrac{t}{\alpha - 1} \left( \dot{x}(t) + \beta \nabla f(x(t)) \right) \ .$$
By differentiation, as before, we obtain that
$$\dot{v}(t) = \dot{x}(t) + \dfrac{1}{\alpha - 1} \left( \dot{x}(t) + \beta \nabla f(x(t)) \right) + \dfrac{t}{\alpha - 1} \left( \ddot{x}(t) + \beta \nabla^2 f(x(t)) \dot{x}(t) \right) .$$
Using the fact that $x=x(t)$ satisfies (\ref{Hessian_VD}), it follows that
$$ \dot{v}(t) = \dot{x}(t) + \dfrac{1}{\alpha-1} \dot{x}(t) + \dfrac{\beta}{\alpha-1} \nabla f(x(t)) + \dfrac{t}{\alpha-1} \left( - \dfrac{\alpha}{t} \dot{x}(t) - \left( 1 + \dfrac{\beta}{t} \right) \nabla f(x(t)) \right) \ . $$
Finally, this equivalent exactly with the second line of (\ref{FirstOrder_VD}).
\end{Proof}

Now, the following theorem is the main result of this subsection. In \cite{Oberman}, the authors showed the equivalence between (\ref{AGM2}) and (\ref{NaG_Velocity}) for the case when $\alpha = 3$. Moreover, we emphasize that their method closely resembles the forward Euler method applied to the dynamical system (\ref{FirstOrder_VD}). This intuitive idea was pointed out in the authors work. On the other hand, since under the gradient we use a predictor $y_n$ and not the true value $x_n$ at the discretization time $t_n$, this method can not be interpreted as a proper discretization to the first order dynamical system. Following this idea, we will show that Nesterov's method (\ref{NaG_Velocity}) can be viewed as a discrete counterpart of (\ref{FirstOrder_VD}) using the idea of Lie-Trotter splitting. In fact, we show that Nesterov's accelerated gradient method is a Lie-Trotter splitting applied to a decomposition consisting of three subproblems at time $t_n = nh$.

\begin{theorem}\label{T2}
Nesterov's accelerated gradient method (\ref{AGM2}) is the numerical discretization of Lie-Trotter splitting applied to two non-autonomous subproblems and one autonomous subproblem at the discretization time $t_n = nh$. Furthermore, for each continuous subproblem, an explicit forward Euler method is applied to each of them.
\end{theorem}

\begin{Proof}
\textbf{(I.)} The first step of the proof is inspired by the work of A.M. Oberman and M. Laborde \cite{Oberman}. We will present an alternative form of the Nesterov's algorithm (\ref{AGM2}) with general vanishing damping with coefficient $\alpha > 1$, that, we will observe, it is very similar to the continuous counterpart (\ref{FirstOrder_VD}). For this, let $\alpha > 1$ and consider the Nesterov's discretization scheme (\ref{AGM2}). At the same time, consider the following alternate form of the Nesterov's accelerated gradient method :
\begin{align}\label{NaG_Velocity}\tag{NaG}
\begin{cases}
x_{n+1} = x_{n} + \dfrac{(\alpha-1)h}{t_n} (v_n-x_n) - \tilde{\beta} \nabla f(y_n) \\
v_{n+1} = v_n - \dfrac{h t_n}{\alpha - 1} \nabla f(y_n)
\end{cases}
\end{align}
Moreover, let $\tilde{\beta} = h^2$ and $y_n = x_n + \dfrac{h(\alpha-1)}{t_n} (v_n-x_n)$. Then, we show that the algorithms (\ref{AGM2}) and (\ref{NaG_Velocity}) are equivalent. Furthermore, $t_n = h(n+\alpha)$ is equivalent to $\alpha_n = \dfrac{n}{n+\alpha}$ and $t_n = hn$ is equivalent to the choice $\alpha_n = \dfrac{n-\alpha}{n}$.
We begin by rewritting the discretization scheme (\ref{NaG_Velocity}) as follows
\begin{align*}
\begin{cases}
x_{n+1} - x_{n} - \dfrac{(\alpha-1)h}{t_n} v_n + \dfrac{(\alpha-1)h}{t_n} x_n = - h^2 \nabla f(y_n) \\ \\
v_{n+1} - v_n = - \dfrac{h t_n}{\alpha - 1} \nabla f(y_n)
\end{cases}
\end{align*}
Multiplying the second equation by $\dfrac{-h(\alpha-1)}{t_n}$ and adding it to the first one, we obtain that :
$$ v_{n+1} = x_n + \dfrac{t_n}{h(\alpha-1)} (x_{n+1}-x_n) \ . $$
Now, using the fact that $y_n = x_n + \dfrac{h(\alpha-1)}{t_n} v_n - \dfrac{h(\alpha-1)}{t_n} x_n$ and that $ v_{n} = x_{n-1} + \dfrac{t_{n-1}}{h(\alpha-1)} (x_{n}-x_{n-1})$, it follows that
$$ y_n = x_n + (x_n - x_{n-1}) \left( \dfrac{t_{n-1}}{t_n} - \dfrac{h(\alpha-1)}{t_n} \right) \ .$$
In this way, we obtain Nesterov's algorithm (\ref{AGM2}), with the inertial coefficient $\alpha_n = \dfrac{1}{t_n} \left( t_{n-1} - h(\alpha-1) \right)$. Finally, if we take the discretization time to be $t_n = nh$, where $h$ is the stepsize, then $\alpha_n = \dfrac{n-\alpha}{n}$. On the other hand, the choice $t_n = h(n+\alpha)$ leads to $\alpha_n = \dfrac{n}{n+\alpha}$. These two represent, in general, the classical choices for the inertial sequence $(\alpha_n)_{n \in \mathbb{N}}$. \\ 
\textbf{(II.)} Now, the second step is to consider splitting the first order dynamical system (\ref{FirstOrder_VD}) in the following subproblems:
\begin{align*}
(P_1):
\begin{cases}
\dot{x}(t) = \dfrac{\alpha-1}{t} (v(t)-x(t)) \\
\dot{v}(t) = - \dfrac{t}{\alpha-1} \nabla f(x(t))
\end{cases}
\quad (P_2) :
\begin{cases}
\dot{x}(t) = - \beta \nabla f(x(t)) \\
\dot{v(t)} = 0
\end{cases}
\end{align*}
Using the notation from the first section, the splitting discrete solution is $\psi(t_n) = \tilde{\phi}^{[2]}(t_n) \circ \tilde{\phi}^{[1]}(t_n)$. We consider $(\tilde{x}_{n+1}, \tilde{v}_{n+1})^{T}$ to be the numerical solution for the first subproblem $(P_1)$ at time $t_n = nh$. Then, applying a forward Euler discretization on $(P_2)$ and taking $\beta = h$, one obtains that
\begin{align*}
\begin{cases}
x_{n+1} = \tilde{x}_{n+1} - h^2 \nabla f(\tilde{x}_{n+1}) \\
v_{n+1} = \tilde{v}_{n+1}
\end{cases}
\end{align*}
where $(x_{n+1}, v_{n+1})^T$ is the current splitting solution at time $t_n = nh$.
Now, also using forward Euler discretization we find the numerical solution $(\tilde{x}_{n+1}, v_{n+1})^{T}$ at time $t_n=nh$ for the continuous subproblem $(P_1)$. As before, we employ splitting the subproblem $(P_1)$ into two non-autonomous subproblems, the first one containing constant velocity $v=v(t_n)$ and the second one preserving the value of $x = x(t_n)$. So, let
\begin{align*}
(P_1)_1 :
\begin{cases}
\dot{x}(t) = \dfrac{\alpha-1}{t} (v(t)-x(t)) \\
\dot{v}(t) = 0
\end{cases}
\quad (P_1)_2 :
\begin{cases}
\dot{x}(t) = 0 \\
\dot{v}(t) = \dfrac{-t}{\alpha-1} \nabla f(x(t))
\end{cases}
\end{align*}
Then, the numerical Lie-Trotter splitting on $(P_1)$ at the discrete time $t_n = nh$ is defined as $\tilde{\phi}^{[1]}(t_n) = (\tilde{\phi}^{[12]} \circ \tilde{\phi}^{[11]})(t_n)$, where $\tilde{\phi}^{[1i]}$  is the forward Euler operator that approximates the flow of the subproblem $(P_1)_i$. If we apply the forward Euler numerical scheme on both the new subproblems, we obtain the following :
\begin{align*}
\begin{cases}
\tilde{x}_{n+1/2} = x_n + \dfrac{h(\alpha-1)}{t_n} (v_n-x_n) \\
\tilde{v}_{n+1/2} = v_n
\end{cases}
\text{ and }
\begin{cases}
\tilde{x}_{n+1} = \tilde{x}_{n+1/2} \\
\tilde{v}_{n+1} = \tilde{v}_{n+1/2} - \dfrac{h t_n}{\alpha-1} \nabla f(\tilde{x}_{n+1/2})
\end{cases}
\end{align*}
Then, it follows that the Lie-Trotter approximation to $(P_1)$ can be written as follows
\begin{align*}
\begin{cases}
\tilde{x}_{n+1} = x_n + \dfrac{h(\alpha-1)}{t_n} (v_n-x_n) \\
\tilde{v}_{n+1} = v_n - \dfrac{h t_n}{\alpha-1} \nabla f \left(  x_n + \dfrac{h(\alpha-1)}{t_n} \left( v_n - x_n \right) \right)
\end{cases}
\end{align*}
where $(x_n,v_n)^{T}$ is the approximation of the solution of the continuous problem at time $t_{n-1} = (n-1)h$ and since we are employing a sequential-type splitting, it is the numerical solution at time $t_{n-1} = (n-1)h$ of $(P_2)$. In this way, putting all together, we obtain that
\begin{align*}
\begin{cases}
x_{n+1} = x_{n} + \dfrac{h(\alpha-1)}{t_n} (v_n - x_n) - h^2 \nabla f \left(  x_n + \dfrac{h(\alpha-1)}{t_n} \left( v_n - x_n \right) \right) \\
v_{n+1} = v_n - \dfrac{ht_n}{\alpha-1} \nabla f \left(  x_n + \dfrac{h(\alpha-1)}{t_n} \left( v_n - x_n \right) \right)
\end{cases}
\end{align*}
By taking $y_n = x_n + \dfrac{h(\alpha-1)}{t_n} (v_n-x_n)$ and using the first step of the proof, we observe that the conclusion follows easily.
\end{Proof}

Finally, we give some observations regarding the Lie-Trotter splitting that was employed in the proof presented before.
\begin{remark}\label{R3}
The numerical solution at time $t_n = nh$ can be written as $\psi(t_n) =  \left( \tilde{\phi}^{[2]} \circ \left( \tilde{\phi}^{[12]} \circ \tilde{\phi}^{[11]} \right) \right) (t_n)$, where the three operators are the approximation to the flows of the three subproblems considered above. Using the ideas that we have presented, we observe that all of the approximations are made at the discrete time $t_n = nh$. Consequently, in contrast to the Strang-Marchuk splitting, we have used the stepsize $h$ and not $h/2$ in the numerical solutions of the subproblems. Finally, we observe that the splitting technique has the advantage that, at a fixed discrete time $t_n$, the same stepsize and the same value of $t_n$ are used in all of the discretized subproblems. \\
Finally, we point out that, in the present proof, we have considered the following technique: we have taken the dynamical system (\ref{FirstOrder_VD}) and then we applied the Lie-Trotter splitting with the choice $\beta = h$. Then, we have showed the equivalence with \ref{NaG_Velocity}. But, we remind that in \ref{NaG_Velocity} we have chosen $\tilde{\beta} = h^2$ and this is valid since this choice does not correspond to a dynamical system, by the fact that it is a parameter from the discretization scheme.
\end{remark}

\subsection{The role of symplectic methods in optimization algorithms}\label{S22}

In this subsection we will make a brief review concerning the role of the sympletic integrators in the setting of unconstrained optimization problems. Furthermore, we point out that many classical inertial algorithms that are based upon dynamical systems with vanishing or with constant damping can be interpreted as a combination between the Lie-Trotter splitting and choices of symplectic methods. Firstly, let's consider the heavy ball with friction model (\ref{HBF}). Also, consider its natural discretization, namely the Polyak's algorithm given by (\ref{PIM}). In \cite{Franca}, it was pointed out that the momentum method can be interpreted as a splitting discretization between conservative and dissipative vector fields. Indeed, consider the subproblems associated to (\ref{HBF}) as follows
\begin{align*}
(P_1) :
\begin{cases}
\dot{x}(t) = 0 \\
\dot{v}(t) = - \gamma v
\end{cases}
\text{ and } (P_2) :
\begin{cases}
\dot{x}(t) = v \\
\dot{v}(t) = - \nabla f(x(t))
\end{cases}
\end{align*}
where the Hamiltonian $H(x,v) = \dfrac{\| v \|^2}{2} + f(x)$ is the energy of the entire dynamical system. The first subproblem $(P_1)$ represents the dissipative system and the second one, i.e. $(P_2)$ is the conservative system. Then, if one applies the explicit forward Euler method on $(P_1)$ and the symplectic Euler method (\ref{SE_2}) on $(P_2)$, it follows that the Lie-Trotter discretization scheme at time $t_n=nh$ , i.e. $\psi(t_n) = \tilde{\phi}^{[2]}(t_n) \circ \tilde{\phi}^{[1]}(t_n)$ becomes 
\begin{align*}
\begin{cases}
y_n = x_n + h (1-h \gamma) v_n \\
x_{n+1} = y_n - h^2 \nabla f(x_n)
\end{cases}
\end{align*}
which is an alternative form of writting the momentum method (\ref{PIM}), using the idea that, from the symplectic discretization of the second subproblem, $v_{n+1} = \dfrac{x_{n+1} - x_n}{h}$ is the discretization of the velocity of the underlying system. \\
Now, we turn our attention to the dynamical system that models Nesterov's accelerated gradient method (\ref{AGM2}). This dynamical system is the true continuous counterpart of Nesterov's algorithm, since under the gradient it contains a continuous version of the inertial term $y_n$, i.e. this inertial term is a discretization of a combination between the exact flow and its first order derivative. For brevity, we recall the idea behind \cite{Jordan}, where the authors considered Nesterov's method as a combination of a Lie-Trotter splitting and a symplectic Euler method. For this, we consider a more generalized vanishing damping system \ref{Extended-AVD}. \\
It is obvious that taking $\gamma = h$, one obtains the Nesterov's algorithm (\ref{AGM2}). Furthermore, (\ref{Extended-AVD}) can be written as a first order evolution equation, using the idea of velocity, as we have pointed out in the case of the heavy ball system : 
\begin{align}
\begin{cases}
\dot{x}(t) = v(t) \\
\dot{v}(t) = -\dfrac{\alpha}{t} v(t) - \nabla f \left( x \left( t \right) + \left( \gamma - \dfrac{\alpha \gamma^2}{t} \right) v(t) \right) 
\end{cases}
\end{align}
By the fact that the above system has no inherent energetic structure, in \cite{Jordan}, the term $\nabla f(x(t))$ was added and substracted, such that one can use, as in the case of the heavy ball with friction system, the dissipative and conservative parts of the system. Thus, the following subproblems were considered :
\begin{align*}
(P_1) :
\begin{cases}
\dot{x}(t) = 0 \\
\dot{v}(t) = - \dfrac{\alpha}{t} v + \nabla f(x(t))  - \nabla f \left( x \left( t \right) + \left( \gamma - \dfrac{\alpha \gamma^2}{t} \right) v(t) \right) 
\end{cases}
\text{ and } (P_2) :
\begin{cases}
\dot{x}(t) = v \\
\dot{v}(t) = - \nabla f(x(t))
\end{cases}
\end{align*}
Here, the first subproblem $(P_1)$ represents the non-potential forces and the second subproblem $(P_2)$ stands for the potential forces. Moreover, as in the case of the heavy ball with friction system, applying a forward Euler method to the first subproblem and the symplectic Euler integrator (\ref{SE_2}) to the second one leads to the Nesterov's method (\ref{AGM2}). \\
Now, we shift our focus to the Hessian-driven damping dynamical system $\left( \text{DIN-AVD}_{\tfrac{\alpha}{t},\beta,1 + \tfrac{\beta}{t}} \right)$. First of all, we will present a natural discretization of this dynamical system that resembles a Polyak-type counterpart to the IGAHD algorithm from \cite{Attouch_Chbani}. Second of all, we show that a discretization resembling the IGAHD algorithm can be obtained, as in the case of Nesterov's method, by a combination of a Lie-Trotter splitting method with a symplectic Euler integrator, but for a Hessian-based system endowed with a perturbation under the gradient as in (\ref{Extended-AVD}). We consider the $\left( \text{DIN-AVD}_{\tfrac{\alpha}{t},\beta,1 + \tfrac{\beta}{t}} \right)$. The natural discretization this system is a Polyak-type method, where at time $t_n = nh$ under the gradient appears the old iteration $x_n$. We approximate $\nabla^2 f(x(t)) \dot{x}(t)$ by the difference of consecutive gradients as in the work of Attouch et. al. \cite{Attouch_Chbani} and we consider the discretization of the continuous system, as follows :
\begin{align*}
& \dfrac{x_{n+1} - 2x_n + x_{n-1}}{h^2} + \dfrac{\alpha}{t_n} \dfrac{x_n - x_{n-1}}{h} + \dfrac{\beta}{h} \left[ \nabla f(x_n) - \nabla f(x_{n-1}) \right] = - \nabla f(x_n) - \dfrac{\beta}{t_n} \nabla f(x_n) \\
& \dfrac{x_{n+1} - 2x_n + x_{n-1}}{h^2} + \dfrac{\alpha}{nh} \dfrac{x_n - x_{n-1}}{h} + \dfrac{\beta}{h} \left[ \nabla f(x_n) - \nabla f(x_{n-1}) \right] = - \nabla f(x_n) - \dfrac{\beta}{nh} \nabla f(x_n) 
\end{align*}
By some easy computations, it follows that
$$ x_{n+1} = x_n + \left( 1 - \dfrac{\alpha}{n} \right) (x_n - x_{n-1}) - \beta h \left[ \nabla f(x_n) - \nabla f(x_{n-1}) \right] - \dfrac{\beta h}{n} \nabla f(x_n) - h^2 \nabla f(x_n) $$
Define $\alpha := 1 - \dfrac{\alpha}{n}$. Then, we obtain that
\begin{align}\label{Polyak-IGAHD}\tag{Polyak-IGAHD}
\begin{cases}
& y_n = x_n + \alpha_n (x_n - x_{n-1}) - \beta h \left[ \nabla f(x_n) - \nabla f(x_{n-1}) \right] - \dfrac{\beta h}{n} \nabla f(x_n) \\
& x_{n+1} = y_n - h^2 \nabla f(x_n) 
\end{cases}
\end{align}
Now, our only result of this subsection consists in showing that a similar algorithm to the IGAHD (see \cite{Attouch_Chbani}) with a predictor $y_n$ under the gradient can be obtained as in the case of (\ref{Extended-AVD}) by applying a sequential splitting method with a symplectic Euler integrator on the second subproblem. Furthermore, the only difference between our algorithm and IGAHD is that the latter one contains $\nabla f(x_{n-1})$ instead $\nabla f(x_n)$. The choice $\nabla f(x_{n-1})$ seems unnatural, by the fact that at the discrete time $t_n = nh$, the gradient of the flow is approximated by $\nabla f(x_{n-1})$ and not with the current gradient value. Nevertheless, by using the same technique of approximation, we can easily obtain a similar discretization scheme with the IGAHD algorithm.
\begin{theorem}\label{T4}
Let $f \in C^2 (\mathbb{R}^N, \mathbb{R})$ and consider the following Hessian driven damping second order system :
\begin{align*}
\hspace*{-0.5cm}
\ddot{x}(t) + \dfrac{\alpha}{t} \dot{x}(t) + \beta \nabla^2 f (x(t)) \dot{x}(t) = - \dfrac{\beta}{t} \nabla f(x(t)) - \nabla f \left( x(t) + \gamma \left( 1-\dfrac{\alpha \gamma}{t} \right) \dot{x}(t) - \beta \gamma^2 \nabla^2 f(x(t)) \dot{x}(t) - \dfrac{\beta \gamma^2}{t} \nabla f(x(t)) \right)
\end{align*}
Then, a Lie-Trotter splitting discretization, with a forward Euler method on the first subproblem and a symplectic integrator on the second one leads to an IGAHD-type algorithm.
\end{theorem}

\begin{Proof}
As in the case of (\ref{Extended-AVD}) we consider two subproblems, one containing the so-called non-potential forces and the other one containing the potential forces :
\begin{align*}
(P_1) :
\begin{cases}
\dot{x} = 0 \\
\dot{v} = - \dfrac{\alpha}{t} v(t) - \beta \nabla^2 f (x(t)) v(t) - \dfrac{\beta}{t} \nabla f(x(t)) - \nabla f \left( \mathcal{A}(t, x(t), v(t)) \right) + \nabla f(x(t))
\end{cases}
\end{align*}
where 
$$ \mathcal{A}(t, x(t), v(t)) = x(t) + \gamma \left( 1-\dfrac{\alpha \gamma}{t} \right) v(t) - \beta \gamma^2 \nabla f(x(t)) v(t) - \dfrac{\beta \gamma^2}{t} \nabla f(x(t)) $$
and the second subproblem as
\begin{align*}
(P_2) :
\begin{cases}
\dot{x} = v \\
\dot{v} = - \nabla f(x(t))
\end{cases}
\end{align*}
Taking $\gamma = h$ and defining $\alpha_n := 1 - \dfrac{\alpha h}{t_n}$ we obtain that
\begin{align*}
\begin{cases}
x_{n+1/2} = x_n \\
v_{n+1/2} = \alpha_n v_n - \beta h \nabla^2 f(x_n) v_n - \dfrac{\beta h}{t_n} \nabla f(x_n) + h \nabla f(x_n) - h \nabla f \left( x_n + h \alpha_n v_n - \beta h^2 \nabla^2 f(x_n) v_n - \dfrac{\beta h^2}{t_n} \nabla f(x_n) \right)
\end{cases}
\end{align*}
Furthermore, for the second subproblem $(P_2)$, we apply the symplectic Euler integrator (\ref{SE_2}) and obtain that
\begin{align*}
\begin{cases}
x_{n+1} = x_{n+1/2} + h v_{n+1} \\
v_{n+1} = v_{n+1/2} - h \nabla f(x_{n+1/2})
\end{cases}
\end{align*}
Using the fact that for $t_{n} = nh$, the splitting solution is the discrete solution of $(P_2)$, namely $(x_{n+1},v_{n+1})^T \in \mathbb{R}^N \times \mathbb{R}^N$, we obtain from the first equation given above, that $v_n = \dfrac{x_n-x_{n-1}}{h}$. By using this and the approximation involving the Hessian, i.e. $\nabla^2 f(x_n) v_n$ at time $t_n = nh$ by two consecutive gradients, it follows that
\begin{align}\label{IGAHD}\tag{IGAHD-type}
\begin{cases}
& y_n = x_n + \alpha_n (x_n - x_{n-1}) - \beta h \left[ \nabla f(x_n) - \nabla f(x_{n-1}) \right] - \dfrac{\beta h}{n} \nabla f(x_n) \\
& x_{n+1} = y_n - h^2 \nabla f(y_n)
\end{cases}
\end{align}
\end{Proof}

We end this section by emphasizing the role of choosing the appropriate first order evolution equation for the Hessian driven damping system in the splitting methods.

\begin{remark}\label{R5}
1.) At the second order dynamical system from Theorem \ref{T4}, under the gradient at $\mathcal{A}(t,x(t),\dot{x}(t))$, we can approximate $\dfrac{\beta \gamma^2}{t} \nabla f(x(t))$ as $\dfrac{\beta h^2}{t_n} \nabla f(x_{n-1})$. In this way, we obtain the classical IGAHD algorithm from \cite{Attouch_Chbani}. \\
2.) The idea behind Theorem \ref{T4} is the fact that we separate the perturbed gradient $-\dfrac{\beta}{t} \nabla f(x(t))$ from the unperturbed gradient value. Moreover, the latter one is similar to the (\ref{Extended-AVD}) since, under the gradient, the value contains a combination of both the exact flow $x=x(t)$ and its first order derivative $\dot{x}=\dot{x}(t)$.\\
3.) In \cite{ShiDuJordan}, Shi et. al. have studied a symplectic Euler method that was applied to a full Hessian driven damping system similar to \ref{DIN-AVD}. The key difference between our analysis and theirs is that our algorithms are also based upon a Lie-Trotter discretization in combination with symplectic methods.
\end{remark}

\subsection{Discussion on some splitting techniques for the extended AVD dynamical system}\label{S23}

Until now, we have used the symplectic Euler method (\ref{SE_2}) but, in the present subsection, we consider a Lie-Trotter splitting discretization scheme by combining a forward Euler method on the first subproblem and the symplectic Euler integrator (\ref{SE_1}) on the second one, on the dynamical system (\ref{Extended-AVD}). In this way we consider the following two subproblems:
\begin{align*}
(P_1) :
\begin{cases}
\dot{x}(t) = 0 \\
\dot{v}(t) = - \dfrac{\alpha}{t} v + \nabla f(x(t))  - \nabla f \left( x \left( t \right) + \left( \gamma - \dfrac{\alpha \gamma^2}{t} \right) v(t) \right) 
\end{cases}
\text{ and } (P_2) :
\begin{cases}
\dot{x}(t) = v \\
\dot{v}(t) = - \nabla f(x(t))
\end{cases}
\end{align*}
By using a forward Euler discretization on $(P_1)$ and setting $\gamma = h$, one obtains that
\begin{align*}
\begin{cases}
& x_{n+1/2} = x_n \\
& v_{n+1/2} = \left( 1- \dfrac{\alpha h}{t_n} \right) v_n - h \nabla f \left( x_n + h \left( 1 - \dfrac{\alpha h}{t_n} \right) v_n \right) + h \nabla f(x_n)
\end{cases}
\end{align*}
Also, applying the symplectic Euler method (\ref{SE_1}) on the subproblem $(P_2)$, it follows that
\begin{align*}
\begin{cases}
& x_{n+1} = x_{n+1/2} + h v_{n+1/2} \\
& v_{n+1} = v_{n+1/2} - h \nabla f(x_{n+1})
\end{cases}
\end{align*}
Now, as before, using the discretization time $t_n = nh$ and denoting $\alpha_n := 1 - \dfrac{\alpha}{n}$, we introduce the auxiliary term $y_n = x_n + h \alpha_n v_n$. So, we get that
\begin{align*}
\begin{cases}
& x_{n+1} = x_n + h \alpha_n v_n - h^2 \nabla f(y_n) + h^2 \nabla f(x_n) \\
& v_{n+1} = \alpha_n v_n - h \nabla f(y_n) + h \nabla f(x_n) - h \nabla f(x_{n+1}) 
\end{cases}
\end{align*}
By multiplying the second equation by $h$ and substracting them, we obtain the perturbed velocity which differs by the one obtained in the case of (\ref{SE_2}) by an extra gradient, namely: 
$$ v_n = \dfrac{x_n - x_{n-1}}{h} - h \nabla f(x_n) \, . $$
Finally, we obtain our first algorithm, i.e.
\begin{align}\label{LT-SE1}\tag{LT-SE1}
\begin{cases}
& y_n = x_n + \alpha_n (x_n - x_{n-1}) - h^2 \alpha_n \nabla f(x_n) \\
& x_{n+1} = y_n - h^2 \nabla f(y_n) + h^2 \nabla f(x_n)
\end{cases}
\end{align}
Now, we consider a similar approach, namely on subproblem $(P_1)$ we apply the forward Euler method, and on $(P_2)$ we consider the Stormer-Verlet algorithm (\ref{SV_2}). As before, for the first subproblem we have that
\begin{align*}
\begin{cases}
& x_{n+1/2} = x_n \\
& v_{n+1/2} = \left( 1- \dfrac{\alpha h}{t_n} \right) v_n - h \nabla f \left( x_n + h \left( 1 - \dfrac{\alpha h}{t_n} \right) v_n \right) + h \nabla f(x_n)
\end{cases}
\end{align*}
Furthermore, using (\ref{SV_2}) on the second sub-problem, it follows that
\begin{align*}
\begin{cases}
& \tilde{v}_{n+1/2} = v_{n+1/2} - \dfrac{h}{2} \nabla f(x_{n+1/2}) \\
& x_{n+1} = x_{n+1/2} + h \tilde{v}_{n+1/2} \\
& v_{n+1} = \tilde{v}_{n+1/2} - \dfrac{h}{2} \nabla f(x_{n+1})
\end{cases}
\end{align*}
As before, using $y_n = x_n + h \alpha_n v_n$, where $t_n$ is the discretization time and $\alpha_n$ represents the intertial momentum coefficient, it follows that
\begin{align*}
\begin{cases}
& x_{n+1} = x_n + h \alpha_n v_n - h^2 \nabla f(y_n) + h^2 \nabla f(x_n) - \dfrac{h^2}{2} \nabla f(x_n) \\
& v_{n+1} = \alpha_n v_n - h \nabla f(y_n) + h \nabla f(x_n) - \dfrac{h}{2} \nabla f(x_n) - \dfrac{h}{2} \nabla f(x_{n+1}) 
\end{cases}
\end{align*}
In a similar spirit with the determination of the algorithm (\ref{LT-SE1}), we obtain our second algorithm, i.e.
\begin{align}\label{LT-SV2}\tag{LT-SV2}
\begin{cases}
& y_n = x_n + \alpha_n (x_n - x_{n-1}) - \dfrac{h^2}{2} \alpha_n \nabla f(x_n) \\
& x_{n+1} = y_n - h^2 \nabla f(y_n) + \dfrac{h^2}{2} \nabla f(x_n)
\end{cases}
\end{align}
We emphasize that our newly introduced algorithms (\ref{LT-SE1}) and (\ref{LT-SV2}) differ by a $\dfrac{h}{2}$ term at the gradient $\nabla f(x_n)$, in the case of the Stormer-Verlet integrator the perturbed velocity contains also a $\dfrac{h}{2}$ term, which in an intuitive sense, it annihilates much more properly large gradient values. \\
Now, we shift our attention to an algorithm that was developed by N.C. Nguyen et.al. in \cite{Nguyen}, that was used as a residual method for solving systems of equations. Using a variant related to the constant stepsize $h$, the residual method is the following:
\begin{align}\label{ARDM}\tag{ARDM}
\begin{cases}
& y_n = x_n + \alpha_n (x_n - x_{n-1}) - h^2 ( 1 + \alpha_n ) \nabla f(x_n) \\
& x_{n+1} = y_n - h^2 \nabla f(y_n)
\end{cases}
\end{align}
where the method reads as follows: \textit{accelerated residual descent method}. Also, N.C. Nguyen et.al. have considered a restarting strategy that it is out of the scope of the present paper. Furthermore, their arguments into developing (\ref{ARDM}) are based upon the stability theory related to finite difference approximations and are also based upon an approach similar to that of Su et.al. (see \cite{Boyd}). Our next purpose is to introduce a new dynamical system, such that the optimization algorithm (\ref{ARDM}) can be considered as a Lie-Trotter splitting for a combination between the forward Euler and the symplectic Euler method (\ref{SE_2}), as in the case of (\ref{Extended-AVD}) dynamical system. For this, we consider the following second order evolution equation:
\begin{align}\label{DynSys-ARDM}\tag{DynSys-ARDM}
\ddot{x}(t) + \dfrac{\alpha}{t} \dot{x}(t) = - \nabla f \left( x \left( t \right) + \left( \gamma - \dfrac{\alpha \gamma^2}{t} \right) \dot{x}(t) - \gamma^2 \left( 2 - \dfrac{\alpha \gamma}{t} \right) \nabla f(x(t)) \right) - \left( 2 - \dfrac{\alpha \gamma}{t} \right) \nabla f(x(t))
\end{align}
As before, we divide this evolution equation into dissipative and conservative parts, i.e. into non-potential and potential forces:
\begin{align*}
\hspace*{-0.5cm}
(P_1) :
\begin{cases}
\dot{x}(t) = 0 \\
\dot{v}(t) = - \dfrac{\alpha}{t} v(t) + \nabla f(x(t))  - \nabla f \left( x \left( t \right) + \left( \gamma - \dfrac{\alpha \gamma^2}{t} \right) v(t) - \gamma^2 \left( 2 - \dfrac{\alpha \gamma}{t} \right) \nabla f(x(t)) \right) -  \left( 2 - \dfrac{\alpha \gamma}{t} \right) \nabla f(x(t))
\end{cases}
\end{align*}
and
\begin{align*}
(P_2) :
\begin{cases}
\dot{x}(t) = v \\
\dot{v}(t) = - \nabla f(x(t))
\end{cases}
\end{align*}
Using a forward Euler explict method on $(P_1)$ and taking $\gamma = h$, we obtain that
\begin{align*}
\begin{cases}
& x_{n+1/2} = x_n \\
& v_{n+1/2} = \left( 1- \dfrac{\alpha h}{t_n} \right) v_n - \left( 2 - \dfrac{\alpha h}{t_n} \right) \nabla f(x_n) - h \nabla f \left( x_n + h \left( 1 - \dfrac{\alpha h}{t_n} \right) v_n - h^2 \left( 2 - \dfrac{\alpha h}{t_n} \right) \nabla f(x_n) \right) + h \nabla f(x_n)
\end{cases}
\end{align*}
On the other hand, using the symplectic Euler method (\ref{SE_2}) on $(P_2)$, we arrive at the following equalities:
\begin{align*}
\begin{cases}
& x_{n+1} = x_{n+1/2} + h v_{n+1} \\
& v_{n+1} = v_{n+1/2} - h \nabla f(x_{n+1/2}) 
\end{cases}
\end{align*}
Now, in contrast with the previous algorithm, at the discrete time $t_n = nh$, we consider the additional iteration $y_n = x_n + h \alpha_n v_n - h^2 (1+\alpha_n) \nabla f(x_n)$. Using this notation, we get the following:
\begin{align*}
\begin{cases}
& x_{n+1} = x_n + h \alpha_n v_n - h^2 \nabla f(y_n) - h^2 (1+\alpha_n) \nabla f(x_n) \\
& v_{n+1} = \alpha_n v_n - h \nabla f(y_n) -
h (1+\alpha_n) \nabla f(x_n)
\end{cases}
\end{align*}
As in our splitting-based algorithms that were presented above, we multiply the second equation by $h$ and by substracting them, we get the algorithm (\ref{ARDM}). Finally, we end this subsection by considered a remark that will represent a key role in our theoretical analysis.

\begin{remark}\label{R6}
So far, we have considered three algorithms, namely (\ref{ARDM}), (\ref{LT-SE1}) and (\ref{LT-SV2}). The difference between the first one and the symplectic methods applied to the (\ref{Extended-AVD}) dynamical system is that the accelerated residual method contains a perturbation $h^2 (1+\alpha_n)$ only in the additional inertial term $y_n$. The sequential Euler and Stormer-Verlet algorithms for the extended dynamical system with asymptotically vanishing damping contains a perturbation at $y_n$ and also at the current iteration value $x_{n+1}$. In our terms, if an element appears at $x_{n+1}$ and it contains the gradient $\nabla f(x_n)$, we say that the algorithm has a \textit{bias correction term}. This is the case for the (\ref{LT-SE1}) and (\ref{LT-SV2}) algorithms. On the other hand, we can say that (\ref{ARDM}) is \textit{bias corrected free}. This difference can be seen in the definition of the inertial momentum term, where at (\ref{ARDM}), $y_n$ contains $h^2 (1+\alpha_n) \nabla f(x_n)$. So, now we consider a general form of inertial type algorithms that contain two gradient evaluations and that encompasses (\ref{ARDM}), (\ref{LT-SE1}) and (\ref{LT-SV2}):
\begin{align}\label{LT-Symplectic}\tag{LT-Symplectic}
\begin{cases}
& y_n = x_n + \alpha_n (x_n - x_{n-1}) - \omega_n \nabla f(x_n) \\
& x_{n+1} = y_n - h^2 \nabla f(y_n) + \gamma_n \nabla f(x_n)
\end{cases}
\end{align}
where $\alpha_n = 1 - \dfrac{\alpha}{n}$ and $t_n = nh$. Also, we point out that $\omega_n$ is the underlying weight that appears in the Lie-Trotter discretization and $\gamma_n$ is the bias correction term. Last, we observe that for (\ref{LT-SE1}) we have that $\omega_n = h^2 \alpha_n$ and $\gamma_n = h^2$. Furthermore, for (\ref{LT-SV2}), these values are halved. On the other hand, (\ref{ARDM}) can be obtained by setting $\omega_n = h^2 (1+\alpha_n)$ and $\gamma_n = 0$, for each $n \in \mathbb{N}$. 
\end{remark}

\subsection{Fixed points, spurious solutions and exploding gradients : Heuristic interpretation}\label{S24}

In this section we present some intuitive explanations regarding the general algorithm of type \ref{LT-Symplectic}. We consider the asymptotic behavior, the critical points and the equilibrium states related to our optimization methods. Also, regarding the phenomenons that we observe, we shall present also some remedies. \\
First of all, we shall start with Nesterov algorithm \ref{AGM2}. Let's suppose that there exists $x^{\ast} \in \mathbb{R}^N$ such that the sequence $(x_n)_{n \in \mathbb{N}}$ is convergent to $x^{\ast}$. Then, let $\lim\limits y_n = y^{\ast}$. We observe that $y^{\ast} = x^{\ast}$ and so $\nabla f(x^{\ast}) = 0$. This is equivalent to the fact that $x^{\ast} \in \argmin f$. At the same time, define the operator $T_h : \mathbb{R}^N \to \mathbb{R}^N$, as $T_h(y) = y - h^2 \nabla f(y)$. From all of this it follows that if $\lim\limits_{n \to \infty} x_n = x^{\ast}$, then we have that $x^{\ast} \in \argmin f$, i.e. $x^{\ast} = T_h(x^{\ast})$. So, the behavior of the Nesterov scheme is related to the fact that if the principal sequence converges, then the limit is also a minimizer of the convex function $f$ and at the same time it is a fixed point of the associated operator $T_h$. \\
Now, we turn our attention to the \ref{IGAHD} algorithm. Keeping in mind our notation from \ref{LT-Symplectic}, let $\omega_n = \dfrac{\beta h}{n}$. In this case $\lim\limits_{n \to \infty} \omega_n = 0$. The fact that the limit of the sequence $(\omega_n)_{n \in \mathbb{N}}$ is $0$ will be crucial in our analysis. We shall use the fact that two consecutive iterates of the gradient represent a natural approximation of the Hessian, and this is related to $\lim\limits_{n \to \infty} \left[ \nabla f(x_n) - \nabla f(x_{n-1}) \right] = 0$, since $\lim\limits_{n \to \mathbb{N}} x_n = x^{\ast}$. As before, we have considered that $x^{\ast}$ the limit of $(x_n)_{n \in \mathbb{N}}$. Again, denoting by $y^{\ast}$ the limit of $(y_n)_{n \in \mathbb{N}}$, we find that, as in the case of \ref{AGM2}, $x^{\ast} \in \argmin f$ (here an important property relies on the fact that the coefficient of the difference between consecutive gradient $\beta h$ is bounded). So, the same conclusion is valid also for the \ref{IGAHD} algorithm, regarding the fixed points of the operator $T_h$ and the minimizers of $f$. So, if the major iteration of the algorithm converges, then the limit is also a minimizer of the convex objective function $f$. \\
We have said that the property that $\lim\limits_{n \to \infty} \omega_n = 0$ is crucial. We will call this phenomenon by the name \emph{annihilating the exploding gradient}, since if the limit was strictly positive then the gradient at $x_n$ would not be disappearing at all, and so we are left with a residual that, as time grows, will tend to the gradient of $x^{\ast}$, which could be large enough if $x^{\ast}$ was not a minimizer. For instance, let $\omega = \lim\limits_{n \in \mathbb{N}} \omega_n$, such that $\omega > 0$. As before, let $\lim\limits_{n \in \mathbb{N}} x_n = x^{\ast}$ and $\lim\limits_{n \in \mathbb{N}} y_n = y^{\ast}$. It is trivial to show that $y^{\ast} = x^{\ast} - \omega \nabla f(x^{\ast})$, and so $x^{\ast}$ verifies the following equation:
\begin{align}\label{eq1}
\nabla f(x^{\ast}) = - \dfrac{h^2}{\omega} \nabla f(x^{\ast} - \omega \nabla f(x^{\ast}))
\end{align}
It is easy to see that a minimizer satisfies \ref{eq1}. On the other hand, if $x^{\ast}$ is not a minimizer then we have the so called \emph{spurious roots}. This means that not all the solutions of \ref{eq1} are minimizers of the objective function $f$. As in the previous cases, by the fact that $x_{n+1} = y_n - h^2 \nabla f(y_n)$, let $T_h$ defined as before. A simple analysis shows that the limit $x^{\ast}$ of $(x_n)_{n \in \mathbb{N}}$ satisfies the relation \ref{eq1}. Also, we have that $y^{\ast} = T_h(y^{\ast})$ if and only if $y^{\ast} = x^{\ast} - \omega \nabla f(x^{\ast}) \in \argmin f$. This is the reason we also call these elements \emph{spurious fixed points}, by the fact that they are fixed points of $T_h$, but they can be represented as a perturbation of the limit point $x^{\ast}$. \\
Now, let's consider the \ref{ARDM} algorithm. Using, the same notation as before, we have that $\omega = 2h^2$, which is strictly positive. By making the same analysis we reach the following equation between $x^{\ast}$ and $y^{\ast}$, i.e
\begin{align}\label{eq2}
\nabla f(x^{\ast}) = - \dfrac{1}{2} \nabla f(x^{\ast} - 2h^2 \nabla f(x^{\ast}))
\end{align}
This comes from \ref{eq1} by taking $\omega = 2 h^2$. Also, as for the case of the previous discussions regarding the modification of $\omega$ in \ref{IGAHD}, we obtain the same conclusion regarding the spurious roots of \ref{ARDM}. \\
The algorithm \ref{LT-SE1} meets the requirement that $\omega = h^2 > 0$. Also, letting $\lim\limits_{n \to \infty} x^{\ast}$ and $\lim\limits_{n \to \infty y_n} = y^{\ast}$, one can show that we obtain the equation
\begin{align}\label{eq3}
\nabla f(x^{\ast}) = - \nabla f(x^{\ast} - h^2 \nabla f(x^{\ast}))
\end{align} 
So, in all of the three cases from above, namely \ref{eq1}, \ref{eq2} and \ref{eq3} we have that the limit point $x^{\ast}$ is not necessarily a minimizer of the convex function $f$. We want also to point out that, for \ref{LT-SE1}, we have that $x_{n+1} = y_n - h^2 \nabla f(y_n) + h^2 \nabla f(x_n)$.In this case, let $T_h : \mathbb{R}^N \times \mathbb{R}^N \to \mathbb{R}^N$ be defined as $T_h(y,x) = y - h^2 \nabla f(y) + h^2 \nabla f(x)$. In this case, by simple computations, we obtain that, for the case when $x^{\ast}$ is the limit of $(x_n)_{n \in \mathbb{N}}$ and $y^\ast$ is the limit of $(y_n)_{n \in \mathbb{N}}$, we have that $T_h(y^{\ast},x^{\ast}) = T_h(x^{\ast} - h^2 \nabla f(x^{\ast}),x^{\ast}) = x^{\ast} - h^2 \nabla f(y^{\ast})$. So, we have that $x^{\ast} = T_h(y^{\ast},x^{\ast})$ is equivalent with the fact that $y^{\ast} = x^{\ast} - h^2 \nabla f(x^{\ast}) \in \argmin f$. Nevertheless, let us notice that $x^{\ast} = T_h(x^{\ast},x^{\ast})$ means that $x^{\ast}$ is a minimizer. So, it is worth pointing out that the limit point $x^\ast$ minimizes $f$ if and only if is a fixed point of $T_h$ on the cartesian product $\mathbb{R}^N \times \mathbb{R}^N$, i.e. a fixed point in the classical sense (i.e. on the diagonal set). \\
Bearing all these in mind, we present a remedy for the algorithms that belong to the class \ref{LT-Symplectic}, such that they preserve fixed points, do not have the property of gradient explosion and do not introduce spurious roots. For \ref{LT-Symplectic}, suppose that $\lim\limits_{n \to \infty} \omega_n = 0$, $\lim\limits_{n \to \infty} \gamma_n = \gamma \geq  0$ and that $\gamma \neq h^2$. As before, let $x^{\ast}$ be the limit of $(x_n)_{n \in \mathbb{N}}$ and $y^{\ast}$ be the limit of $(y_n)_{n \in \mathbb{N}}$. In this case, by the fact that $\omega_n$ converges to $0$ such that the gradient $\nabla f(x^{\ast})$ is annihilated, we obtain that $y^{\ast} = x^{\ast}$. This leads to the fact that $(\gamma-h^2) \nabla f(x^{\ast}) = 0$, so $x^{\ast} \in \argmin f$. In this case, $T_h(y,x) = y - h^2 \nabla f(y) + \gamma \nabla f(x)$. In this case, $T_h(y^{\ast},x^{\ast}) = x^{\ast}$ is equivalent with $x^{\ast} \in \argmin f$. So, with the above conditions on the coefficients, we have a natural extension of the asymptotic behavior of Nesterov algorithm \ref{AGM2}, where the operator $T_h(y)$ is extended to the cartesian product by the mean of $T_h(y,x)$. \\
We end this section, by considered another remedy but a little bit more particular, namely we modify explicitly \ref{LT-SE1}, such that it has the properties presented above. In the case of \ref{LT-SE1}, by our symplectic approach, we have the perturbed velocity $v_n = \dfrac{x_n - x_{n-1}}{h} - h \nabla f(x_n)$. This is different from the symplectic approach of Nesterov's algorithm \ref{AGM2} by the fact that, in our proposed scheme, $\lim\limits_{n \to \infty} v_n = - h \nabla f(x^{\ast})$. Regarding the \ref{Extended-AVD} dynamical system, the velocity $v$ satisfies $\lim\limits_{t \to \infty} v(t) = 0$, i.e. the velocity of the continuous system approaches $0$ as time grows. We want to reflect this property also to our discretization method. So, we want to have that $v_n = \dfrac{x_n - x_{n-1}}{h} - h \theta_{n-1} \nabla f(x_n)$, with $\lim\limits_{n \to \infty} \theta_n = 0$. We recall that in the construction of \ref{LT-SE1}, we have used the technique of adding and subtracting $\nabla f(x(t))$ in order to split the potential and non-potential forces. We shall make this trick, but for $\theta(t) \nabla f(x(t))$, where $\theta$ vanishes at $\infty$. Furthermore, we consider the discretized version, i.e. $\theta(t_n) \approx \theta_n$. \\
As for \ref{LT-SE1}, we consider a Lie-Trotter splitting discretization scheme by combining a forward Euler method on the first subproblem and the symplectic Euler integrator (\ref{SE_1}) on the second one, on the dynamical system (\ref{Extended-AVD}). Let the following two subproblems:
\begin{align*}
(P_1) :
\begin{cases}
\dot{x}(t) = 0 \\
\dot{v}(t) = - \dfrac{\alpha}{t} v + \theta(t) \nabla f(x(t))  - \nabla f \left( x \left( t \right) + \left( \gamma - \dfrac{\alpha \gamma^2}{t} \right) v(t) \right) 
\end{cases}
\text{ and } (P_2) :
\begin{cases}
\dot{x}(t) = v \\
\dot{v}(t) = - \theta(t) \nabla f(x(t))
\end{cases}
\end{align*}
By using a forward Euler discretization on $(P_1)$, it follows that
\begin{align*}
\begin{cases}
& x_{n+1/2} = x_n \\
& v_{n+1/2} = \left( 1- \dfrac{\alpha h}{t_n} \right) v_n - h \nabla f \left( x_n + h \left( 1 - \dfrac{\alpha h}{t_n} \right) v_n \right) + h \theta_{n} \nabla f(x_n)
\end{cases}
\end{align*}
Also, applying the symplectic Euler method (\ref{SE_1}) on the subproblem $(P_2)$, it follows that
\begin{align*}
\begin{cases}
& x_{n+1} = x_{n+1/2} + h v_{n+1/2} \\
& v_{n+1} = v_{n+1/2} - h \theta_{n} \nabla f(x_{n+1})
\end{cases}
\end{align*}
Using the discrete time $t_n = nh$ and keeping in mind that $\alpha_n := 1 - \dfrac{\alpha}{n}$ and that $y_n = x_n + h \alpha_n v_n$, we get:
\begin{align*}
\begin{cases}
& x_{n+1} = x_n + h \alpha_n v_n - h^2 \nabla f(y_n) + h^2 \theta_{n} \nabla f(x_n) \\
& v_{n+1} = \alpha_n v_n - h \nabla f(y_n) + h \theta_{n} \nabla f(x_n) - h \theta_{n} \nabla f(x_{n+1}) 
\end{cases}
\end{align*}
By multiplying the second equation by $h$ and substracting them, we obtain the modified velocity, i.e.
$$ v_n = \dfrac{x_n - x_{n-1}}{h} - h \theta_{n-1} \nabla f(x_n) \, . $$
In this way, we obtain the third symplectic-type method:
\begin{align}\label{LT-SE3}\tag{LT-SE3}
\begin{cases}
& y_n = x_n + \alpha_n (x_n - x_{n-1}) - h^2 \alpha_n \theta_{n-1} \nabla f(x_n) \\
& x_{n+1} = y_n - h^2 \nabla f(y_n) + h^2 \theta_{n} \nabla f(x_n)
\end{cases}
\end{align}
In this case, we have that $\lim\limits_{n \to \infty} \omega_n = \lim\limits_{n \to \infty} h^2 \alpha_n \theta_{n-1} = 0$ and also $\lim\limits_{n \to \infty} \gamma_n = \lim\limits_{n \to \infty} h^2 \theta_{n} = 0$. Finally, we observe that we are in the case of our previous discussion regarding the remedy concerning the exploding gradients and the spurious roots.

\section{Proposed algorithms and their rate of convergence}\label{S3}

Even though we have introduced new algorithms, namely \ref{LT-Symplectic} via an intuitive geometric approach by the means of symplectic methods, we will focus our attention on a generalized Hessian driven damping system and not on the extended AVD system. The discussion about our choice is postponed to Remark \ref{R11}. Hence, our aim is to consider a theoretical and an empirical convergence analysis regarding Nesterov's algorithm (\ref{AGM2}), IGAHD algorithm from \cite{Attouch_Chbani} and our approach via symplectic integrators, i.e. algorithms of the form \ref{LT-S-IGAHD}.

\subsection{Modified Hessian-driven damping continuous and discrete systems}\label{S31}

In the previous sections, we have presented the key idea behind the splitting technique applied to gradient-type dynamical systems. We have considered \ref{LT-Symplectic} and \ref{LT-SE3} algorithms that have a geometric interpretation that is related to symplectic integrators. We observe that all these optimization algorithms have two coefficients in front of the gradient value at $x_n$, namely $\gamma_n$ and $\omega_n$, for each $n \in \mathbb{N}$. On the other hand, there are the Hessian-driven damping discretizations like \ref{IGAHD}, which contain the difference of two consecutive gradients, i.e $\nabla f(x_n) - \nabla f(x_{n-1})$ along with $\nabla f(x_n)$, under the gradient of the inertial element $y_n$. By the fact that the effect of the Hessian is to dampen the oscillations, it is more natural to consider a Hessian-driven damping algorithm. So, in the present subsection we shall modify some of the previous dynamical systems in order to make consistent modification into our proposed algorithms. At the same time, our algorithms that we propose shadow the influence of the spurious roots and exploding gradients (explained in the Subsection \ref{S24}), preserve the Hessian structure and are derived from a symplectic perspective as the other methods derived above. In order to do this, for $n \geq 1$ we consider the following generalization of \ref{IGAHD}, namely the \textit{Lie-Trotter-Symplectic-IGAHD}, i.e.
\begin{align}\label{LT-S-IGAHD}\tag{LT-S-IGAHD}
\begin{cases}
& y_n = x_n + \alpha_n (x_n - x_{n-1}) - \lambda_n \left[ \nabla f(x_n) - \nabla f(x_{n-1}) \right] - \omega_n \nabla f(x_n) \\
& x_{n+1} = y_n - h^2 \nabla f(y_n) + \gamma_n \nabla f(x_n) 
\end{cases}
\end{align}
where $\lambda_n$ and the coefficients $\omega_n$ and $\gamma_n$ are consistent with the notations from \ref{LT-Symplectic}. Furthermore, for $\gamma_n = 0$, $\lambda_n = \beta h$ and $\omega_n = \dfrac{\beta h}{n}$, we obtain as a particular case \ref{IGAHD}. On the other hand, taking into account Subsection \ref{S24}, it is naturally to consider $\omega_n$ such as $\lim\limits_{n \to \infty} \omega_n = 0$. Also, for the case of \ref{IGAHD} this property is easily verifiable, with $\omega_n = O \left(  \dfrac{1}{n} \right)$ as $n \to \infty$. Now, we present our result of this subsection, in which we introduce a new dynamical system such that \ref{LT-S-IGAHD} is a Lie-Trotter splitting discretization.

\begin{theorem}\label{T7}
Let $f \in C^2 (\mathbb{R}^N, \mathbb{R})$ and consider the following Hessian driven damping second order system :
\begin{align*}
\hspace*{-0.5cm}
\ddot{x}(t) + \dfrac{\alpha}{t} \dot{x}(t) + \dfrac{\lambda(t)}{\gamma} \nabla f^2 (x(t)) \dot{x}(t) = - \dfrac{1}{\gamma^2} \left( \omega(t) - \tilde{\gamma}(t) \right) \nabla f(x(t)) - \nabla f \left( \mathcal{A}(t, x(t), \dot{x}(t)) \right)
\end{align*}
where $\tilde{\gamma} = \tilde{\gamma}(t)$ and $\omega = \omega(t)$ are given functions. Also, for a pair of the form $(t, x(t), y(t))$ we have used the notation 
$$ \mathcal{A}(t, x(t), y(t)) := x(t) + \gamma \left( 1-\dfrac{\alpha \gamma}{t} \right) y(t) - \gamma \lambda(t) \nabla^2 f(x(t)) y(t) - \omega(t) \nabla f(x(t)) $$
Then, a Lie-Trotter splitting discretization, with a forward Euler method on the first subproblem and a symplectic integrator on the second one leads to \ref{LT-S-IGAHD} algorithm.
\end{theorem}

\begin{Proof}
As in the case of (\ref{IGAHD}) we consider two subproblems, one containing the so-called non-potential forces and the other one containing the potential forces :
\begin{align*}
\hspace*{-0.5cm}
(P_1) :
\begin{cases}
\dot{x}(t) = 0 \\
\dot{v}(t) = - \dfrac{\alpha}{t} v(t) - \dfrac{\lambda(t)}{\gamma} \nabla f^2 (x(t)) v(t) - \dfrac{\omega(t)}{\gamma^2} \nabla f(x(t)) + \dfrac{\tilde{\gamma}(t)}{\gamma^2} \nabla f(x(t)) - \nabla f \left( \mathcal{A}(t, x(t), v(t)) \right) + \nabla f(x(t))
\end{cases}
\end{align*}
and the second subproblem as
\begin{align*}
(P_2) :
\begin{cases}
\dot{x}(t) = v(t) \\
\dot{v}(t) = - \nabla f(x(t))
\end{cases}
\end{align*}
Taking $\gamma = h$, approximating $\tilde{\gamma}(t_n) \approx \gamma_n$, $\lambda(t_n) \approx \lambda_n$ and $\omega(t_n) \approx \omega_n$, and defining $\alpha_n := 1 - \dfrac{\alpha h}{t_n}$ we obtain that
\begin{align*}
\begin{cases}
x_{n+1/2} = x_n \\
v_{n+1/2} = \alpha_n v_n - \lambda_n \nabla^2 f(x_n) v_n - \dfrac{\omega_n}{h} \nabla f(x_n) + \dfrac{\gamma_n}{h} \nabla f(x_n) + h \nabla f(x_n) - h \nabla f \left( y_n \right)
\end{cases}
\end{align*}
where we have the notation for the auxiliary iteration element $y_n := \mathcal{A}(t_n, x_n, v_n)$, i.e. $y_n = x_n + h \alpha_n v_n - h \lambda_n \nabla^2 f(x_n) v_n - \omega_n \nabla f(x_n)$.
Furthermore, for the second subproblem $(P_2)$, we apply the symplectic Euler integrator (\ref{SE_2}) and obtain that
\begin{align*}
\begin{cases}
x_{n+1} = x_{n+1/2} + h v_{n+1} \\
v_{n+1} = v_{n+1/2} - h \nabla f(x_{n+1/2})
\end{cases}
\end{align*}
Using the fact that for $t_{n} = nh$, the splitting solution is the discrete solution of $(P_2)$, namely $(x_{n+1}, v_{n+1})^T \in \mathbb{R}^N \times \mathbb{R}^N$ , we obtain, from the definition of $x_{n+1}$, that $v_n = \dfrac{x_n-x_{n-1}}{h}$. Using the approximation technique for $\nabla^2 f(x_n) v_n$ as in Theorem \ref{T4}, it follows that our algorithm is identical with \ref{LT-S-IGAHD}.
\end{Proof}

Now, we will briefly give a chain of observations that will give insights into the construction of \ref{LT-S-IGAHD} type algorithms.

\begin{remark}\label{R8}
The same principle as derived in Theorem \ref{T7} can be applied to the \ref{DynSys-ARDM}. It is easy to see that if we add a perturbation of the form $\dfrac{\tilde{\gamma}(t)}{\gamma^2} \nabla f(x(t))$, the we will obtain an additional element $\gamma_n \nabla f(x_n)$ in the discrete version of the exact solution, i.e. $x_{n+1}$. Furthermore, we will explain in detail the technical reason why the new added iteration is significant in the process of creating new splitting-type algorithms. For example, we refer to the particular case of the dynamical system from Theorem \ref{T7}, for the case when $\omega(t) = \dfrac{\beta \gamma^2}{t}$. This is the system corresponding to the IGAHD dynamical system from Theorem \ref{T4}, but with an additional gradient with coefficient $\dfrac{\tilde{\gamma}(t)}{\gamma^2}$. Even though we have two gradients in the right hand side, namely $\dfrac{\beta}{t} \nabla f(x(t))$ and $\dfrac{\tilde{\gamma}(t)}{\gamma^2} \nabla f(x(t))$, the first one is proportional to the one under the gradient of $\mathcal{A}(t, x(t), \dot{x}(t))$, i.e. $ \dfrac{\beta \gamma^2}{t} \nabla f(x(t))$. This basic observation (also valid for the term involving the Hessian, which also has $\gamma^2$ as a coefficient of proportionality) leads to the crucial idea in which we can construct algorithms of the form \ref{LT-S-IGAHD}. This is equivalent to the fact that if a term proportional to the gradient at the current value does not appear under the gradient of $\mathcal{A}(t,x(t),x, \dot{x}(t))$, then this term  will not appear in the inertial term $y_n$ but will appear in $x_{n+1}$, i.e. it is added to $y_n$ and $\nabla f(y_n)$.
\end{remark}

\begin{remark}\label{R9}
In the Subsection \ref{S24} we have pointed out that in order to make the algorithms have the natural behavior in the sense that they do not possess exploding gradients or spurious fixed points, a natural technique is to make the division into potential and non-potential forces and add and substract $\theta(t) \nabla f(x(t))$ as in the case of \ref{LT-SE3}. This method was used for the \ref{LT-SE3} algorithm by the fact that we have employed the \ref{SE_1} symplectic integrator. In the case of \ref{LT-S-IGAHD} (as in the case of Nesterov's method \ref{AGM2}), the symplectic numerical scheme \ref{SE_2} does not require such a process.
\end{remark}

\begin{remark}\label{R10}
One might argue that our Hessian-type algorithms \ref{LT-S-IGAHD} have a major drawback, namely we require the computation of $\nabla f(x_n)$ at $x_{n+1}$, along with $\nabla f(y_n)$. This is not true since $\nabla f(x_n)$ was already computed in order to find the value of the inertial term $y_n$ so, the algorithms require the same computational time.
\end{remark}

\begin{remark}\label{R11}
One might wonder why we do not analyzed theoretically \ref{LT-Symplectic} type algorithms and why we have introduced an additional gradient value $\nabla f(x_{n-1})$. The motivation is two-fold. First of all, the product $\nabla^2 f(x(t)) \dot{x}(t)$ is equal to $\dfrac{d}{dt} \left( \nabla f(x(t)) \right)$ and so it is natural to consider the effects of the Hessian perturbation. Also, the Hessian-type systems lead to a much clearer interpretation of numerical methods in optimization (for this see Subsection 2.1). Second of all, the previous value under the gradient $\nabla f(x_{n-1})$ leads to an intuitive discrete Lyapunov function as we will see in the main result of the next subsection. Without this term, through tedious computations, one can see that the Lyapunov function can not be naturally chosen.
\end{remark}

\begin{remark}\label{R12}
Recently, in \cite{ChenLuo}, L. Chen and H. Luo have presented some Hessian-driven damping dynamical systems. They have incorporated the effects of a dynamically changing damping coefficient. Also, they have considered alternative proofs concerning the continuous and the discrete energy functions using the so-called strong Lyapunov property. On the other hand, their optimization methods are in connection to the first order system associated to the second order evolution equation containing the Hessian. This idea is related to the technique of A.M. Oberman and M. Laborde \cite{Oberman}, in which Nesterov type accelerated gradient methods can be interpreted as explicit discretizations of Hessian systems, but with an extra gradient descent step. We point out that the algorithms developed in \cite{ChenLuo} come from an interpretation different to that developed by us in Subsection \ref{S21}. Our algorithms have a beautiful geometric interpretation related to the idea of symplectic integrators that almost conserve the energy and splitting methods that separate the potential and non-potential forces. Also, they have presented that Nesterov's method is an explicit method for a Hessian type system but with an addition gradient step. This can be seen in contrast with the present work where we actually show in a rigorous manner how \ref{AGM2} can be constructed through the lens of the Lie-Trotter splitting.
\end{remark}

\subsection{The rate of convergence of the optimization methods}\label{S32}

In the present section, we shall consider some necessary tools for the convergence analysis of the Hessian-driven damping discretization \ref{LT-S-IGAHD}. For this, we recall that an important assumption concerning the theoretical rate of convergence of optimization algorithms that the obhective function has a Lipschitz continuous gradient. An interesting consequence is the so-called descent lemma (see inequality 2.1.6 from page 56 of \cite{Nesterov_Book}), which is presented below.

\begin{lemma}[Descent lemma]\label{L13}
Let $f : \mathbb{R}^N \to \mathbb{R}$ be F\' rechet differentiable, such that $\nabla f : \mathbb{R}^N \to \mathbb{R}^N$ is $L$ - Lipschitz continuous. Then, one has that
\begin{align}\label{DL}\tag{DL}
f(y) \leq f(x) + \langle \nabla f(x), y-x \rangle + \dfrac{L}{2} \| y-x \|^2 \quad , \quad \forall x,y \in \mathbb{R}^N \quad .
\end{align}
\end{lemma}

It is worth noticing that if the objective function is convex, then we have actually an equivalence between \ref{DL} and the $L$ - Lipschitz continuity of $\nabla f$. Furthermore, Lemma \ref{DL} presents a quadratic upper bound of the underlying objective function. Now, we shall briefly recall from \cite{Attouch_Chbani} extended variant of Lemma \ref{L13}, where a modified variant of \ref{DL} contains also the norm of the difference between two values of $\nabla f$ (as a side note, see also inequality 2.17 from page 57 of \cite{Nesterov_Book}).

\begin{lemma}[Extended descent lemma]\label{L14}
Let $f : \mathbb{R}^N \to \mathbb{R}$ be convex and Fr\' echet differentiable, such that $\nabla f : \mathbb{R}^N \to \mathbb{R}^N$ is $L$ - Lipschitz continuous. Consider a parameter $s$, such that $0 < s \leq \dfrac{1}{L}$. Then one has that
\begin{align}\label{EDL}\tag{EDL}
f(y - s \nabla f(y)) \leq f(x) + \langle \nabla f(y), y-x \rangle - \dfrac{s}{2} \| \nabla f(y) \|^2 - \dfrac{s}{2} \| \nabla f(x) - \nabla f(y) \|^2 \quad , \quad \forall x,y \in \mathbb{R}^N \quad .
\end{align}
\end{lemma}

In the following sequel, we present a technical lemma that will play a crucial role in our main result regarding the asymptotic convergence of the Hessian-type algorithms. So, it is the time to 'extend' the 'Extended descent lemma\ (Lemma \ref{L14}). This represents our first result of this section.

\begin{lemma}[Three points extended lemma]\label{L15}
Let $f : \mathbb{R}^N \to \mathbb{R}$ be a convex and Fr\' echet differentiable function, such that $\nabla f:\mathbb{R}^N \to \mathbb{R}^N$ is $L$ - Lipschitz continuous. Let $s \in \left( 0 , \dfrac{1}{L} \right]$. For simplicity, consider the following notations:
$ \mathcal{A}_1 (s) := s $, $ \mathcal{A}_2 (s) := - \mathcal{A}_1 (s) $, $ \mathcal{A}_3 (s ; \gamma, L) := - \gamma (1- L s) $, $ \mathcal{A}_4 (s) := \mathcal{A}_1 (s) - \dfrac{s}{2} $ and $ \mathcal{A}_5 (s ; \gamma) := - \gamma^2 / (2s) $, where $\gamma \geq 0$. Then $\forall x,y,z \in \mathbb{R}^N$, one has that $f$ satisfies the following inequality:
\begin{align}\label{EEDL}\tag{E-EDL}
f(y - s \nabla f(y) + \gamma \nabla f(z)) &\leq f(x) + \langle \nabla f(y), y-x \rangle - \mathcal{A}_1 (s) \| \nabla f(y) \|^2 \\
&- \left[ \mathcal{A}_2 (s) \langle \nabla f(y) , \nabla f(x) \rangle + \mathcal{A}_3 (s ; \gamma, L) \langle \nabla f(y), \nabla f(z) \rangle \right] \nonumber \\
& - \left[ \mathcal{A}_4 (s) \| \nabla f(x) \|^2 + \mathcal{A}_5 (s ; \gamma) \| \nabla f(z) \|^2 \right] \nonumber \quad .
\end{align}
\end{lemma}

\begin{Proof}
Let $x, y, z \in \mathbb{R}^N$. For simplicity, consider $y^+ := y - s \nabla f(y) + \gamma \nabla f(z)$. Applying \ref{DL} for $y^+$ and $y$, it follows that 
\begin{align*}
f(y^+) \leq f(y) + \langle \nabla f(y), y^+ - y \rangle + \dfrac{L}{2} \| y^+ - y \|^2 \quad .
\end{align*}
This leads to
\begin{align*}
& f(y^+) \leq f(y) + \langle \nabla f(y), - s \nabla f(y) + \gamma \nabla f(z) \rangle + \dfrac{L}{2} \| - s \nabla f(y) + \gamma \nabla f(z) \|^2 \\
& f(y^+) \leq f(y) - s \| \nabla f(y) \|^2 + \gamma \langle \nabla f(y), \nabla f(z) \rangle + \dfrac{L}{2} \langle -s \nabla f(y) + \gamma \nabla f(z) , -s \nabla f(y) + \gamma \nabla f(z) \rangle \\
& f(y^+) \leq f(y) - s \| \nabla f(y) \|^2 + \gamma \langle \nabla f(y), \nabla f(z) \rangle + \dfrac{L}{2} s^2 \| \nabla f(y) \|^2 + \dfrac{L}{2} \gamma^2 \| \nabla f(z) \|^2 - Ls \gamma \langle \nabla f(y), \nabla f(z) \rangle
\end{align*}
So, it follows that
\begin{align}\label{eq1_EEDL}
& f(y^+) \leq f(y) - \dfrac{s}{2} \left( 2 - Ls \right) \| \nabla f(y) \|^2 + \dfrac{L}{2} \gamma^2 \| \nabla f(z) \|^2 + \gamma (1-Ls) \langle \nabla f(y), \nabla f(z) \rangle
\end{align}
Taking into account that $\nabla f$ is $L$ - Lipschitz continuous is equivalent to the fact that the conjugate function $f^\ast$ is $\dfrac{1}{L}$ - strongly convex and applying this property for the pair $(\nabla f(y), \nabla f(x))$, we get that
\begin{align*}
f^\ast (\nabla f(y)) \geq f^\ast (\nabla f(x)) + \langle \nabla f^\ast (\nabla f(x)) , \nabla f(y) - \nabla f(x) \rangle + \dfrac{1}{2L} \| \nabla f(y) - \nabla f(x) \|^2
\end{align*}
Using the fact that $(\nabla f)^{-1} = \nabla f^\ast$, it simplifies to
\begin{align*}
f^\ast (\nabla f(y)) \geq f^\ast (\nabla f(x)) + \langle x , \nabla f(y) - \nabla f(x) \rangle + \dfrac{1}{2L} \| \nabla f(y) - \nabla f(x) \|^2
\end{align*}
Using Fenchel identity for $y$, namely $f(y) = \langle \nabla f(y), y \rangle - f^\ast (\nabla f(y))$, we get
\begin{align*}
& \langle \nabla f(y) , y \rangle - f(y) \geq  f^\ast (\nabla f(x)) + \langle x , \nabla f(y) - \nabla f(x) \rangle + \dfrac{1}{2L} \| \nabla f(y) - \nabla f(x) \|^2
\end{align*}
So, this is equivalent to
\begin{align*}
& f(y) \leq \langle \nabla f(y), y \rangle - f^\ast (\nabla f(x)) - \langle x , \nabla f(y) - \nabla f(x) \rangle - \dfrac{1}{2L} \| \nabla f(y) - \nabla f(x) \|^2 \\
& f(y) \leq \left[ \langle x , \nabla f(x) \rangle - f^\ast (\nabla f(x)) \right] + \langle \nabla f(y) , y - x \rangle - \dfrac{1}{2L} \| \nabla f(y) - \nabla f(x) \|^2
\end{align*}
Again using Fenchel identity but applied on $x$, namely $f(x) = \langle x , \nabla f(x) \rangle - f^\ast (\nabla f(x))$, it follows that
\begin{align}\label{eq2_EEDL}
f(y) \leq f(x) + \langle \nabla f(y), y - x \rangle  - \dfrac{1}{2L} \| \nabla f(y) - \nabla f(x) \|^2 \quad .
\end{align}
Taking into account inequalities \ref{eq1_EEDL} and \ref{eq2_EEDL}, thus we obtain
\begin{align}\label{eq3_EEDL}
f(y^+) & \leq f(x) + \langle \nabla f(y), y-x \rangle - \dfrac{1}{2L} \| \nabla f(y) - \nabla f(x) \|^2 - \dfrac{s}{2} \left( 2 - Ls \right) \| \nabla f(y) \|^2 \\
&+ \dfrac{L}{2} \gamma^2 \| \nabla f(z) \|^2 + \gamma (1-Ls) \langle \nabla f(z) , \nabla f(y) \rangle \nonumber \quad .
\end{align}
From the theorem's assumptions, we know that $ \dfrac{s}{2} \leq s - \dfrac{Ls^2}{2}$, $\dfrac{1}{L} \geq s$ and $- \dfrac{1}{2L} \leq - \dfrac{s}{2}$. From the inequality \ref{eq3_EEDL} and from the condition on $L$, one obtains that
\begin{align*}
f(y^+) & \leq f(x) + \langle \nabla f(y) , y - x \rangle - \dfrac{s}{2} \| \nabla f(y) - \nabla f(x) \|^2 - \dfrac{s}{2} \| \nabla f(y) \|^2 \\ 
&+ \dfrac{L \gamma^2}{2} \| \nabla f(z) \|^2 + \gamma (1-Ls) \langle \nabla f(z) , \nabla f(y) \rangle
\end{align*}
Through some simple computations, it follows
\begin{align*}
f(y^+) & \leq f(x) + \langle \nabla f(y) , y - x \rangle - \dfrac{s}{2} \| \nabla f(y) \|^2 - \dfrac{s}{2} \| \nabla f(x) \|^2 \\
& + s \langle \nabla f(y), \nabla f(x) \rangle - \dfrac{s}{2} \| \nabla f(y) \|^2
+ \dfrac{L \gamma^2}{2} \| \nabla f(z) \|^2 + \gamma (1-Ls) \langle \nabla f(z) , \nabla f(y) \rangle
\end{align*}
Hence,
\begin{align*}
f(y^+) &\leq f(x) + \langle \nabla f(y), y-x \rangle - \left[ \dfrac{s}{2} + \dfrac{s}{2} \right] \| \nabla f(y) \|^2 - \left[ \dfrac{s}{2} \| \nabla f(x) \|^2 - \dfrac{L \gamma^2}{2} \| \nabla f(z) \|^2 \right] \\
& - \left[ - s \langle \nabla f(y), \nabla f(x) \rangle - \gamma (1-Ls) \langle \nabla f(y), \nabla f(z) \rangle) \right]
\end{align*}
Using the fact $\dfrac{L}{2} \leq \dfrac{1}{2s} $, the last inequality leads to
\begin{align*}
f(y^+) &\leq f(x) + \langle \nabla f(y), y-x \rangle - s \| \nabla f(y) \|^2 - \left[ \dfrac{s}{2} \| \nabla f(x) \|^2 - \dfrac{ \gamma^2}{2s}  \| \nabla f(z) \|^2 \right] \\
& - \left[ - s \langle \nabla f(y), \nabla f(x) \rangle - \gamma (1-Ls) \langle \nabla f(y), \nabla f(z) \rangle \right]
\end{align*}
Using the notations for $\mathcal{A}_1 (s)$, $\mathcal{A}_2 (s)$, $\mathcal{A}_3 (s; \gamma, L)$, $\mathcal{A}_4 (s)$ and $\mathcal{A}_5 (s; \gamma)$, the inequality \ref{EEDL} follows easily. 
\end{Proof}

\begin{corollary}\label{C16}
Consider $x, y, z \in \mathbb{R}^N$. By the fact that the stepsize $s$ lies in  $ \left( 0 ,\dfrac{1}{L} \right]$, it follows that the major bound \ref{EEDL} for $f(y^+)$ reduces to \ref{EDL}, with $\gamma = 0 $.
\end{corollary}

Before embarking into our main result of the present paper, we consider some technical lemmas that will be used in the forthcoming sequel. These results involve inequalities that are related to quadratic inequalities with nonconstant coefficients. Even though they are elementary results, for brevity we give them along with their proofs.

\begin{lemma}\label{L17}
Let $a,b,c : \mathbb{R} \to \mathbb{R}$ be given functions. If there exists $x \in \mathbb{R}$, such that the following inequalities are satisfied: $a(x)>0$ and $b^2(x) - 4 a(x) c(x) \leq 0$, then $x$ satisfies $a(x) x^2 + b(x) x + c(x) \geq 0$. 
\end{lemma}
 
\begin{Proof}
First of all, we consider the following computations:
\begin{align*}
a(x) \left( x + \dfrac{b(x)}{2a(x)} \right)^2 - \dfrac{b^2(x) - 4 a(x) c(x)}{4a(x)} &= a(x) \left[ x^2 + \dfrac{b^2(x)}{4a^2(x)} + x \dfrac{b(x)}{a(x)} \right] - \dfrac{b^2(x)}{4a(x)} + c(x) \\
&= a(x) x^2 + \dfrac{b^2(x)}{4a^2(x)} a(x) + x \dfrac{b(x)}{a(x)} a(x) - \dfrac{b^2(x)}{4a(x)} + c(x) \, .
\end{align*}
This is equivalent to the fact that
$$ a(x) x^2 + b(x) x + c(x) = a(x) \left( x + \dfrac{b(x)}{2a(x)} \right)^2 - \dfrac{b^2(x)-4a(x)c(x)}{4a(x)} \, . $$
Now, $a(x)x^2 + b(x) x + c(x) \geq 0 \Longleftrightarrow a(x) \left( x + \dfrac{b(x)}{2a(x)} \right)^2 \geq \dfrac{b^2(x)-4a(x)c(x)}{4a(x)}$. Taking into account that $b^2(x) - 4a(x)c(x) \geq 0$ and dividing by $a(x) > 0$, it follows that $\left( x + \dfrac{b(x)}{2a(x)} \right)^2 \geq \dfrac{b^2(x)-4a(x)c(x)}{4a^2(x)}$. Likewise, let $\left( x + \dfrac{b(x)}{2a(x)} \right)^2 \geq \dfrac{b^2(x)-4a(x)c(x)}{4a^2(x)}$. As before, we observe that these terms are positive and multiplying be $a(x) > 0$, we obtain that $a(x) x^2 + b(x) x + c(x) \geq 0$. So, wrapping things up, we need to show that $\left( x + \dfrac{b(x)}{2a(x)} \right)^2 \geq \dfrac{b^2(x)-4a(x)c(x)}{4a^2(x)}$, since this is equivalent to the quadratic inequality with nonconstant coefficients. Because $b^2(x) - a(x) c(x) \leq 0$, then it follows that $\dfrac{b^2(x) - 4 a(x) c(x)}{4 a^2 (x)} \leq 0$. Combining this with the fact that the left hand side of the inequality is positive, the conclusion stands up. 
\end{Proof} 
 
\begin{lemma}\label{L18}
Let $a,b,c : \mathbb{R} \to \mathbb{R}$ be given functions. If there exists $x \in \mathbb{R}$, such that the following inequalities are satisfied: $a(x) > 0$, $ b^2 (x) -  4a(x)c(x) \geq 0$ and $x  \in I_1 \cup I_2$, where: $I_1 := \left( - \infty, \dfrac{-b(x) - \sqrt{b^2(x) - 4 a(x)c(x)}}{2a(x)} \right]$ and $I_2 := \left[ \dfrac{-b(x) + \sqrt{b^2(x) - 4 a(x)c(x)}}{2a(x)}, + \infty \right)$, then $x$ satisfies $a(x) x^2 + b(x) x + c(x) \geq 0$.
\end{lemma} 
 
\begin{Proof}
As before, we must show that $\left( x + \dfrac{b(x)}{2a(x)} \right)^2 \geq \dfrac{b^2(x) - 4 a(x) c(x)}{4 a^2 (x)}$, so we must analyze two cases:
\begin{itemize}
\item[\textbf{1)}] $x \in I_1$, that is equivalent to $x \leq - \dfrac{b(x)}{2a(x)} - \dfrac{\sqrt{b^2(x) - 4 a(x)c(x)}}{2a(x)}$. Thus, we have that $\dfrac{\sqrt{b^2(x) - 4 a(x)c(x)}}{2a(x)} \leq - x - \dfrac{b(x)}{2a(x)} = - \left( x + \dfrac{b(x)}{2a(x)} \right)$. But, it is easy to see that $x \leq - \dfrac{b(x)}{2a(x)}$, i.e. $- \left( x + \dfrac{b(x)}{2a(x)} \right) \geq 0$. So, it follows that $\dfrac{\sqrt{b^2(x) - 4 a(x)c(x)}}{2a(x)} \leq - \left( x + \dfrac{b(x)}{2a(x)} \right) = \Big| x + \dfrac{b(x)}{2a(x)} \Big|$. Squaring the positive terms, it follows that $\dfrac{b^2(x) - 4 a(x)c(x)}{4 a^2(x)} \leq \left( x + \dfrac{b(x)}{2a(x)} \right)^2$.
\item[\textbf{2)}] $x \in I_2$, i.e. $x \geq - \dfrac{b(x)}{2a(x)} + \dfrac{\sqrt{b^2(x) - 4 a(x)c(x)}}{2a(x)}$, which is equivalent to $\dfrac{\sqrt{b^2(x) - 4 a(x)c(x)}}{2a(x)} \leq x + \dfrac{b(x)}{2a(x)}$. But, since $x + \dfrac{b(x)}{2a(x)} \geq \dfrac{\sqrt{b^2(x) - 4 a(x)c(x)}}{2a(x)}$, it follows that $x + \dfrac{b(x)}{2a(x)} \geq 0$, so $x + \dfrac{b(x)}{2a(x)} = \Big| x + \dfrac{b(x)}{2a(x)} \Big|$. Hence, using that $\dfrac{\sqrt{b^2(x) - 4 a(x)c(x)}}{2a(x)} \leq \Big| x + \dfrac{b(x)}{2a(x)} \Big|$ and squaring the positive terms, we get the conclusion that $\dfrac{b^2(x) - 4 a(x)c(x)}{4a^2(x)} \leq \left( x + \dfrac{b(x)}{2a(x)} \right)^2$.
\end{itemize}
\end{Proof} 
 
Until our main result, we consider some observations regarding the functions that appear in Lemma \ref{L17} and Lemma \ref{L18}.
\begin{remark}\label{R19}
Related to Lemma \ref{L17} and Lemma \ref{L18}, we have the following:
\begin{itemize}
\item For a given pair of functions $(a(x),b(x),c(x))$, if the element $x \in \mathbb{R}$ exists, then the quadratic inequality takes place. But this does not mean that such a point exists for all pair of functions $(a(x), b(x), c(x))$.
\item It is easy to observe that, in Lemma \ref{L17}, the conditions $a(x) > 0$ and $b^2(x) - 4 a(x) c(x) \leq 0$ implies that $c(x) \geq 0$. Furthermore, we do not have a positivity-type inequality on the function $b = b(x)$.
\item In Lemma \ref{L18}, we observe that $a(x) > 0$ and if $c(x) \leq 0$, then it implies that $4a(x)c(x) \leq b^2 (x)$ is trivially satisfied. As before, we do not need the positivity of $b=b(x)$.
\end{itemize}
\end{remark}

In what follows we have our main result of this section. Given a convex objective function endowed with Lipschitz continuous gradient, we show that the rate of convergence of the algorithms given by the general discretization \ref{LT-S-IGAHD} is $O \left( \dfrac{1}{n^2} \right)$ as $n \to \infty$, under some suitable conditions over the parameters. Given $x^\ast \in \argmin f$ (under the assumption that $\argmin f \neq \emptyset$), the Lyapunov analysis that we employ is based upon the discrete energy defined as:
\begin{align}\label{Energy}\tag{Discrete-Energy}
& \mathcal{E}_n := t_{n}^{2} \left( f(x_n) - f(x^\ast) \right) + \dfrac{1}{2s} \| z_n \|^2 \\
& z_n := (x_{n-1} - x^\ast) + t_{n} (x_n - x_{n-1}) + \lambda_n t_{n+1} \nabla f(x_{n-1}) \, , \nonumber
\end{align}
where $t_{n+1} = \dfrac{n}{\alpha - 1}$, with $\alpha \geq 3$. We shall also use $\alpha_n = \dfrac{n-\alpha}{n}$, so it follows that $t_n - 1 = \alpha_n t_{n+1}$. Moreover, our technique for the proof is inspired by the ones given in \cite{Attouch_Chbani} and \cite{AttouchCabot}.

\begin{theorem}\label{T20}
Let $f : \mathbb{R}^N \to \mathbb{R}$ be a convex and Fr\' echet differentiable function whose gradient $\nabla f : \mathbb{R}^N \to \mathbb{R}^N$ is $L$ - Lipschitz continuous. Also, suppose that $\argmin f \neq \emptyset$. Let $x_n$ be generated by the algorithm \ref{LT-S-IGAHD} for every $n \geq 1$ and, for simplicity, denote $h^2$ as $s$, which satisfies the condition $s \in \left( 0, \dfrac{1}{L} \right]$. Suppose that the sequences $(\gamma_n)_{n \in \mathbb{N}}$ and $(\lambda_n)_{n \in \mathbb{N}}$ are convergent. Likewise, assume that the following assumptions are satisfied: 
\begin{itemize}
\item[i)] there exists $N^\prime \in \mathbb{N}$ such that $\left[ s(1-\gamma_n L) - \dfrac{n+1}{n} \lambda_{n+1} \right]^2 < s^2 - \gamma_n^2 \, , \forall n > max \lbrace \alpha - 1, N^\prime \rbrace$; 
\item[ii)] $\gamma_n = (\lambda_n + \omega_n) - \dfrac{n+1}{n} \lambda_{n+1} \, , \forall n \geq \alpha - 1$.
\end{itemize}
For $x^\ast \in \argmin f$, the sequence $(\mathcal{E})_{n \in \mathbb{N}}$ defined in \ref{Energy} is non-increasing and the following converge rate is satisfied:
$$ f(x_n) - f(x^\ast) = O \left( \dfrac{1}{n^2} \right) \text{ as } n \to \infty$$
\end{theorem}

\begin{Proof}
\textbf{( Step I )} The first step of the proof is to apply the so-called Three points extended lemma, namely Lemma \ref{L15} and then construct the difference $\mathcal{E}_{n+1} - \mathcal{E}_n$ of the discrete energy \ref{Energy}. \\
We apply the inequality \ref{EEDL} with the pair $x = z = x_n$, $y = y_n$ and $\gamma = \gamma_n$. Then, we obtain that:
\begin{align*}
f(y_n - s \nabla f(y_n) + \gamma_n \nabla f(x_n)) &\leq f(x_n) + \langle \nabla f(y_n), y_n - x_n \rangle - \mathcal{A}_1 (s) \| \nabla f(y_n) \|^2 \\
&- \left[ \mathcal{A}_2 (s) \langle \nabla f(y_n), \nabla f(x_n) \rangle + \mathcal{A}_3 (s; \gamma_n, L) \langle \nabla f(y_n), \nabla f(x_n) \rangle \right] \\
&- \left[ \mathcal{A}_4 (s) \| \nabla f(x_n) \|^2 + \mathcal{A}_5 (s; \gamma_n) \| \nabla f(x_n) \|^2 \right] \, .
\end{align*}
This means that
\begin{align}\label{ineq1:Proof}
f(x_{n+1}) &\leq f(x_n) + \langle \nabla f(y_n), y_n - x_n \rangle - \mathcal{A}_1 (s) \| \nabla f(y_n) \|^2 \\
&- \left[ \mathcal{A}_2 (s) + \mathcal{A}_3 (s; \gamma_n, L) \right] \langle \nabla f(y_n), \nabla f(x_n) \rangle  \nonumber \\
&- \left[ \mathcal{A}_4 (s) + \mathcal{A}_5 (s; \gamma_n) \right] \| \nabla f(x_n) \|^2 \, . \nonumber
\end{align}
As before, we will apply \ref{EEDL}, but with the choices $x = x^\ast$, $y = y_n$, $z = x_n$ and $\gamma = \gamma_n$. It follows that:
\begin{align*}
f(y_n - s \nabla f(y_n) + \gamma_n \nabla f(x_n)) &\leq f(x^\ast) + \langle \nabla f(y_n), y_n - x^\ast \rangle - \mathcal{A}_1 (s) \| \nabla f(y_n) \|^2 \\
&- \left[ \mathcal{A}_2 (s) \langle \nabla f(y_n), \nabla f(x^\ast) \rangle + \mathcal{A}_3 (s; \gamma_n, L) \langle \nabla f(y_n), \nabla f(x_n) \rangle \right] \\
&- \left[ \mathcal{A}_4 (s) \| \nabla f(x^\ast) \|^2 + \mathcal{A}_5 (s; \gamma_n) \| \nabla f(x_n) \|^2 \right] \, .
\end{align*}
Using the fact that $x^\ast \in \argmin f$, so $\nabla f(x^\ast) = 0$, we obtain
\begin{align}\label{ineq2:Proof}
f(x_{n+1}) &\leq f(x^\ast) + \langle \nabla f(y_n), y_n - x^\ast \rangle - \mathcal{A}_1 (s) \| \nabla f(y_n) \|^2 \\
&- \mathcal{A}_3 (s; \gamma_n, L) \langle \nabla f(y_n), \nabla f(x_n) \rangle  - \mathcal{A}_5 (s; \gamma_n) \| \nabla f(x_n) \|^2 \, . \nonumber
\end{align}
Using the definition of $t_{n+1}$, we consider 
\begin{align}\label{n1}
n \geq \alpha-1 \, .
\end{align}
hence $t_{n+1} - 1 \geq 0$. Now, we multiply \ref{ineq1:Proof} by $(t_{n+1} - 1)$ and then add it to \ref{ineq2:Proof}. One obtains that
\begin{align*}
(t_{n+1}-1)f(x_{n+1}) + f(x_{n+1}) &\leq (t_{n+1}-1)f(x_n) + f(x^\ast) + (t_{n+1}-1) \langle \nabla f(y_n), y_n-x_n \rangle + \langle \nabla f(y_n), y_n - x^\ast \rangle\\
& -(t_{n+1}-1) \left[ \mathcal{A}_2 (s) + \mathcal{A}_3 (s;\gamma_n,L) \right] \langle \nabla f(y_n), \nabla f(x_n) \rangle - \mathcal{A}_1 (s) (t_{n+1}-1) \| \nabla f(y_n) \|^2 \\
& - (t_{n+1}-1) \left[ \mathcal{A}_4 (s) + \mathcal{A}_5 (s;\gamma_n) \right] \| \nabla f(x_n) \|^2 - \mathcal{A}_1 (s) \| \nabla f(y_n) \|^2 \\
&- \mathcal{A}_3 (s;\gamma_n,L) \langle \nabla f(y_n), \nabla f(x_n) \rangle - \mathcal{A}_5 (s; \gamma_n) \| \nabla f(x_n) \|^2 \, .
\end{align*}
By grouping up the terms, it follows that
\begin{align*}
t_{n+1}f(x_{n+1}) - (t_{n+1}-1)f(x_n) - f(x^\ast) &\leq (t_{n+1}-1) \langle \nabla f(y_n), y_n-x_n \rangle + \langle \nabla f(y_n), y_n - x^\ast \rangle\\
& - \left[ (t_{n+1}-1) \mathcal{A}_2 (s) + t_{n+1} \mathcal{A}_3 (s;\gamma_n,L) \right] \langle \nabla f(y_n), \nabla f(x_n) \rangle \\
& - \left[ (t_{n+1}-1) \mathcal{A}_4 (s) + t_{n+1} \mathcal{A}_5 (s;\gamma_n) \right] \| \nabla f(x_n) \|^2 - \mathcal{A}_1 (s) t_{n+1} \| \nabla f(y_n) \|^2 \, .
\end{align*}
We rewrite the left hand side of the inequality from above as follows:
\begin{align*}
t_{n+1}f(x_{n+1}) - (t_{n+1}-1)f(x_n) - f(x^\ast) &= t_{n+1} f(x_{n+1}) - (t_{n+1}-1) f(x_n) - (t_{n+1} - t_{n+1} + 1) f(x^\ast) \\
&= t_{n+1} f(x_{n+1}) - (t_{n+1}-1) f(x_n) - t_{n+1} f(x^\ast) - (1-t_{n+1}) f(x^\ast) \\
&= t_{n+1} \left[ f(x_{n+1}) - f(x^\ast) \right] - (t_{n+1}-1) \left[ f(x_{n}) - f(x^\ast) \right]
\end{align*}
Then, the above inequality takes the following form:
\begin{align}\label{ineq3:Proof}
t_{n+1} \left[ f(x_{n+1}) - f(x^\ast) \right] &\leq (t_{n+1}-1) \left[ f(x_{n}) - f(x^\ast) \right] + \langle \nabla f(y_n), (t_{n+1}-1)(y_n-x_n) + (y_n - x^\ast) \rangle\\ \nonumber
& - \left[ (t_{n+1}-1) \mathcal{A}_2 (s) + t_{n+1} \mathcal{A}_3 (s;\gamma_n,L) \right] \langle \nabla f(y_n), \nabla f(x_n) \rangle \\ \nonumber
& - \left[ (t_{n+1}-1) \mathcal{A}_4 (s) + t_{n+1} \mathcal{A}_5 (s;\gamma_n) \right] \| \nabla f(x_n) \|^2 - \mathcal{A}_1 (s) t_{n+1} \| \nabla f(y_n) \|^2 \, . \nonumber
\end{align}
Using the fact that $t_{n+1} \geq 0$, for every $n \geq 0$, hence for every $n$ satisfying \ref{n1}, we multiply the inequality \ref{ineq3:Proof} by $t_{n+1}$ and it follows that:
\begin{align*}
t_{n+1}^2 \left[ f(x_{n+1}) - f(x^\ast) \right] &\leq (t_{n+1}^2-t_{n+1}) \left[ f(x_{n}) - f(x^\ast) \right] + t_{n+1} \langle \nabla f(y_n), (t_{n+1}-1)(y_n-x_n) + (y_n - x^\ast) \rangle\\ \nonumber
& - t_{n+1} \left[ (t_{n+1}-1) \mathcal{A}_2 (s) + t_{n+1} \mathcal{A}_3 (s;\gamma_n,L) \right] \langle \nabla f(y_n), \nabla f(x_n) \rangle \\ \nonumber
& - t_{n+1} \left[ (t_{n+1}-1) \mathcal{A}_4 (s) + t_{n+1} \mathcal{A}_5 (s;\gamma_n) \right] \| \nabla f(x_n) \|^2 - \mathcal{A}_1 (s) t_{n+1}^2 \| \nabla f(y_n) \|^2 \, . \nonumber
\end{align*}
Since $(t_{n+1}^2-t_{n+1}) \left[ f(x_{n}) - f(x^\ast) \right] = (t_{n+1}^2 - t_{n+1} - t_n^2) \left[ f(x_{n}) - f(x^\ast) \right] + t_n^2 \left[ f(x_{n}) - f(x^\ast) \right]$ and using the fact that $t_{n+1}^2 - t_{n+1} - t_n^2 \leq 0$ and that $f(x_n) - f(x^\ast) \geq 0$ (since $x^\ast \in \argmin f$) for every $n \geq 0$ and thus for each $n$ satisfying \ref{n1}, it follows
\begin{align*}
t_{n+1}^2 \left[ f(x_{n+1}) - f(x^\ast) \right] &\leq t_n^2 \left[ f(x_{n}) - f(x^\ast) \right] + t_{n+1} \langle \nabla f(y_n), (t_{n+1}-1)(y_n-x_n) + (y_n - x^\ast) \rangle\\ \nonumber
& - t_{n+1} \left[ (t_{n+1}-1) \mathcal{A}_2 (s) + t_{n+1} \mathcal{A}_3 (s;\gamma_n,L) \right] \langle \nabla f(y_n), \nabla f(x_n) \rangle \\ \nonumber
& - t_{n+1} \left[ (t_{n+1}-1) \mathcal{A}_4 (s) + t_{n+1} \mathcal{A}_5 (s;\gamma_n) \right] \| \nabla f(x_n) \|^2 - \mathcal{A}_1 (s) t_{n+1}^2 \| \nabla f(y_n) \|^2 \, . \nonumber
\end{align*}
Using the definition of $\mathcal{E}_n$ given in \ref{Energy}, one can easily see that
\begin{align}\label{EnergyLHS}
\mathcal{E}_{n+1} - \mathcal{E}_n &\leq t_{n+1} \langle \nabla f(y_n), (t_{n+1}-1)(y_n-x_n) + (y_n - x^\ast) \rangle \\ \nonumber
& - t_{n+1} \left[ (t_{n+1}-1) \mathcal{A}_2 (s) + t_{n+1} \mathcal{A}_3 (s;\gamma_n,L) \right] \langle \nabla f(y_n), \nabla f(x_n) \rangle \\ \nonumber
& - t_{n+1} \left[ (t_{n+1}-1) \mathcal{A}_4 (s) + t_{n+1} \mathcal{A}_5 (s;\gamma_n) \right] \| \nabla f(x_n) \|^2 - \mathcal{A}_1 (s) t_{n+1}^2 \| \nabla f(y_n) \|^2 \nonumber \\
& + \dfrac{1}{2s} \| z_{n+1} \|^2 - \dfrac{1}{2s} \| z_n \|^2 \, . \nonumber
\end{align}
In computing the last term involving $z_n$ and $z_{n+1}$, we shall use the following well-known equality:
\begin{align}\label{NormEquality}
\dfrac{1}{2} \| z_{n+1} \|^2 - \dfrac{1}{2} \| z_n \|^2 = \langle z_{n+1} - z_{n}, z_{n+1} \rangle - \dfrac{1}{2} \| z_{n+1} - z_n \|^2 \, . 
\end{align}
Armed with the definition of $z_n$ given in \ref{Energy}, we have the following estimations:
\begin{align*}
z_{n+1} - z_n &= x_n - x^\ast + t_{n+1} (x_{n+1} - x_n) + \lambda_{n+1} t_{n+2} \nabla f(x_n) - x_{n-1} + x^\ast - t_n (x_n - x_{n-1}) - \lambda_n t_{n+1} \nabla f(x_{n-1}) \\
&= (x_n - x_{n-1}) + t_{n+1} (x_{n+1} - x_n) - t_n (x_n - x_{n-1}) + \left[ \lambda_{n+1} t_{n+2} \nabla f(x_n) - \lambda_n t_{n+1} \nabla f(x_{n-1}) \right] \\
&= t_{n+1} (x_{n+1} - x_n) - (t_n - 1) (x_n - x_{n-1}) + \left[ \lambda_{n+1} t_{n+2} \nabla f(x_n) - \lambda_n t_{n+1} \nabla f(x_{n-1}) \right] \, .
\end{align*}
But, we know that $t_n - 1 = \alpha_n t_{n+1}$ for each $n \geq 0$, consequently for every $n$ satisfying \ref{n1}, we get that
\begin{align*}
z_{n+1} - z_n &= t_{n+1} (x_{n+1} - x_n) - \alpha_n t_{n+1} (x_n - x_{n-1}) + \left[ \lambda_{n+1} t_{n+2} \nabla f(x_n) - \lambda_n t_{n+1} \nabla f(x_{n-1}) \right] \\
&= t_{n+1} \left[ x_{n+1} - (x_n + \alpha_n (x_n - x_{n-1})) \right] + \left[ \lambda_{n+1} t_{n+2} \nabla f(x_n) - \lambda_n t_{n+1} \nabla f(x_{n-1}) \right] \, .
\end{align*}
From the definition of \ref{LT-S-IGAHD}, it follows that $x_n + \alpha_n (x_n - x_{n-1}) = y_n + \lambda_n \left[ \nabla f(x_n) - \nabla f(x_{n-1}) \right] + \omega_n \nabla f(x_n) \, .$
Hence, it follows that
\begin{align*}
z_{n+1} - z_n &= t_{n+1} \left[ x_{n+1} - y_n - \lambda_n \left( \nabla f(x_n) - \nabla f(x_{n-1}) \right) - \omega_n \nabla f(x_n) \right] + \left[ \lambda_{n+1} t_{n+2} \nabla f(x_n) - \lambda_n t_{n+1} \nabla f(x_{n-1}) \right] \\
z_{n+1} - z_n &= t_{n+1} (x_{n+1} - y_n) + \left[ - \lambda_n t_{n+1} - \omega_n t_{n+1} + \lambda_{n+1} t_{n+2} \right] \nabla f(x_n) \, .
\end{align*}
But, from \ref{LT-S-IGAHD}, we have that $x_{n+1} - y_n = - s \nabla f(y_n) + \gamma_n \nabla f(x_n)$. This leads to
\begin{align}\label{eq4:Proof}
z_{n+1} - z_{n} = - s t_{n+1} \nabla f(y_n) + \left[ \gamma_n t_{n+1} + \lambda_{n+1} t_{n+2} - \lambda_n t_{n+1} - \omega_n t_{n+1} \right] \nabla f(x_n) \, . 
\end{align}
It is easy to remark that assumption (ii) is equivalent to $\gamma_n = (\lambda_n + \omega_n) - \dfrac{t_{n+2}}{t_{n+1}} \lambda_{n+1}$ for every $n \geq \alpha - 1$, by taking into account that $t_{n+1} = \dfrac{n}{\alpha-1}$ and so $t_{n+2} = \dfrac{n+1}{\alpha-1}$. Now, multiplying by $t_{n+1}$, it follows that $\gamma_n t_{n+1} = (\lambda_n + \omega_n) t_{n+1} - \lambda_{n+1} t_{n+2}$. Using this relationship, \ref{n1} and \ref{eq4:Proof}, for $n \geq \alpha - 1$, we have that
\begin{align}\label{eq5:Proof}
z_{n+1} - z_n = - s t_{n+1} \nabla f(y_n)
\end{align}
At the same time, by using \ref{NormEquality} and \ref{eq5:Proof}, we get
\begin{align*}
\dfrac{1}{2s} \| z_{n+1} \|^2 - \dfrac{1}{2s} \| z_n \|^2 = - \langle t_{n+1} \nabla f(y_n), (x_n - x^\ast) + t_{n+1} (x_{n+1}-x_n) + \lambda_{n+1} t_{n+2} \nabla f(x_n) \rangle - \dfrac{s}{2} t_{n+1}^2 \| \nabla f(y_n) \|^2 
\end{align*}
This means that
\begin{align}\label{VelocityComputation}
\dfrac{1}{2s} \| z_{n+1} \|^2 - \dfrac{1}{2s} \| z_n \|^2 &= - \lambda_{n+1} t_{n+1} t_{n+2} \langle \nabla f(y_n), \nabla f(x_n) \rangle -\dfrac{s}{2} t_{n+1}^2 \| \nabla f(y_n) \|^2 \\
&- t_{n+1} \langle \nabla f(y_n), (x_n - x^\ast) + t_{n+1} (x_{n+1} - x_n) \rangle \, . \nonumber
\end{align}
Combining \ref{EnergyLHS} with \ref{VelocityComputation}, we obtain
\begin{align*}
\mathcal{E}_{n+1} - \mathcal{E}_n &\leq t_{n+1} \langle \nabla f(y_n), (t_{n+1}-1)(y_n-x_n) + (y_n - x^\ast) \rangle \\ \nonumber
& - t_{n+1} \left[ (t_{n+1}-1) \mathcal{A}_2 (s) + t_{n+1} \mathcal{A}_3 (s;\gamma_n,L) \right] \langle \nabla f(y_n), \nabla f(x_n) \rangle \\ \nonumber
& - t_{n+1} \left[ (t_{n+1}-1) \mathcal{A}_4 (s) + t_{n+1} \mathcal{A}_5 (s;\gamma_n) \right] \| \nabla f(x_n) \|^2 - \mathcal{A}_1 (s) t_{n+1}^2 \| \nabla f(y_n) \|^2 \\
&- \lambda_{n+1} t_{n+1} t_{n+2} \langle \nabla f(y_n), \nabla f(x_n) \rangle -\dfrac{s}{2} t_{n+1}^2 \| \nabla f(y_n) \|^2 - t_{n+1} \langle \nabla f(y_n), (x_n - x^\ast) + t_{n+1} (x_{n+1} - x_n) \rangle \, . 
\end{align*}
On the other hand, we make the following estimations:
\begin{align*}
& \langle \nabla f(y_n), (t_{n+1}-1)(y_n-x_n) + (y_n - x^\ast) \rangle - \langle \nabla f(y_n), (x_n - x^\ast) + t_{n+1} (x_{n+1} - x_n) \rangle = \\
& \langle \nabla f(y_n), (t_{n+1}-1)(y_n-x_n) + (y_n - x^\ast) - (x_n - x^\ast) - t_{n+1} (x_{n+1} - x_n) \rangle = \\
& \langle \nabla f(y_n), t_{n+1} (y_n - x_{n+1}) \rangle \, .
\end{align*}
From \ref{LT-S-IGAHD}, we know that $y_n - x_{n+1} = s \nabla f(y_n) - \gamma_n \nabla f(x_n)$ and this leads to
\begin{align*}
& \langle \nabla f(y_n), (t_{n+1}-1)(y_n-x_n) + (y_n - x^\ast) \rangle - \langle \nabla f(y_n), (x_n - x^\ast) + t_{n+1} (x_{n+1} - x_n) \rangle = \\
& \langle \nabla f(y_n), s t_{n+1} \nabla f(y_n) - \gamma_n t_{n+1} \nabla f(x_n) \rangle = s t_{n+1} \| \nabla f(y_n) \|^2 - \gamma_n t_{n+1} \langle \nabla f(y_n) , \nabla f(x_n) \rangle \, .
\end{align*}
So, it follows that
\begin{align*}
\mathcal{E}_{n+1} - \mathcal{E}_n &\leq - \mathcal{A}_1 (s) t_{n+1}^2 \| \nabla f(y_n) \|^2 + s t_{n+1}^2 \| \nabla f(y_n) \|^2 -\dfrac{s}{2} t_{n+1}^2 \| \nabla f(y_n) \|^2 \\
&- t_{n+1} \left[ (t_{n+1}-1) \mathcal{A}_2 (s) + t_{n+1} \mathcal{A}_3 (s;\gamma_n,L) \right] \langle \nabla f(y_n), \nabla f(x_n) \rangle \\
& - t_{n+1} \left[ (t_{n+1}-1) \mathcal{A}_4 (s) + t_{n+1} \mathcal{A}_5 (s;\gamma_n) \right] \| \nabla f(x_n) \|^2 \\
&- \lambda_{n+1} t_{n+1} t_{n+2} \langle \nabla f(y_n), \nabla f(x_n) \rangle - \gamma_n t_{n+1}^2 \langle f(y_n), f(x_n) \rangle \, . 
\end{align*}
Using the fact that $\mathcal{A}_1 (s) = s$, one obtains
\begin{align*}
\mathcal{E}_{n+1} - \mathcal{E}_n &\leq - t_{n+1} \left[ (t_{n+1}-1) \mathcal{A}_2 (s) + t_{n+1} \mathcal{A}_3 (s;\gamma_n,L) + \lambda_{n+1} t_{n+2} + \gamma_n t_{n+1} \right] \langle \nabla f(y_n), \nabla f(x_n) \rangle \\
& - t_{n+1} \left[ (t_{n+1}-1) \mathcal{A}_4 (s) + t_{n+1} \mathcal{A}_5 (s;\gamma_n) \right] \| \nabla f(x_n) \|^2 - \dfrac{s}{2} t_{n+1}^2 \| \nabla f(y_n) \|^2 \, .
\end{align*}
Recalling the fact that we have already pointed out that assumption (ii) is equivalent to $ \lambda_{n+1} t_{n+2} + \gamma_n t_{n+1} = (\lambda_n + \omega_n) t_{n+1}$ for every $n \geq \alpha - 1$, we find that
\begin{align*}
\mathcal{E}_{n+1} - \mathcal{E}_n &\leq - t_{n+1} \left[ (t_{n+1}-1) \mathcal{A}_2 (s) + t_{n+1} \mathcal{A}_3 (s;\gamma_n,L) + (\lambda_n + \omega_n) t_{n+1} \right] \langle \nabla f(y_n), \nabla f(x_n) \rangle \\
& - t_{n+1} \left[ (t_{n+1}-1) \mathcal{A}_4 (s) + t_{n+1} \mathcal{A}_5 (s;\gamma_n) \right] \| \nabla f(x_n) \|^2 - \dfrac{s}{2} t_{n+1}^2 \| \nabla f(y_n) \|^2 \, .
\end{align*}
So, the difference in the discrete energy functional can be written as
\begin{align}\label{EnergyDecrease}
\mathcal{E}_{n+1} - \mathcal{E}_n \leq - t_{n+1} \tilde{F}_n \, ,
\end{align}
where
\begin{align}\label{Fn}
\tilde{F}_n &:= \dfrac{s}{2} t_{n+1} \| \nabla f(y_n) \|^2 + \left[ (t_{n+1}-1) \mathcal{A}_2 (s) + t_{n+1} \mathcal{A}_3 (s; \gamma_n, L) + (\lambda_n + \omega_n) t_{n+1} \right] \langle \nabla f(y_n), \nabla f(x_n) \rangle \\ \nonumber
&+ \left[ (t_{n+1}-1) \mathcal{A}_4 (s) + t_{n+1} \mathcal{A}_5 (s; \gamma_n) \right] \| \nabla f(x_n) \|^2 \, .
\end{align}
\textbf{( Step II )} The second step of our proof consists in applying Lemma \ref{L17} and Lemma \ref{L18} to two quadratic inequalities with nonconstant coefficients and show the conditions on the index $n$ such that they are satisfied. \\
For simplicity, let's denote $Y := \nabla f(y_n)$ and $X := \nabla f(x_n)$. Hence, we obtain that 
\begin{align*}
\tilde{F}_n = \dfrac{s}{2} t_{n+1} \| Y \|^2 &+ \left[ (t_{n+1}-1) \mathcal{A}_2 (s) + t_{n+1} \mathcal{A}_3 (s; \gamma_n, L) + (\lambda_n+\omega_n) t_{n+1} \right] \langle Y, X \rangle \\
&+ \left[ (t_{n+1}-1) \mathcal{A}_4 (s) + t_{n+1} \mathcal{A}_5 (s; \gamma_n) \right] \| X \|^2 \, . 
\end{align*}
Using Cauchy-Schwarz inequality, it follows that
\begin{align}\label{eq1:Quadratic}
\tilde{F}_n \geq \dfrac{s}{2} t_{n+1} \| Y \|^2 &- \left[ (t_{n+1}-1) \mathcal{A}_2 (s) + t_{n+1} \mathcal{A}_3 (s; \gamma_n, L) + (\lambda_n+\omega_n) t_{n+1} \right] \| Y \| \| X \| \\
&+ \left[ (t_{n+1}-1) \mathcal{A}_4 (s) + t_{n+1} \mathcal{A}_5 (s; \gamma_n) \right] \| X \|^2 := F_n \, . \nonumber
\end{align}
The next step of the present proof is to apply Lemma \ref{L17} and to show that for $n$ large enough, $F_n \geq 0$. Using Lemma \ref{L17}, we observe that the coefficient of $\| Y \|^2$ is positive, namely $\dfrac{s}{2} t_{n+1} > 0$ for $n \geq \alpha - 1$. Now, the only thing left for us is to show that
\begin{align}\label{CondQuadratic}
\left[ (1-t_{n+1}) \mathcal{A}_2 (s) - t_{n+1} \mathcal{A}_3 (s; \gamma_n, L) - (\lambda_n + \omega_n) t_{n+1} \right]^2 \leq 2 s t_{n+1} \left[ (t_{n+1}-1) \mathcal{A}_4 (s) + t_{n+1} \mathcal{A}_5 (s; \gamma_n) \right] \, .
\end{align}
This means that
\begin{align*}
& (1-t_{n+1})^2 \mathcal{A}_2^2 (s) + t_{n+1}^2 \mathcal{A}_3^2 (s; \gamma_n, L) - 2 t_{n+1} (1-t_{n+1}) \mathcal{A}_2 (s) \mathcal{A}_3 (s; \gamma_n, L) + (\lambda_n + \omega_n)^2 t_{n+1}^2 \\
&- 2(\lambda_n + \omega_n) t_{n+1} \left[ (1-t_{n+1}) \mathcal{A}_2 (s) - t_{n+1} \mathcal{A}_3 (s; \gamma_n, L) \right] \leq 2 s t_{n+1} \left[ (t_{n+1}-1) \mathcal{A}_4 (s) + t_{n+1} \mathcal{A}_5 (s; \gamma_n) \right ], .
\end{align*}
Expanding the terms, it follows that
\begin{align*}
& t_{n+1}^2 \mathcal{A}_2^2 (s) + \mathcal{A}_2^2 (s) - 2 t_{n+1} \mathcal{A}_2^2 (s) + t_{n+1}^2 \mathcal{A}_3^2 (s; \gamma_n, L) - 2 t_{n+1} \mathcal{A}_2 (s) \mathcal{A}_3 (s; \gamma_n, L) + 2 t_{n+1}^2 \mathcal{A}_2 (s) \mathcal{A}_3 (s; \gamma_n, L) \\
&+ (\lambda_n + \omega_n)^2 t_{n+1}^2 -2 \mathcal{A}_2 (s) (\lambda_n + \omega_n) t_{n+1} + 2 \mathcal{A}_2 (s) (\lambda_n + \omega_n) t_{n+1}^2 + 2(\lambda_n + \omega_n) t_{n+1}^2 \mathcal{A}_3 (s; \gamma_n, L) \\
&\leq 2 s \mathcal{A}_4 (s) t_{n+1}^2 - 2 s \mathcal{A}_4 (s) t_{n+1} + 2 s t_{n+1}^2 \mathcal{A}_5 (s; \gamma_n) \, . 
\end{align*}
Grouping by $t_{n+1}$ and $t_{n+1}^2$ implies that
\begin{align*}
& t_{n+1}^2 [ \mathcal{A}_2^2 (s) + \mathcal{A}_3^2 (s; \gamma_n, L) + 2 \mathcal{A}_2 (s) \mathcal{A}_3 (s; \gamma_n, L) +(\lambda_n + \omega_n)^2 + 2 \mathcal{A}_2(s) (\lambda_n + \omega_n) + 2(\lambda_n + \omega_n) \mathcal{A}_3 (s; \gamma_n, L) \\
&- 2 s \mathcal{A}_4 (s) - 2s \mathcal{A}_5 (s; \gamma_n) ] + t_{n+1} [ - 2 \mathcal{A}_2^2 (s) - 2 \mathcal{A}_2 (s) \mathcal{A}_3 (s; \gamma_n, L) - 2 \mathcal{A}_2 (s) (\lambda_n + \omega_n) + 2 s \mathcal{A}_4 (s) ] + \mathcal{A}_2^2 (s) \leq 0 \, .
\end{align*}
Finally, one has that
\begin{align}\label{GroupingTerms}
& t_{n+1}^2 \left[ ( \mathcal{A}_2 (s) + \mathcal{A}_3 (s; \gamma_n, L) )^2 - 2s (\mathcal{A}_4(s) + \mathcal{A}_5(s; \gamma_n)) + (\lambda_n+\omega_n)^2 + 2(\lambda_n+\omega_n) (\mathcal{A}_2(s) + \mathcal{A}_3 (s; \gamma_n, L)) \right] + \\
& t_{n+1} \left[ - 2 \mathcal{A}_2^2 (s) - 2 \mathcal{A}_2 (s) \mathcal{A}_3 (s; \gamma_n, L) - 2(\lambda_n+\omega_n) \mathcal{A}_2 (s) + 2 s \mathcal{A}_4 (s) \right] + \mathcal{A}_2^2 (s) \leq 0 \, . \nonumber
\end{align}
We rewrite \ref{GroupingTerms} as
\begin{align}\label{eq2:Quadratic}
G_n t_{n+1}^2 - H_n t_{n+1} - I_n \geq 0 \, ,
\end{align}
where the coefficients are
\begin{align*}
\begin{cases}
G_n &:= - (\mathcal{A}_2 (s) + \mathcal{A}_3 (s; \gamma_n, L))^2 + 2s (\mathcal{A}_4 (s) + \mathcal{A}_5 (s; \gamma_n)) - (\lambda_n + \omega_n)^2 - 2(\lambda_n + \omega_n) (\mathcal{A}_2 (s) + \mathcal{A}_3 (s; \gamma_n, L)) \\
H_n &:= -2 \mathcal{A}_2^2 (s) - 2 \mathcal{A}_2 (s) \mathcal{A}_3 (s; \gamma_n, L) - 2 (\lambda_n + \omega_n) \mathcal{A}_2 (s) + 2s \mathcal{A}_4 (s) \\
I_n &:= \mathcal{A}_2^2 (s)
\end{cases}
\end{align*}
Since $I_n = s^2 > 0$, in order to use Lemma \ref{L18} and taking into account Remark \ref{R19}, we only need to show that for $n$ large enough $G_n > 0$, hence $H_n^2 + 4 G_n I_n \geq 0$. Now, this means that \ref{eq2:Quadratic} is equivalent to the fact that $(-G_n) t_{n+1}^2 + H_n t_{n+1} + I_n \leq 0$. Based upon  \ref{GroupingTerms}, we compute $-G_n$ as follows:
\begin{align*}
-G_n &= \left[ s + \gamma_n (1-Ls) \right]^2 - s^2 + \gamma_n^2 + (\lambda_n + \omega_n)^2 + 2 (\lambda_n + \omega_n)(-s-\gamma_n(1-Ls)) \\
&= \left[ s + \gamma_n (1-Ls) \right]^2 + (\gamma_n^2 - s^2) + (\lambda_n+\omega_n)^2 - 2 (\lambda_n + \omega_n)(s+\gamma_n(1-Ls)) \\
&= \left[ s + \gamma_n (1-Ls) \right]^2 + (\gamma_n^2 - s^2) + (\lambda_n+\omega_n)^2 - (\lambda_n + \omega_n)(s+\gamma_n(1-Ls)) - (\lambda_n + \omega_n)(s+\gamma_n(1-Ls)) \\
&= \left[ s + \gamma_n(1-Ls) \right] \cdot \left[ s + \gamma_n(1-Ls) - (\lambda_n+\omega_n) \right] + (\gamma_n^2 - s^2) + (\lambda_n+\omega_n) \left[ (\lambda_n+\omega_n) - (s+\gamma_n(1-Ls)) \right] \\
&= (\gamma_n^2 - s^2) + \left[ s+\gamma_n(1-Ls) - (\lambda_n+\omega_n) \right] \cdot \left[ s + \gamma_n(1-Ls) \right] - (\lambda_n+\omega_n) \left[ s + \gamma_n (1-Ls) - (\lambda_n+\omega_n) \right] \\
&= (\gamma_n^2 - s^2) + \left[ (s+\gamma_n(1-Ls)) - (\lambda_n+\omega_n) \right] \cdot \left[ (s+\gamma_n(1-Ls)) - (\lambda_n + \omega_n) \right] \\
&= (\gamma_n^2 - s^2) + \left[ (s+\gamma_n(1-Ls)) - (\lambda_n+\omega_n) \right]^2 \, .
\end{align*}
Considering the above computations for $G_n$ and the form of $\lambda_n + \omega_n$ from assumption (ii), we obtain, after some basic evaluations, exactly assumption (i), where it is satisfied for each $n > max \lbrace \alpha - 1 , N^\prime \rbrace$. \\
\textbf{( Step III )} In the last step of the present proof we recollect the results from above and apply Lemma \ref{L17} and Lemma \ref{L18} and actually show the $O \left( \dfrac{1}{n^2} \right)$ asymptotic rate of converge. At the same time, we actually determine the index $N \in \mathbb{N}$, such that for $n \geq N$, the rate of convergence in the objective function is obtained. \\
From now on, we consider only when $n > max \lbrace N_1, N_2, N^\prime \rbrace$, where $N_1 := \alpha - 1$, $N_2 := (\alpha-1) \dfrac{H_n + \sqrt{H_n^2 + 4 G_n I_n}}{2 G_n}$ and $N^\prime$ the index from the hypothesis of the theorem. Now, we know from the theorem's hypothesis that $(\gamma_n)_{n \in \mathbb{N}}$ and $(\lambda_n)_{n \in \mathbb{N}}$ are convergent. By using assumption (ii), it follows that $(\omega_n)_{n \in \mathbb{N}}$ is also convergent. This implies that also $(H_n)_{n \in \mathbb{N}}$ and $(G_n)_{n \in \mathbb{N}}$ are convergent. Since these two sequences are convergent, they are bounded. Hence, the sequence $\left( (\alpha-1) \dfrac{H_n + \sqrt{H_n^2 + 4 G_n I_n}}{2 G_n} \right)_{n \in \mathbb{N}}$ is also bounded. So, let $M > 0$ such that  $(\alpha-1) \dfrac{H_n + \sqrt{H_n^2 + 4 G_n I_n}}{2 G_n} < M$. Then, we can take $n \geq M > N_2$ (in this sense $n > N_2$ is well defined). \\
Consider \ref{eq2:Quadratic}, with the nonconstant coefficients $G_n$, $H_n$ and $I_n$. We apply Lemma \ref{L18} with $x = t_{n+1}$, $a(x)=a(t_{n+1}) = G_n$, $b(x)=b(t_{n+1}) = - H_n$ and $c(x)=c(t_{n+1}) = - I_n$. Since $n > N_2$, then we consider find an infinity of points $x = t_{n+1} \geq \dfrac{H_n + \sqrt{H_n^2 + 4 G_n I_n}}{2 G_n}$ that obviously satisfy $c(t_{n+1}) = - s^2 < 0$, since $s > 0$. Also, the condition $b^2(t_{n+1}) - 4 a(t_{n+1}) c(t_{n+1}) \geq 0$ is equivalent to $H_n^2 + 4 G_n I_n \geq 0$. Because $n > max \lbrace N_1, N_2, N^\prime \rbrace$, assumption (i) is in fact $G_n > 0 $, i.e. $a(t_{n+1}) > 0$ (this stands true for $n \geq \max \lbrace \alpha - 1 , N^\prime \rbrace$). This means that indeed we have $H_n^2 + 4G_n I_n \geq 0$, for each $n > max \lbrace N_1, N_2, N^\prime \rbrace$. \\
After these computations, now we can consider the quadratic inequality with nonconstant coefficients $F_n \geq 0$, with $F_n$ given in \ref{eq1:Quadratic} and where $n > \max \lbrace N_1, N_2, N^\prime \rbrace$. We apply Lemma \ref{L17} with $x = \| \nabla f(y_n) \|$, $a(x) = a(\| \nabla f(y_n) \|) = \dfrac{s}{2} t_{n+1}$, $b(x) = b(\| \nabla f(y_n) \|) = \left[ (1-t_{n+1}) \mathcal{A}_2 (s) - t_{n+1} \mathcal{A}_3 (s; \gamma_n, L) - (\lambda_n + \omega_n) t_{n+1} \right] \| \nabla f(x_n) \|$ and $c(x) = c(\| \nabla f(y_n) \|) = \left[ (t_{n+1} - 1) \mathcal{A}_4 (s) + t_{n+1} \mathcal{A}_5 (s; \gamma_n) \right] \| \nabla f(x_n) \|^2$. So, in this case we consider $x = t_{n+1}$ with $n > \max \lbrace N_1, N_2, N^\prime \rbrace$. The condition $a(\| \nabla f(y_n) \|) > 0$ is obviously satisfied for $n > \max \lbrace N_1, N_2, N^\prime \rbrace$. Furthermore, the condition is equivalent to \ref{CondQuadratic}. Moreover, this is in fact \ref{eq2:Quadratic}, that we have shown that is valid for $n > max \lbrace N_1, N_2, N^\prime \rbrace$. We conclude that $F_n \geq 0$ for every $\nabla f(y_n)$ with $n > max \lbrace N_1, N_2, N^\prime \rbrace$. \\
By \ref{EnergyDecrease}, we find that $\mathcal{E}_{n+1} - \mathcal{E}_n \leq 0$, so $\mathcal{E}_n \leq \mathcal{E}_N$, where $N := \max \lbrace N_1, N_2, N^\prime \rbrace$. By the fact that $t_n = \dfrac{n-1}{\alpha-1}$ and utilizing the definition \ref{Energy}, we have that $t_n^2 \left[ f(x_n) - f(x^\ast) \right] \leq \mathcal{E}_N$. Finally, we obtain that $f(x_n) - f(x^\ast) \leq \dfrac{\mathcal{E}_N (\alpha-1)^2}{(n-1)^2}$, for $n > N$ and the proof is over. 
\end{Proof}

\begin{remark}\label{R21}
Take the case when the sequence $(\lambda_n)_{n \in \mathbb{N}}$ is constant, with $\lambda_n = \beta \sqrt{s}$. 
Then, the value of $z_n$ given in \ref{Energy} is $z_n = (x_{n-1} - x^\ast) + t_n(x_n - x_{n-1}) + \beta \sqrt{s} t_{n+1} \nabla f(x_{n-1})$. This discrete value differs significantly from that given in \cite{Attouch_Chbani}, where $z_n = (x_{n-1} - x^\ast) + t_n \left[ (x_n - x_{n-1}) + \beta \sqrt{s} \nabla f(x_{n-1}) \right]$, by the fact that we have used $\nabla f(x_n)$ and not $\nabla f(x_{n-1})$ in \ref{LT-S-IGAHD} (see Remark \ref{R5} and Subsection \ref{S22}). Last of all, we point out that even though \ref{Energy} is similar to that used in \cite{Attouch_Chbani}, there are subtle differences consisting in the fact that our choice of the discrete functional is optimal, since in $z_{n+1}$ we got rid of $\nabla f(x_{n-1})$. This means that every other types of discrete Lyapunov functionals that can be constructed, contain $\nabla f(x_{n-1})$, that can't be eliminated by the fact that they have a nonconstant term in front of them.
\end{remark}

\begin{remark}\label{R22}
We observe that assumption (ii) of Theorem \ref{T20} reduces to the fact that $\gamma^2_n < s^2$, namely $\gamma_n \in (-s, s)$ for each $n > max \lbrace \alpha-1, N^\prime \rbrace$. Furthermore, following the last step of the previous proof and considering Lemma \ref{L17} and Remark \ref{R19}, the following sequence of equivalences stand true:
\begin{align*}
\gamma_n^2 < s^2 \left( 1 - \dfrac{\alpha-1}{n} \right) & \Longleftrightarrow s - \dfrac{\gamma_n^2}{s} > \dfrac{s(\alpha-1)}{n} \\
&\Longleftrightarrow \dfrac{n}{\alpha-1} \left( \dfrac{s}{2} - \dfrac{\gamma_n^2}{2s} \right) > \dfrac{s}{2} 
\Longleftrightarrow t_{n+1} \left[ \mathcal{A}_4 (s) + \mathcal{A}_5 (s;\gamma_n) \right] > \mathcal{A}_4 (s) \, ,
\end{align*}
So, in fact we have a stronger bound for $\gamma_n$, namely $ - s^2 \left( 1 - \dfrac{\alpha-1}{n} \right) < \gamma_n <  s^2 \left( 1 - \dfrac{\alpha-1}{n} \right) $, and this is valid for all $n > max \lbrace \alpha - 1, N^\prime \rbrace$.
\end{remark}

\begin{remark}\label{R23}
In Theorem \ref{T20} we have considered the strict inequality belonging to assumption (ii) to be valid for each $n > \max \lbrace \alpha - 1, N^\prime \rbrace$. But, if we look closer at the third step of the proof from the same theorem, we observe that we could have taken $n > \max \lbrace \alpha - 1, N_2, N^\prime \rbrace$, where $N_2$ was defined through $H_n$, $G_n$ and $I_n$. Hence, in this way, we obtain a more relaxed bound on the restriction on the index $n$. \\
On the other hand, we observe that the assumptions that the sequences $(\gamma_n)_{n \in \mathbb{N}}$ and $(\lambda_n)_{n \in \mathbb{N}}$ are bounded, so also $(\omega_n)_{n \in \mathbb{N}}$, can be linked to our reasoning from Subsection \ref{S24}, where we have supposed the convergence of the sequences in our heuristic analysis.
\end{remark}

We end this section by giving some examples of algorithms of tpye \ref{LT-S-IGAHD} that satisfy the assumptions from Theorem \ref{T20} and we also show that indeed exists an index $N^\prime$ for each of them. In our first example, we consider a general class of Hessian driven optimization algorithms with asymptotic vanishing damping, that include three parameters $a$, $\mu$ and $b$ that will play a crucial role in the numerical behavior, in the sense that, from a scientific computing point of view, we have a much larger freedom and flexibility into the hypertuning of the parameters.
\begin{example}\label{E24}
We consider the choice of the coefficients as follows: \\
$\gamma_n = s \sqrt{\dfrac{\alpha-1}{n+a}}$, $\lambda_{n+1} = s \dfrac{n}{n+1} + \mu \dfrac{n}{(n+1)(n+b)}$ and $\omega_n = s \sqrt{\dfrac{\alpha-1}{n+a}} + \dfrac{s}{n} + \mu \left[ \dfrac{1}{n+b} - \dfrac{n-1}{n(n+b-1)} \right]$, where $a \geq 0$, $b \geq 0$ and $\mu \geq 0$. Now, we validate our algorithm as follows: \\
With the choices of the sequences $(\gamma_n)_{n \in \mathbb{N}}$ and $(\lambda_n)_{n \in \mathbb{N}}$, we obtain the following chain of equivalences:
\begin{align*}
\omega_{n} &= s\sqrt{\dfrac{\alpha-1}{n+a}} + \dfrac{s}{n}+\mu\left[ \dfrac{1}{n+b}-\dfrac{n-1}{n(n+b-1)}\right] \\
\omega_{n} &= s\sqrt{\dfrac{\alpha-1}{n+a}} - s\left(\dfrac{n-1}{n}-1\right)+\mu\left[ \dfrac{1}{n+b}-\dfrac{n-1}{n(n+b-1)}\right] \\
\omega_{n} &= s\sqrt{\dfrac{\alpha-1}{n+a}} - s\dfrac{n-1}{n}-\mu \dfrac{n-1}{n(n+b-1)}+s+\dfrac{\mu}{n+b}\\
\omega_{n} &= \gamma_{n}-\lambda_{n}+\dfrac{n+1}{n}\lambda_{n+1}
\end{align*}
so assumption (ii) of Theorem \ref{T20} is satisfied for each $n \geq 1$. \\
We must show assumption (i) of Theorem \ref{T20}. For this, we have the following equivalences:
\begin{align*}
\left[s(1-\gamma_{n}L)-\dfrac{n+1}{n}\lambda_{n+1}\right]^2<s^2-\gamma_{n}^2 & \Longleftrightarrow \left[s\gamma_{n}L+\dfrac{\mu}{n+b}\right]^2<s^2-\gamma_{n}^2 \\
& \Longleftrightarrow s^2\gamma_{n}^2L^2+ 2s\mu L \gamma_{n}\dfrac{1}{n+b}+\dfrac{\mu^2}{(n+b)^2}< s^2-\gamma_{n}^2 \\
& \Longleftrightarrow s^4\dfrac{\alpha-1}{n+a}L^2+2s^2\mu L\dfrac{\sqrt{\alpha-1}}{(n+b)\sqrt{n+a}}+\dfrac{\mu^2}{(n+b)^2} < s^2-s^2\dfrac{\alpha-1}{n+a}
\end{align*}
Multiplying by $1/s^2$, we get that
\begin{align}\label{eq1:Example1}
s^2\dfrac{\alpha-1}{n+a}L^2+2\mu L\dfrac{\sqrt{\alpha-1}}{(n+b)\sqrt{n+a}}+\left(\dfrac{\mu}{s}\right)^2\dfrac{1}{(n+b)^2}+\dfrac{\alpha-1}{n+a}<1 \, .
\end{align}
Let $\sqrt{n+a}<n+b$. So $n+a-n^2-2bn-b^2<0$, which means that $n^2 + (2b-1)n + (b^2-a) > 0$, with $\Delta=(2b-1)^2-4(b^2-a)$. We consider the following cases :
\begin{itemize}
\item If $\Delta \geq 0$, namely $b-a \leq \dfrac{1}{4}$ then we can take $n > \dfrac{1-2b+\sqrt{4(a-b)+1}}{2} $, so the inequality is satisfied.
\item If $\Delta < 0$, i.e. $b - a > \dfrac{1}{4}$, then we can take $n \geq 1$.
\end{itemize}
Since
$$ s^2\dfrac{\alpha-1}{n+a}L^2 + 2\mu L\dfrac{\sqrt{\alpha-1}}{(n+b) \sqrt{n+a}}+\left(\dfrac{\mu}{s}\right)^2 \dfrac{1}{(n+b)^2} + \dfrac{\alpha-1}{n+a} < s^2\dfrac{\alpha-1}{n+a}L^2+2\mu L\dfrac{\sqrt{\alpha-1}}{n+a}+\left(\dfrac{\mu}{s}\right)^2\dfrac{1}{n+a}+\dfrac{\alpha-1}{n+a} \, , $$
and taking into account \ref{eq1:Example1}, then it is enough to show that:
\begin{align*}
s^2\dfrac{\alpha-1}{n+a}L^2+2\mu L\dfrac{\sqrt{\alpha-1}}{n+a}+\left(\dfrac{\mu}{s}\right)^2\dfrac{1}{n+a}+\dfrac{\alpha-1}{n+a}<1 &\Longleftrightarrow \dfrac{1}{n+a}\left[(\alpha-1)(s^2L^2+1)+2\mu L\sqrt{\alpha-1}+\left(\dfrac{\mu}{s}\right)^2\right] <1 \\
&\Longleftrightarrow n > (\alpha-1)(s^2L^2+1)+2\mu L\sqrt{\alpha-1}+\left(\dfrac{\mu}{s}\right)^2-a  \, .
\end{align*}
Finally, if $b - a > \dfrac{1}{4}$ then we can consider $$N^\prime := (\alpha-1)(s^2L^2+1)+2\mu L\sqrt{\alpha-1}+\left(\dfrac{\mu}{s}\right)^2-a \, ,$$ and if $b - a \leq \dfrac{1}{4}$, then we can take $$N^\prime := \max \Big\lbrace \dfrac{1-2b+\sqrt{4(a-b)+1}}{2}, (\alpha-1)(s^2L^2+1)+2\mu L\sqrt{\alpha-1}+\left(\dfrac{\mu}{s}\right)^2-a \Big\rbrace.$$
\end{example}

The second example is also based upon three coefficient $\beta$, $b$ and $\mu$. It is worth noticing that for the particular case when $\mu = 0$ we obtain the \ref{IGAHD} algorithm that differs from that given in \cite{Attouch_Chbani} (see Remark \ref{R5} and Subsection \ref{S22}).
\begin{example}\label{E25}
We consider the following choice of the coefficients: \\
$\gamma_n = 0$, $\lambda_{n+1} = \beta \sqrt{s} + \dfrac{\mu}{n+b}$ and $\omega_n = \dfrac{\beta \sqrt{s}}{n} + \mu \left[ \dfrac{n+1}{n(n+b)} - \dfrac{1}{n+b-1} \right]$, where $0 < \beta < 2 \sqrt{s}$, $b > 0$ and $\mu \geq 0$. Now, we validate our algorithm as follows: \\
With the choices of the sequences $(\gamma_n)_{n \in \mathbb{N}}$ and $(\lambda_n)_{n \in \mathbb{N}}$, we obtain the following chain of equivalences:
\begin{align*}
\omega_{n} &= \dfrac{\beta\sqrt{s}}{n}+\mu\left[\dfrac{n+1}{n(n+b)}-\dfrac{1}{n+b-1}\right] \\
\omega_{n} &= \dfrac{n+1}{n}\beta\sqrt{s}+\mu\dfrac{n+1}{n(n+b)}-\beta\sqrt{s}-\dfrac{\mu}{n+b-1}\\
\omega_{n} &= \dfrac{n+1}{n}\lambda_{n+1}-\lambda_{n}\\
\omega_{n} &= \gamma_{n}-\lambda_{n}+\dfrac{n+1}{n}\lambda_{n+1}
\end{align*}
We must show assumption (i) of Theorem \ref{T20}. For this, we have the following equivalences:
\begin{align*}
\left[s(1-\gamma_{n}L)-\dfrac{n+1}{n}\lambda_{n+1}\right]^2<s^2-\gamma_{n}^2 &\Longleftrightarrow \left[s-\dfrac{n+1}{n}\lambda_{n+1}\right]^2<s^2 \\
&\Longleftrightarrow \left[s-\dfrac{n+1}{n}\beta\sqrt{s}-\mu\dfrac{n+1}{n(n+b)}\right]^2<s^2 \\
&\Longleftrightarrow -\left[\dfrac{n+1}{n}\beta\sqrt{s}+\mu\dfrac{n+1}{n(n+b)}\right]\left[2s-\dfrac{n+1}{n}\beta\sqrt{s}-\mu\dfrac{n+1}{n(n+b)}\right]< 0
\end{align*}
By the fact that $\dfrac{n+1}{n}\beta\sqrt{s}+\mu\dfrac{n+1}{n(n+b)}$ is positive and since $\beta > 0$, this term can't be $0$. So, the inequality from above is satisfied if and only if $ \dfrac{n+1}{n} \beta\sqrt{s} + \mu \dfrac{n+1}{n(n+b)}< 2 s $.
Multiplying by $\dfrac{n}{\sqrt{s}(n+1)}$, it follows that $\beta +\dfrac{\mu}{\sqrt{s}(n+b)}<2\sqrt{s}\dfrac{n}{n+1}$. Hence 
\begin{align}\label{eq1:Example2}
\beta +\dfrac{\mu}{\sqrt{s}(n+b)}<2\sqrt{s}\dfrac{n}{n+1} \, .
\end{align}
Since $\beta < 2 \sqrt{s}$, there exists $ \tau>0$ such that $2 \sqrt{s} - \beta > \tau$. We consider $\epsilon = \dfrac{\tau}{\dfrac{\mu}{\sqrt{s}}+2\sqrt{s}} > 0$. From the Archimedean property, we know that for $\epsilon>0$ there exists $N \in \mathbb{N}$ such that $\epsilon>\dfrac{1}{N}$, thus the strict inequality $2\sqrt{s}-\beta>\left(\dfrac{\mu}{\sqrt{s}}+2\sqrt{s}\right)\epsilon$ implies that $ 2\sqrt{s}-\beta>\left(\dfrac{\mu}{\sqrt{s}}+2\sqrt{s}\right)\dfrac{1}{N}$. For $n+1 \geq N$ and $n+b \geq N$ which are equivalent to $\dfrac{1}{N} \geq \dfrac{1}{n+1}$ and $\dfrac{1}{N} \geq \dfrac{1}{n+b}$, we obtain that $$\left(\dfrac{\mu}{\sqrt{s}}+2\sqrt{s}\right)\dfrac{1}{N} \geq \dfrac{\mu}{\sqrt{s}} \dfrac{1}{n+b} + 2\sqrt{s} \dfrac{1}{n+1} \Rightarrow2 \sqrt{s}-\beta > \dfrac{\mu}{\sqrt{s}} \dfrac{1}{n+b} + 2\sqrt{s} \dfrac{1}{n+1}.$$
Second to last, we can take $N^\prime := max \lbrace N - 1, N - b \rbrace$ and we obtain that \ref{eq1:Example2} is satisfied for every $n \geq N^\prime$. Finally, one can find in an explicit way the index $N^\prime$. In order to do this, observe that \ref{eq1:Example2} leads to a quadratic inequality, namely: $$ (2 \sqrt{s} - \beta) n^2 + \left[ 2 \sqrt{s} b - \beta(b+1) - \dfrac{\mu}{\sqrt{s}} \right] n - \left( \beta b + \dfrac{\mu}{\sqrt{s}} \right) > 0.$$ In this case, one can see that $N^\prime = \dfrac{\left[ \beta(b+1) + \dfrac{\mu}{\sqrt{s}} - 2 \sqrt{s} b \right] + \sqrt{\Delta}}{2(2 \sqrt{s} - \beta)}$, where $\Delta := \left( 2 b \sqrt{s} - \beta (b+1) - \dfrac{\mu}{\sqrt{s}} \right)^2 + 4 (2 \sqrt{s} - \beta) \left( \beta b + \dfrac{\mu}{\sqrt{s}} \right) \geq 0 \, .$
\end{example}

Our last example consists in choosing $\gamma_n$ to take negative values, but $\lambda_n$ to have only positive values. Thus, $\omega_n$ is determined through assumption (ii) of Theorem \ref{T20}.
\begin{example}\label{E26}
We consider the choice of the coefficients as follows: \\
$\gamma_n = - \dfrac{s}{n+a}$, $\lambda_{n+1} = s \dfrac{n}{n+1} + \mu \dfrac{n}{(n+1)(n+b)}$ and $\omega_n = - \dfrac{s}{n+a} + \dfrac{s}{n} + \mu \left[ \dfrac{1}{n+b} - \dfrac{n-1}{n(n+b-1)} \right]$, where $a \geq 0$, $b \geq 0$ and $\mu \geq 0$. Now, we validate our algorithm as follows: \\
With the choices of the sequences $(\gamma_n)_{n \in \mathbb{N}}$ and $(\lambda_n)_{n \in \mathbb{N}}$, we obtain the following chain of equivalences:
\begin{align*}
\omega_{n} &= - \dfrac{s}{n+a} + \dfrac{s}{n}+\mu\left[ \dfrac{1}{n+b}-\dfrac{n-1}{n(n+b-1)}\right] \\
\omega_{n} &= - \dfrac{s}{n+a} - s\left(\dfrac{n-1}{n}-1\right)+\mu\left[ \dfrac{1}{n+b}-\dfrac{n-1}{n(n+b-1)}\right] \\
\omega_{n} &= - \dfrac{s}{n+a} - s\dfrac{n-1}{n}-\mu \dfrac{n-1}{n(n+b-1)}+s+\dfrac{\mu}{n+b}\\
\omega_{n} &= \gamma_{n}-\lambda_{n}+\dfrac{n+1}{n}\lambda_{n+1}
\end{align*}
so assumption (ii) of Theorem \ref{T20} is satisfied for each $n \geq 1$. \\
Now, regarding assumption (i), we consider the following equivalences:
\begin{align*}
\left[s(1-\gamma_{n}L)-\dfrac{n+1}{n}\lambda_{n+1}\right]^2<s^2-\gamma_{n}^2 & \Longleftrightarrow \left[s\gamma_{n}L+\dfrac{\mu}{n+b}\right]^2<s^2-\gamma_{n}^2 \\
& \Longleftrightarrow s^2\gamma_{n}^2L^2+ 2s\mu L \gamma_{n}\dfrac{1}{n+b}+\dfrac{\mu^2}{(n+b)^2}< s^2-\gamma_{n}^2 \\
& \Longleftrightarrow s^4L^2 \dfrac{1}{(n+a)^2} - 2s^2\mu L\dfrac{1}{(n+a)(n+b)}+\dfrac{\mu^2}{(n+b)^2} < s^2-s^2\dfrac{1}{(n+a)^2}
\end{align*}
Multiplying by $1/s^2$, we get that
\begin{align}\label{eq1:Example3}
\dfrac{(sL)^2+1}{(n+a)^2} - \dfrac{2 \mu L}{(n+a)(n+b)} + \left( \dfrac{\mu}{s} \right)^2 \dfrac{1}{(n+b)^2} < 1 \, .
\end{align}
Now, we analyze two cases: \\
\textbf{1)} Suppose that $a \leq b$. Then, $\dfrac{1}{n+b} \leq \dfrac{1}{n+a}$. We consider $n > N^\prime := \sqrt{(sL)^2 + 1 + \left( \dfrac{\mu}{s} \right)^2} - a$. This leads to the fact that $(n+a)^2 >(sL)^2 + 1 + \left( \dfrac{\mu}{s} \right)^2$. Thus, $\dfrac{(sL)^2 + 1}{(n+a)^2} + \left( \dfrac{\mu}{s} \right)^2 \dfrac{1}{(n+a)^2} < 1$. At the same time, we have the following inequalities:
\begin{align*}
\dfrac{(sL)^2 + 1}{(n+a)^2} - \dfrac{2 \mu L}{(n+a)(n+b)} + \left( \dfrac{\mu}{s} \right)^2 \dfrac{1}{(n+b)^2} & \leq \dfrac{(sL)^2+1}{(n+a)^2} + \left( \dfrac{\mu}{s} \right)^2 \dfrac{1}{(n+b)^2} \leq \dfrac{(sL)^2+1}{(n+a)^2} + \left( \dfrac{\mu}{s} \right)^2 \dfrac{1}{(n+a)^2} < 1 \, . 
\end{align*}
\textbf{2)} Suppose that $b \leq a$. Then, $\dfrac{1}{n+a} \leq \dfrac{1}{n+b}$. We consider $n > N^\prime := \sqrt{(sL)^2 + 1 + \left( \dfrac{\mu}{s} \right)^2} - b$. This leads to the fact that $(n+b)^2 >(sL)^2 + 1 + \left( \dfrac{\mu}{s} \right)^2$. Thus, $\dfrac{(sL)^2 + 1}{(n+b)^2} + \left( \dfrac{\mu}{s} \right)^2 \dfrac{1}{(n+b)^2} < 1$. At the same time, we have the following inequalities:
\begin{align*}
\dfrac{(sL)^2 + 1}{(n+a)^2} - \dfrac{2 \mu L}{(n+a)(n+b)} + \left( \dfrac{\mu}{s} \right)^2 \dfrac{1}{(n+b)^2} & \leq \dfrac{(sL)^2+1}{(n+a)^2} + \left( \dfrac{\mu}{s} \right)^2 \dfrac{1}{(n+b)^2} \leq \dfrac{(sL)^2+1}{(n+b)^2} + \left( \dfrac{\mu}{s} \right)^2 \dfrac{1}{(n+b)^2} < 1 \, . 
\end{align*}
Finally, we end this example emphasizing that we can take $N^\prime = \sqrt{3 + \left( \dfrac{\mu}{s} \right)^2} - a$ if $a \leq b$ and $N^\prime = \sqrt{3 + \left( \dfrac{\mu}{s} \right)^2} - b$ if $b \leq a$, and in this way we obtain an index that does not depend on the Lipschitz constant $L$.
\end{example}

\section{Numerical experiments involving convex functions}\label{S4}

The role of the present section is to validate our theoretical results through some computational simulations, involving some convex test functions. Furthermore, we consider the particular cases of \ref{LT-S-IGAHD} given in Example \ref{E24}, Example \ref{E25} and Example \ref{E26} and we study the role of the parameters for each of them. Also, our codes for the optimization algorithms are implemented in \texttt{Matlab} and the analysis of the convexity of the chosen functions is based upon the a fundamental theory that can be found in \cite{Beck}.

\subsection{Algorithmic description and objective functions}\label{S41}

We consider the algorithms of the form \ref{LT-S-IGAHD} (with the stepsize $s = h^2$), starting from $n \geq 1$. We initialize the starting point $x_0 \in \mathbb{R}^2$ for functions $f : \mathbb{R}^2 \to \mathbb{R}$. Furthermore, even though in many simulations $x_1 \in \mathbb{R}^2$ is considered arbitrary, for consistency we choose $x_1 = x_0 - s \nabla f(x_0)$,  namely the first iteration at $n=1$ is given by a step of gradient descent. Throughout the present section we count the total number of iterations starting from $n=0$ in order to encompass the initial value $x_0$ (even though our algorithms start from $n=1$). Furthermore, since for $x_2$ we need only the values $x_0$, $x_1$ and $y_1$, for brevity we set $y_0 = x_0$. Using tolerance values of low order $\epsilon = 10^{-10}$, we employ the stopping criteria as $\| f(x_n) - f(x_{n-1}) \|_2 \leq \epsilon$ or as $\| f(x_n) - f(x^\ast) \|_2 \leq \epsilon$ in the case when the objective function admits a unique minimum point $x^\ast \in \mathbb{R}^2$. For each of the functions that we consider we compute the Lipschitz constant in order to present the behavior of the discretizations with respect to the parameters. Furthermore, we consider the stepsize $s < \dfrac{1}{L}$ by the fact that in the optimal case when $s = \dfrac{1}{L}$ no interpretation can be given since the algorithms stop just after a few iterations. We consider the following convex functions: \\
\textbf{1)} \underline{\textit{The case of convex, but not strongly convex function}} (see \cite{Attouch_Chbani}), namely $f : \mathbb{R}^2 \to \mathbb{R}$, defined as
\begin{align}\label{F1}\tag{F1}
f(x_1, x_2) = (x_1 + x_2)^2 \, , \, \forall (x_1, x_2) \in \mathbb{R}^2 \, ,
\end{align}
It is easy to observe that $\nabla f(x_1, x_2) = \begin{pmatrix}
2(x_1 + x_2) \\ 
2(x_1 + x_2)
\end{pmatrix}$. Furthermore, the Hessian of $f$ is $\nabla^2 f(x_1, x_2) = 
\begin{pmatrix}
2 & 2 \\
2 & 2
\end{pmatrix}$.
Taking $z = \begin{pmatrix}
z_1 \\
z_2
\end{pmatrix}$,
we infer that $z^T \nabla^2 f(x_1, x_2) z = 2 (z_1 + z_2)^2 \geq 0$, so the Hessian matrix $\nabla^2 f \succeq 0$. Using the second order characterization of convexity (see Theorem 7.12 from \cite{Beck}), we find that the function is convex. Furthermore, computing the solution of the equation $\nabla f(x_1, x_2) = \begin{pmatrix}
0 \\
0
\end{pmatrix}$ and using the fact from above that the objective function is convex, by using Theorem 7.9 of \cite{Beck}, we obtain that  the set of stationary points is given by $\argmin f = \lbrace (x_1, x_2) \in \mathbb{R}^2 \, / \, x_1 = - x_2 \rbrace$. By the fact that we use the $l_2$ norm on the gradient, we aim at evaluating the subordinate matrix norm. By the fact that $\| \nabla^2 f(x_1, x_2) \|_2$ consists of the largest square root of the absolute value of the matrix eigenvalues, we find that $\| \nabla^2 f(x_1, x_2) \| = 4$, so that the Lipschitz constant of the gradient $\nabla f$ is also $L = 4$. On the other hand, since the $\argmin$ set contains an infinity of elements, we use the stopping criteria $\| f(x_n) - f(x_{n-1}) \|_2 \leq \epsilon$, for a given tolerance $\epsilon > 0$ sufficiently small. Moreover, it is also natural to adapt the \ref{Energy} function for the case when we do not actually now to what minimizer our algorithm converges. For this, we consider 
\begin{align}
& \mathcal{E}_n := t_{n}^{2} \left( f(x_n) - f(x_M) \right) + \dfrac{1}{2s} \| \tilde{z}_n \|^2 \\
& \tilde{z}_n := (x_{n-1} - x_M) + t_{n} (x_n - x_{n-1}) + \lambda_n t_{n+1} \nabla f(x_{n-1}) \, , \nonumber
\end{align}
where $M$ is determined by the stopping condition, i.e. when we reach the tolerance factor $\epsilon$, we consider the last iteration value $x_n = x_M$ by the fact that for $\epsilon$ small enough, $x_M$ is close to $x^\ast \in \argmin f$. \\
\textbf{2)} \underline{\textit{The hybrid norm function}} $f : \mathbb{R}^2 \to \mathbb{R}$, defined as
\begin{align}\label{F2}\tag{F2}
f(x_1, x_2) = \sqrt{1 + x_1^2} + \sqrt{1 + x_2^2} \, \, , \, \forall (x_1, x_2) \in \mathbb{R}^2 \, .
\end{align}
The gradient of the objective function has the form $\nabla f(x_1, x_2) = \begin{pmatrix}
\dfrac{x_1}{\sqrt{1 + x_1^2}} \\ 
\dfrac{x_2}{\sqrt{1 + x_2^2}}
\end{pmatrix}$. On the other hand, the Hessian matrix associated to $f$ is given by $\nabla^2 f(x_1, x_2) = 
\begin{pmatrix}
\dfrac{1}{(1+x_1^2)^{3/2}} & 0 \\
0 & \dfrac{1}{(1+x_2^2)^{3/2}}
\end{pmatrix}$.
Taking $z = \begin{pmatrix}
z_1 \\
z_2
\end{pmatrix}$
we obtain that $z^T \nabla^2 f(x_1, x_2) z = z_1^2 \dfrac{1}{(1+x_1^2)^{3/2}} + z_2^2 \dfrac{1}{(1+x_2^2)^{3/2}} \geq 0$, so the Hessian satisfies $\nabla^2 f \succeq 0$. As for the case of \ref{F1}, we use Theorem 7.12 from \cite{Beck} and we reach the convexity of \ref{F2}. Furthermore, computing the solution of the equation $\nabla f(x_1, x_2) = \begin{pmatrix}
0 \\
0
\end{pmatrix}$ and using the fact from above that the objective function is convex, we obtain that  the set of stationary points is given by $\argmin f = \lbrace (0, 0) \rbrace$ (this follows from Theorem 7.9 of \cite{Beck}, namely the necessity and sufficiency of stationary points). By the fact that we use the $l_2$ norm on the gradient, we aim at evaluating the subordinate matrix norm. It is easy to observe that $\| \nabla^2 f(x_1, x_2) \|_1 = \| \nabla^2 f(x_1, x_2) \|_\infty = \max \Big\lbrace \dfrac{1}{(1+x_1^2)^{3/2}}, \dfrac{1}{(1+x_2^2)^{3/2}} \Big\rbrace \leq 1$. Since $\nabla^2 f(x_1, x_2) \in \mathbb{R}^{2 \times 2}$, by elementary algebra computations, we have that $\| \nabla^2 f(x_1, x_2) \|_2 \leq \sqrt{2} \| \nabla^2 f(x_1, x_2) \|_\infty$, so the tightest Lipschitz constant of the gradient is $L = \sqrt{2}$. At the same time, since the $\argmin$ set contains just one global minimum point, we use the stopping criteria $\| f(x_n) - f(x^\ast) \|_2 \leq \epsilon$, where $x^\ast$ is the unique element of $\argmin f$. Also, we consider the \ref{Energy} for the Lyapunov functional. \\
%\vspace*{0.25cm}
\begin{figure}[ht!]
\centering
\hspace*{-1.6cm}
\begin{subfigure}{.5\textwidth}
  \centering
  \includegraphics[width=9cm]{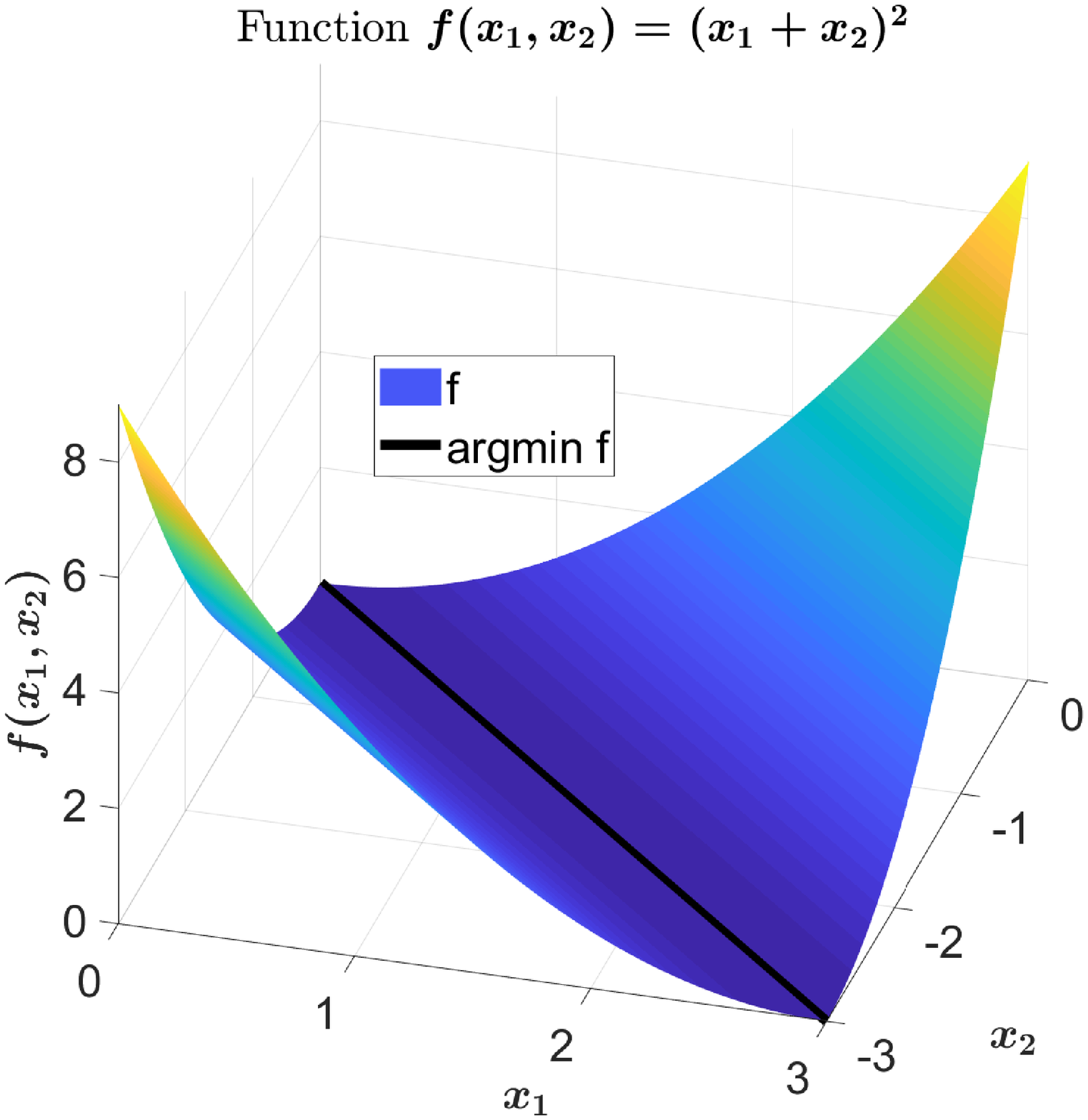}
  \caption{The representation of \ref{F1} along with \\ the set of minimum points $\argmin f$}
  \label{fig:1_functions}
\end{subfigure}%
\hspace*{0.8cm}
\begin{subfigure}{.5\textwidth}
  \centering
  \includegraphics[width=9cm]{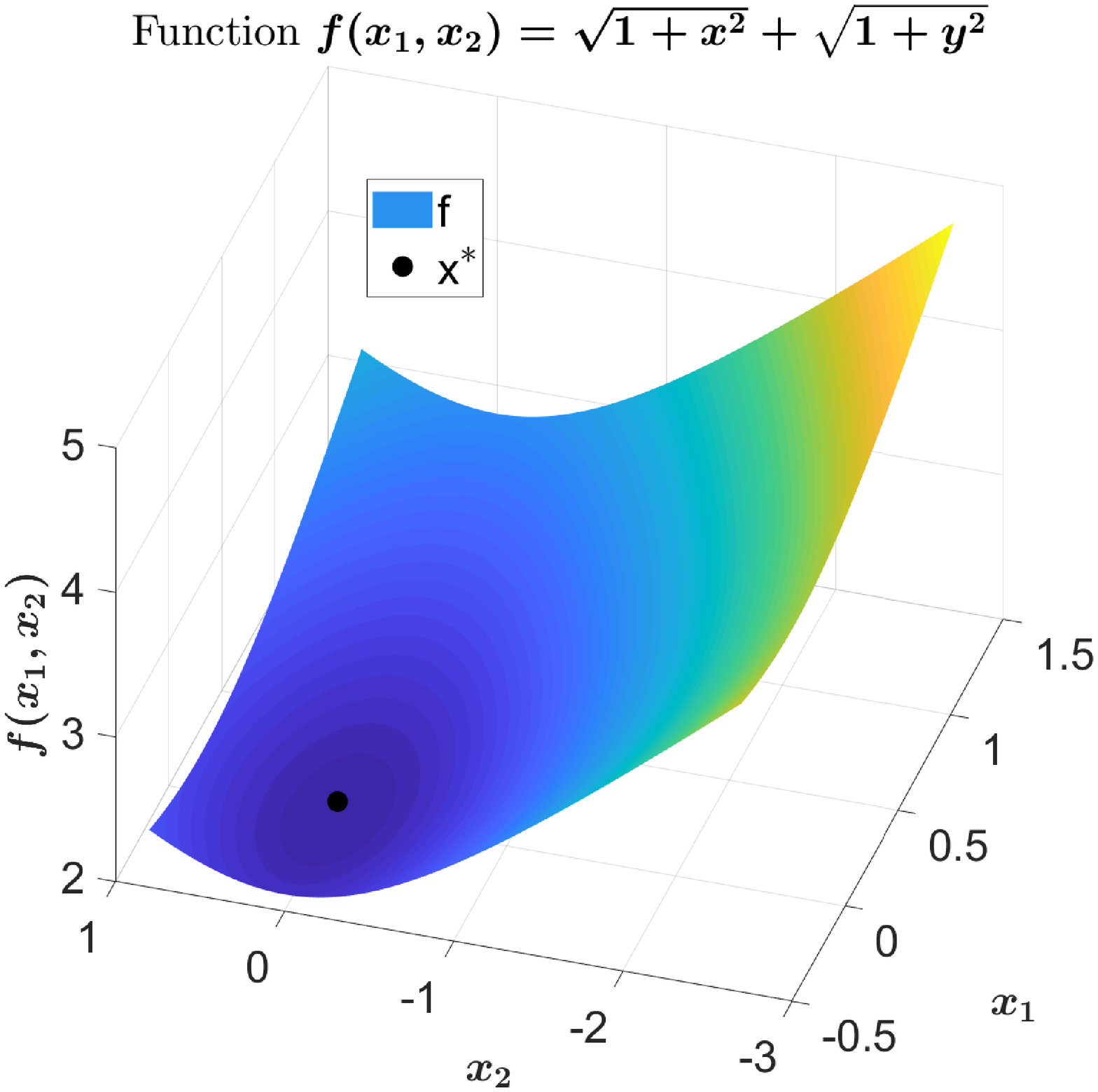}
  \caption{The representation of \ref{F2} along with \\ the set of minimum points $\argmin f$}
  \label{fig:2_functions}
\end{subfigure}
\caption{The convex objective functions \ref{F1} and \ref{F2}}
\label{fig:functions}
\end{figure}
Following the proof of Theorem \ref{T20}, we consider $N_1 = \alpha - 1$ and $N_2$ given by $G_n$, $H_n$ and $I_n$ as in the proof of the aforementioned theorem. Now we turn our attention to some choices of the parameters that appear in the three examples from the previous section, concerning the index $N^\prime$. However, we specify that even though $G_n$, $H_n$ and $I_n$ depend on the Lipschitz constant of the gradient $\nabla f$ (this is also valid for the case of Nesterov's method \ref{AGM2}), the empirical convergence is in fact determined by the index $N^\prime$ that contain some predetermined values of the inherent parameters. 

\subsection{Some discussions over the algorithms parameters}\label{S42}

In this subsection we consider some situations over the parameters referring to Example \ref{E24}, Example \ref{E25} and Example \ref{E26}. For each of them, we analyze in an empiric manner the key role of the underlying parameters over the convergence of the optimization algorithms of type \ref{LT-S-IGAHD}. \\
Likewise, for two parameters $A$ and $B$, we use the heuristic notations "$A \ll B$" for the situation when $A$ is much smaller than $B$ and "$A \gg B$" for the case $A$ is much larger than $B$, respectively. Therefore, for our analysis, we consider a heuristic analysis for our algorithms. \\ \\
\textbf{The case involving Example \ref{E24})} We aim at considering two subcases: \\
\textbf{A.)} When $\dfrac{\mu}{s} \ll 1$. This is, in fact, the practical case when the parameters have small value. Because, in general, $s$ is relatively small, this allow us  to take $\mu \ll 1$. For this, we consider the following situations:
\begin{itemize}
	\item[\textbf{A1.)}] When $a + \dfrac{1}{4} < b$. Now, if $a \gg 1$ ($a$ has a big value), then we can take $a \gg \left( \dfrac{\mu}{s} \right)^2$, so $N^\prime \ll 0$. On the other hand, when $a \ll 1$, then we can take $a \approx\left( \dfrac{\mu}{s} \right)^2$, so $N^\prime > 0$, i.e. $N^\prime = O \left( \mu L \right)$, but $\mu \ll 1$ annihilates the Lipschitz constant $L$ in the index value $N^\prime$.   
	\item[\textbf{A2.)}] When $b \leq a + \dfrac{1}{4}$. If $b \gg 1$, then $a \gg 1$, so we may take $a \gg \left( \dfrac{\mu}{s} \right)^2$. This is in fact the same situation as before when $a$ can take a big value, because the difference under the square root, namely $a-b$, can be close to $0$, in the sense that for $b \approx a$, $N^\prime \approx \dfrac{1}{2} - b$, i.e. $N^\prime \approx 0$. On the other hand, for the case when $b \ll 1$ and $a \gg 1$, we find that $N^\prime \approx \dfrac{\sqrt{\Delta}}{2}$, hence $N^\prime = O(\sqrt{a})$, if we can take $a$ large enough such that $a \gg \left( \dfrac{\mu}{s} \right)^2$.
\end{itemize}
\textbf{B.)} When $\dfrac{\mu}{s} \gg 1$. By the fact that in many situations $s \ll 1$, we can take $\mu > 1$ or $\mu \gg 1$. For this, we consider the following situations:
\begin{itemize}
\item[\textbf{B1.)}] When $a + \dfrac{1}{4} < b$. If $a \gg 1$, then $b \gg 1$. We can take $a \approx \left( \dfrac{\mu}{s} \right)^2$, so $N^\prime = O(\mu L)$, which can be eventually large. On the other hand, if $b \gg 1$ and $a \ll 1$, we can take $a \ll \left( \dfrac{\mu}{s} \right)^2$, so $N^\prime = O \left( \left( \dfrac{\mu}{s} \right)^2 \right)$, i.e. $N^\prime \gg 1$.
\item[\textbf{B2.)}] When $b \leq a + \dfrac{1}{4}$. If $a \gg 1$, thus $b \gg 1$, then we are in a similar situations as presented above, since we can take $a$ close to $b$ such that the difference $a-b$ in the square root in the definition of $N^\prime$ is close to $0$. On the other hand, if $a \gg 1$ and $b \ll 1$, then we can take $a \approx \left( \dfrac{\mu}{s} \right)^2$, thus $N^\prime = O \left( \sqrt{a} \right)$, which is equivalent to $N^\prime \gg 1$.
\end{itemize}
\vspace*{0.5cm}
\textbf{The case involving Example \ref{E25})} For this example, we consider the only valid situation when $2 \sqrt{s} > \beta$: \\
Since the stepsize $s$ has small values in many practical applications, we must consider $\beta < 1$ or $\beta \approx 0$. For small values of $\beta$, we observe that $N^\prime = O \left( \dfrac{ \dfrac{\mu}{\sqrt{s}} - 2 \sqrt{s} b + \sqrt{\Delta} }{2(2 \sqrt{s}-\beta)} \right)$. Now, we consider the following subcases: \\
\textbf{C1.)} When $\dfrac{\mu}{\sqrt{s}} \gg 1$ and $b$ is small or has moderate values, we have that 
$$N^\prime = O \left( \dfrac{ \dfrac{\mu}{\sqrt{s}} + \sqrt{(2 \sqrt{s}-\beta) \dfrac{\mu}{\sqrt{s}} + \left( \dfrac{\mu}{\sqrt{s}} \right)^2}}{2(2\sqrt{s}-\beta)} \right) \gg 1.$$
\textbf{C2.)} If we are in one of the cases when $\dfrac{\mu}{\sqrt{s}} < 1$ or $\dfrac{\mu}{\sqrt{s}} \ll 1$, then, by the analysis made above, we obtain that $N^\prime \approx 0$, since $\mu < \sqrt{s}$ or $\mu \ll 1$. \\
\textbf{C3.)} If $b \gg 0$ and dominates the other parameters, then it is easy to see that $N^\prime = O \left( (2 \sqrt{s} - \beta) \beta \sqrt{b} \right)$, i.e. $N^\prime \gg 1$. Furthermore, in order to have $N^\prime \approx 0$, then we must take $\beta \ll b$. \\ \\ \\
\textbf{The case involving Example \ref{E26})} This is similar with the first example involving the non-negative coefficients $a$, $b$ and $\mu$. For brevity, we consider the following cases for our heuristic interpretation: \\
\textbf{D.)} When $\dfrac{\mu}{s} \ll 1$. Since many times $s \ll 1$, we have that $\mu \ll 1$. Here, we consider the following subcases:
\begin{itemize}
\item[\textbf{D1.)}] When $a \leq b$. If $a \ll 1$ and $b \gg 1$, i.e. $a$ small and $b$ has large values, we find that $N^\prime \approx \sqrt{3 + \left( \dfrac{\mu}{s} \right)^2} = O \left( \left( \dfrac{\mu}{s} \right)^2 \right) \approx 0$. On the other hand, if $a \gg 1$ and $b \gg 1$, we obtain that $N^\prime \approx 0$.
\item[\textbf{D2.)}] When $b \leq a$. If $a \gg 1$, then $b \gg 1$ and, as before, $N^\prime \approx 0$. On the other hand, if $a \gg 1$ and $b \ll 1$, it follows that $N^\prime \approx \sqrt{3 + \left( \dfrac{\mu}{s} \right)^2}$, i.e. the index $N^\prime$ is small. 
\end{itemize}
\textbf{E.)} When $\dfrac{\mu}{s} \gg 1$, then $\sqrt{(sL)^2 + 1 + \left( \dfrac{\mu}{s} \right)^2} \gg 0$. Likewise, in practical situations we are obliged to take $\mu \gg 1$. We consider the following subcases:
\begin{itemize}
\item[\textbf{E1.)}] When $a \leq b$. In order to obtain the optimal case when $N^\prime \leq 0$, we can take $a \gg 0$, more precisely $a \geq \sqrt{3 + \left( \dfrac{\mu}{s} \right)^2} \geq \sqrt{(sL)^2 + 1 + \left( \dfrac{\mu}{s} \right)^2}$. It is trivial to see that since $a$ is large, then $b$ is also large.
\item[\textbf{E2.)}] When $b \leq a$. As in the previous case, we consider the optimal case when $N^\prime \leq 0$, so we can take easily take $b \geq \sqrt{3 + \left( \dfrac{\mu}{s} \right)^2}$. In this situation $b \gg 1$, so $a \gg 1$.
\end{itemize}

\subsection{A computational perspective over the optimization algorithms}\label{S43}

\begin{figure}[ht!]
\centering
\begin{subfigure}{.5\textwidth}
  \centering
  \includegraphics[width=8.5cm]{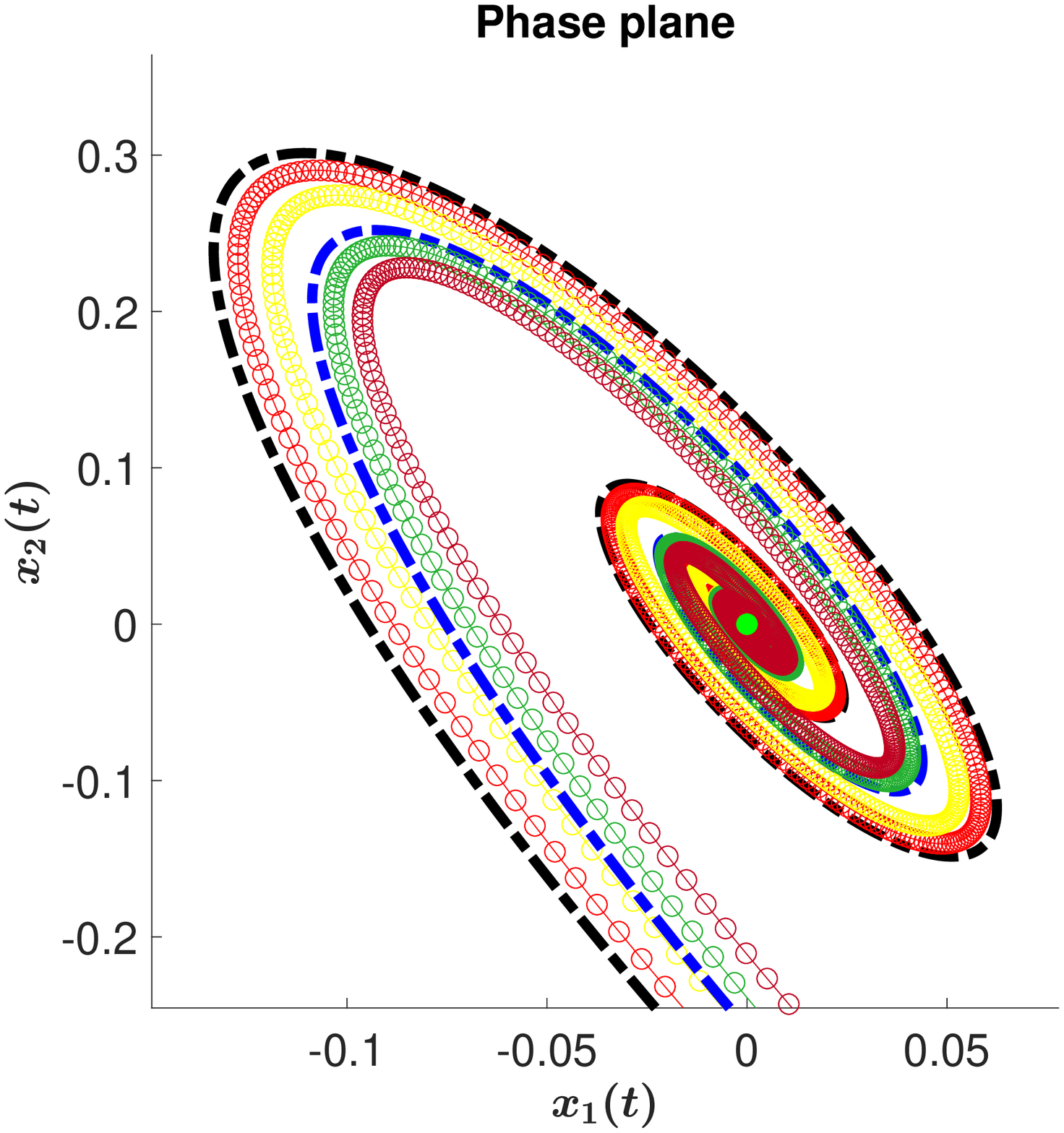}
  \caption{Zoomed trajectories of IGAHD (from \cite{Attouch_Chbani}) and \ref{LT-S-IGAHD} type algorithms from Example \ref{E25}}
  \label{fig:1_phaseplane}
\end{subfigure}%
\hspace*{0.8cm}
\begin{subfigure}{.5\textwidth}
  \centering
  \includegraphics[width=8.5cm]{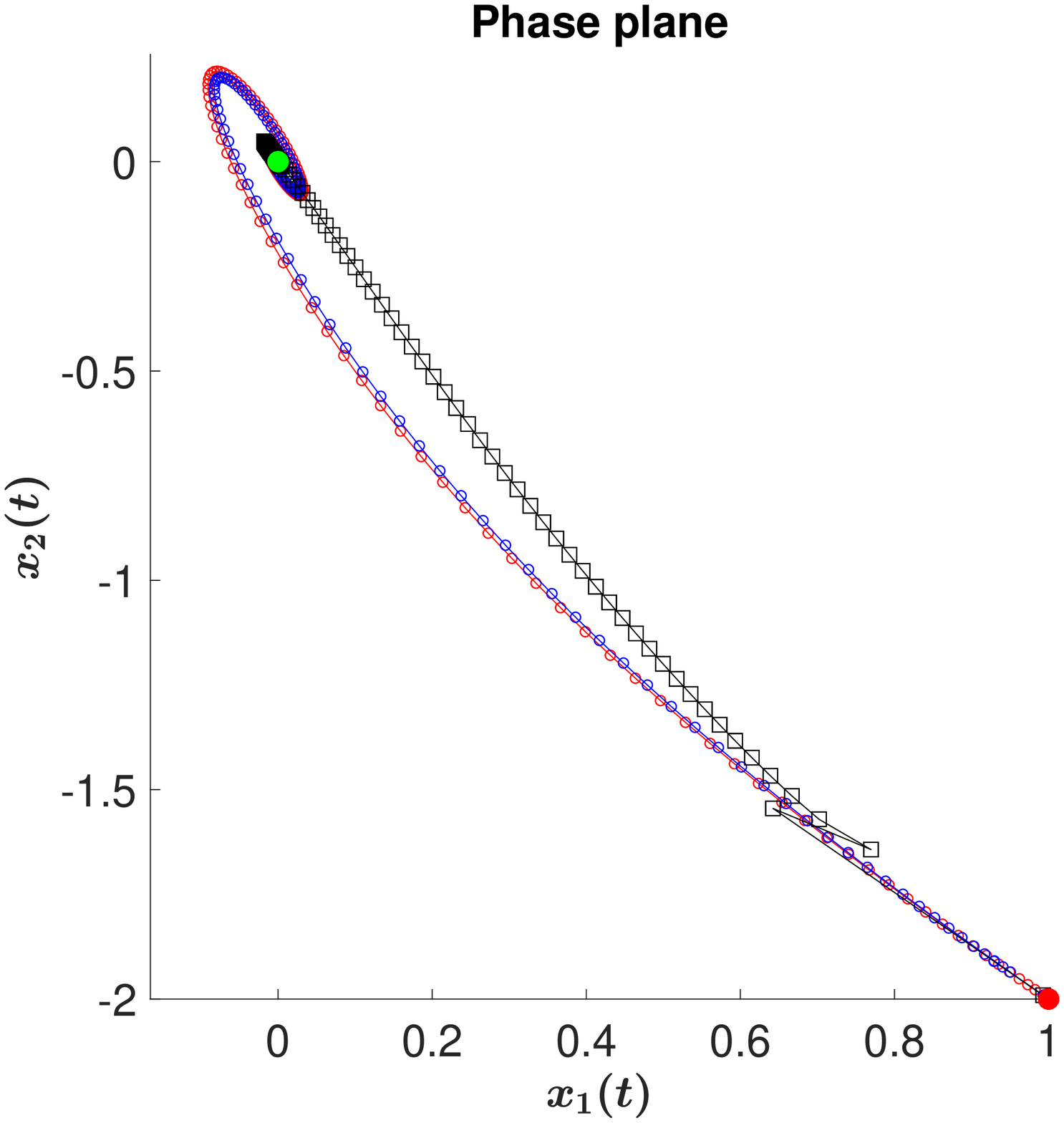}
  \caption{The convergence of the trajectories for \ref{LT-S-IGAHD} of Example \ref{E26}}
  \label{fig:2_phaseplane}
\end{subfigure}
\caption{The phase plane trajectories for \ref{LT-S-IGAHD}-type algorithms for the function \ref{F2}}
\label{fig:phaseplane}
\end{figure}

In the present subsection, we consider some numerical experiments over the algorithms of type \ref{LT-S-IGAHD}, taking into account the theoretical rate of convergence in the convex objective function obtained in the proof of Theorem \ref{T20}. At the same time, we precisely consider the cases and the subcases presented in Subsection \ref{S42}. \\
In Figure \ref{fig:1_phaseplane}, we have considered the convex objective function \ref{F2} along with the initial condition 
$x_0 = \begin{pmatrix}
1 \\
-2
\end{pmatrix}$
Also, we have set $\alpha = 3$ and the tolerance $\epsilon = 10^{-10}$. We have the behavior of the discretization schemes along the trajectories when the stepsize is taken as $s = 0.001$ in order to simulate the underyling dynamical systems. In the Subfigure \ref{fig:1_phaseplane},  we have considered the IGAHD algorithm $\beta = 10^{-2}$ and $b=1$ (black color) and the case when $\beta = 10^{-1}$ and $b=1$ (blue color). Furthermore, we have considered also the \ref{LT-S-IGAHD} from Example \ref{E25}, by the fact that when $\mu \approx 0$, then the algorithms becomes \ref{IGAHD}. We have chosen the following parameters : $\mu = 10^{-1}$, $\beta = 10^{-2}$, $b=1$ (red color), $\mu = 10^{-1}$, $\beta = 10^{-1}$, $b=1$ (saddle brown color, $[139, 69, 19]$ in RGB), $\mu = 0.25$, $\beta = 10^{-2}$, $b=1$ (yellow color) and $\mu = 0.25$, $\beta = 10^{-1}$, $b=1$ (rustic red color, $[75, 0, 13]$ in RGB). With this choices of coefficients, we observe that our algorithms have trajectories very close to the ones of IGAHD. This is valided by the theory where when $\mu$ is taken small, then the only difference between the algorithms gave in Example \ref{E25} and IGAHD is between $\nabla f(x_n)$ and $\nabla f(x_{n-1})$. On the other hand, in Subfigure \ref{fig:2_phaseplane}, we have considered the algorithm from Example \ref{E26}, with the following choices of the algorithm's coefficients: $a=2$, $\mu=10^{-15}$, $b=0.1$ (red color), $a=2$, $\mu=10^{-1}$, $b=0.1$ (blue color) and finally, $a=25$, $\mu=2$, $b=3$ (black color). In this situation we observe the phenomenon of the high amplification of the coefficients, namely when the parameters have large values, the trajectories have a high acceleration toward the critical point. From this subfigure and also from our preliminary numerical simulations, we have observed that there are some optimal choices of the parameters that lead to a faster convergence, that is given by a low number of iterations for a given tolerance.         

\begin{figure}[ht!]
\centering
\begin{subfigure}{.5\textwidth}
  \centering
  \includegraphics[width=8.5cm]{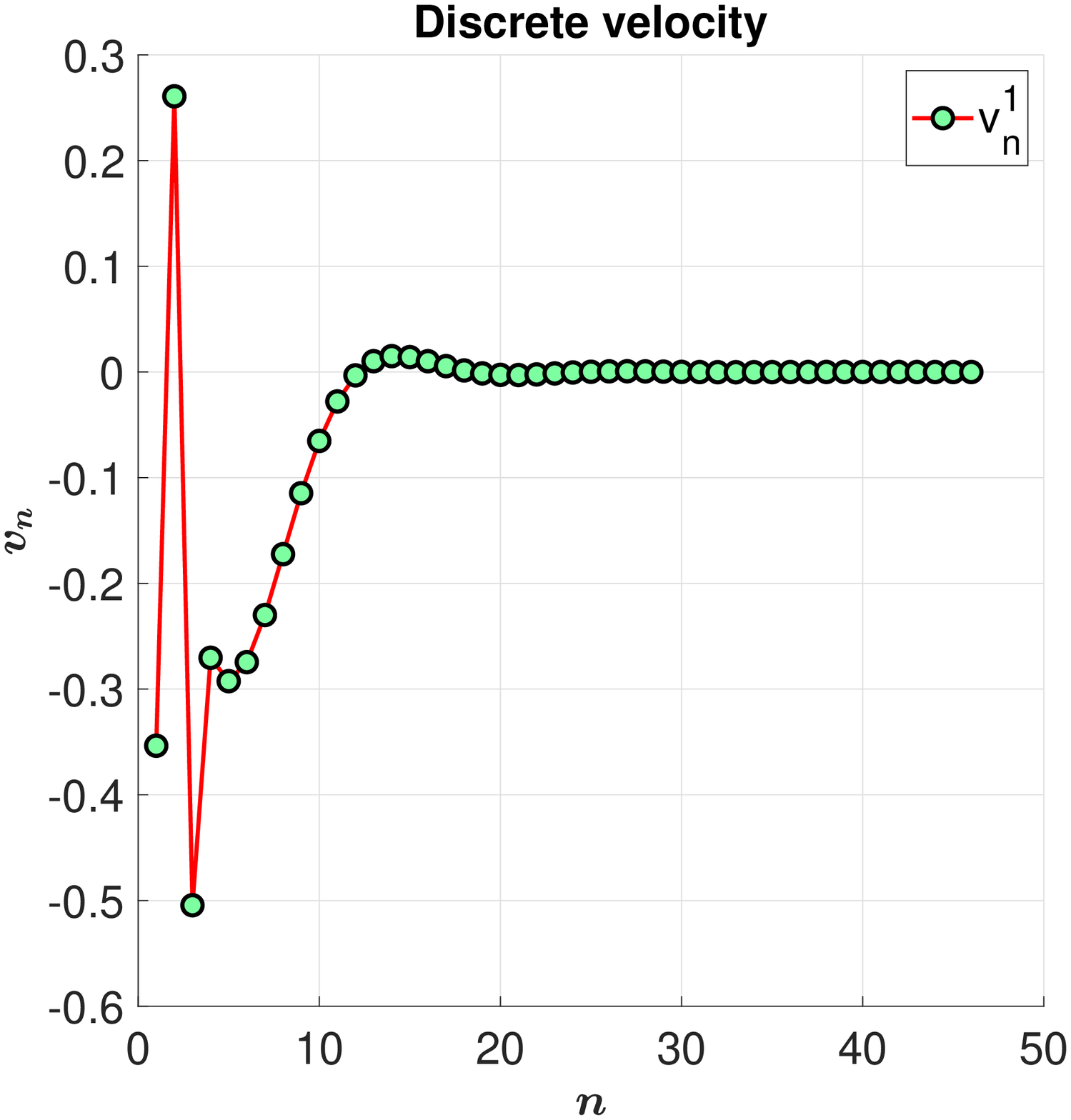}
  \caption{The case when the stepsize has the value $s=0.25$}
  \label{fig:1_velocity}
\end{subfigure}%
\hspace*{0.8cm}
\begin{subfigure}{.5\textwidth}
  \centering
  \includegraphics[width=8.5cm]{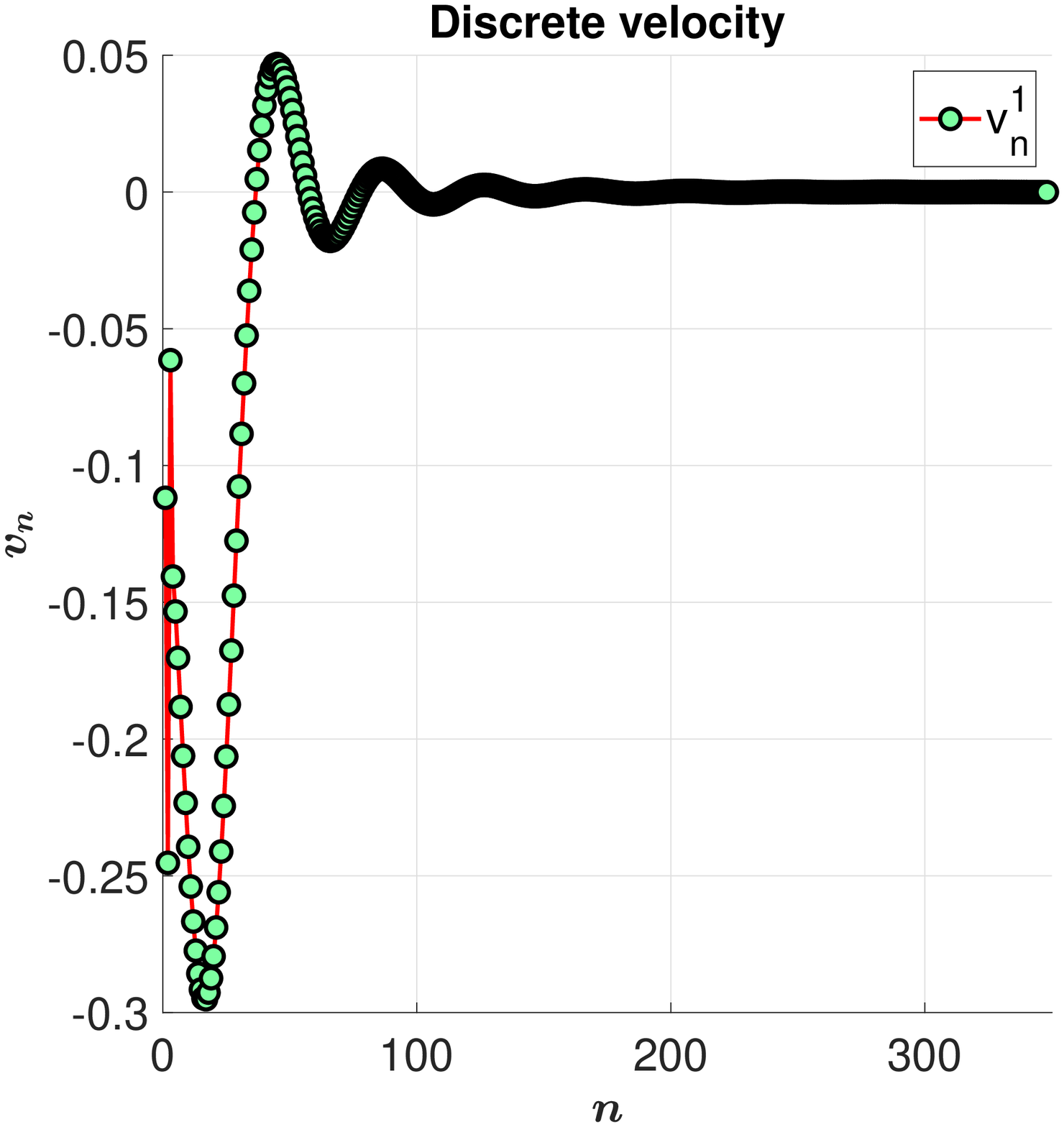}
  \caption{The case when the stepsize has the value $s=0.025$ }
  \label{fig:2_velocity}
\end{subfigure}
\caption{The first component of the discrete velocity of algorithms \ref{LT-S-IGAHD} from Example \ref{E25} for the function \ref{F2} }
\label{fig:velocity}
\end{figure}

Now, we focus on the discrete velocity $v_n = (x_n - x_{n-1})/\sqrt{s}$, where $\sqrt{s} = h > 0$. This is considered in Figure \ref{fig:velocity}. We have taken $\epsilon = 10^{-10}$, $\alpha = 3$ and 
$x_0 = \begin{pmatrix}
1 \\
-2
\end{pmatrix}$
For the algorithm \ref{LT-S-IGAHD} from Example \ref{E25}, we have set the values for the parameters as follows: $\mu = 0.5$, $\beta=10^{-5}$ and $b=2$. In Subfigre \ref{fig:1_velocity} we have plotted the first component of velocity in $\mathbb{R}^2$ for $s = 0.25$, while in Subfigure \ref{fig:2_velocity} we have taken a smaller stepsize $s = 0.025$. We observe that, in the beginning, the velocity presents high oscillations and then it stabilizes through the value $0$, and this is related to the fact that the velocity of the underlying evolution equation has the same behavior. The same stands true also for the second component of velocity and for other choices of the algorithm's parameters. Moreover, we observe that the stabilization process acts much faster than the convergence of the numerical scheme, when the tolerance is of low value. \\
Second to last, we turn our attention to Table \ref{tab:Ex24Fct1}, Table \ref{tab:Ex24Fct2}, Table \ref{tab:Ex26Fct1} and Table \ref{tab:Ex26Fct2}. Here, we have validated all of the cases from Subsection \ref{S42} involving both convex functions \ref{F1} and \ref{F2}, respectively. First thing to note is that \textit{error} represents the difference in two consecutive values in the objective function that is less than the given tolerance $\epsilon$ and along with this, it satisfies and additional stopping criteria, namely $n > N$, where $N = \max \lbrace N_1, N_2, N^\prime \rbrace$. We observe that our heuristic analysis stands very close to the simulations that we have made. For this, observe the values of $N^\prime$ given for each of the cases. An important difference is in the case $A_2$ where the value has an error by the fact that the coefficient $a$ is not much larger than $1$. Furthermore, another difference is for the case $B_2$ when we do not have taken $a-b \approx 0$ and this shows the fact that the heuristic interpretation that we gave leads to instability in the coefficients, in the sense that the values must take extreme values (see also the situation $B_1$). This is linked to the fact that in Subsection \ref{S42} we have considered an asymptotic analysis when some coefficients have very small values while others have very large values. The cases $A_2$, $B_1$ and $B_2$ justifies the difference between simulations and the heuristic analysis. On the other hand, for all of the other cases presented in the four tables, $N^\prime$ is very close to the actual theoretical value. Likewise, it is interesting to note that the empirical values $N_2$ and $N^\prime$ have an empiric relation, namely in the case of moderate or of high value coefficients, $N^\prime$ can be much greater than $N_2$, so in this sense, in the majority of cases, the index $N^\prime$ has a large influence over the convergence of the algorithms of type \ref{LT-S-IGAHD}. Nevertheless, there are some cases when $N_2$ has a larger influence over the rate of convergence, i.e. when the coefficients have extremely have numerical values, then $N_2$ increases (with respect to the values presented above in the other cases) and $N^\prime$ decreases. Even though $N^\prime$ is larger than $N_2$, the difference is lower than in the previous cases. Finally, we point out that in our numerical simulations, the coefficient $G_n $ is strictly positive, when $n > \max \lbrace N_1, N^\prime \rbrace$. This remains valid for all of the cases and for all of the three examples that we gave and it is consistent with the arguments from the third step of the proof of Theorem \ref{T20}. 

\begin{table}[!ht]
    \centering
    \captionof{table}{Analysis of the cases $A_1$, $A_2$, $B_1$ and $B_2$ of Subsection \ref{S42}, \\ for the Example \ref{E24} and for the function \ref{F1}}
    \label{tab:Ex24Fct1}
    \begin{tabular}{|| c || c | c | c | c | c | c | c || c ||}
\hline\hline
\diagbox[width=7em]{$Cases$}{$Params.$} & $error$ & $\varepsilon$ & $\mu$ & $a$ & $b$ & $N_2$ & $N^\prime$ & $N$  \\ \hline
$A_1$ & 1.22e-11 & 1e-10 & 1e-02 & 4.0 & 10.0 & 3.91 & -3.56 & 3.91 \\ \hline
$A_1$ & 3.31e-11 & 1e-10 & 1e-02 & 1e-02 & 10.0 & 3.99 & 0.43 & 3.99 \\ \hline
$A_2$ & 2.82e-11 & 1e-10 & 1e-02 & 10.0 & 4.0 & 3.80 & -1.00 & 3.80 \\ \hline
$A_2$ & 1.50e-12 & 1e-10 & 1e-05 & 3.0 & 1e-01 & 3.92 & 2.17 & 3.92 \\ \hline
$B_1$ & 0.0 & 1e-10 & 1.0 & 100.0 & 102.0 & 3.74 & 11.63 & 11.63 \\ \hline
$B_1$ & 0.0 & 1e-10 & 2.0 & 5e-01 & 6.0 & 3.48 & 422.45 & 422.45 \\ \hline
$B_2$ & 0.0 & 1e-10 & 1.5 & 2.0 & 1.75 & 3.58 & 240.29 & 240.29 \\ \hline
$B_2$ & 0.0 & 1e-10 & 1.5 & 225 & 1e-02 & 8.52 & 17.29 & 17.29 \\ \hline\hline
\end{tabular}
\end{table}

\begin{table}[!ht]
    \centering
    \captionof{table}{Analysis of the cases $A_1$, $A_2$, $B_1$ and $B_2$ of Subsection \ref{S42}, \\ for the Example \ref{E24} and for the function \ref{F2}}
    \label{tab:Ex24Fct2}
    \begin{tabular}{|| c || c | c | c | c | c | c | c || c ||}
\hline\hline
\diagbox[width=7em]{$Cases$}{$Params.$} & $error$ & $\varepsilon$ & $\mu$ & $a$ & $b$ & $N_2$ & $N^\prime$ & $N$  \\ \hline
$A_1$ & 7.03e-11  & 1e-10 & 1e-02 & 4.0 & 10.0 & 3.38 & -3.91 & 3.38 \\ \hline
$A_1$ & 1.92e-12 & 1e-10 & 1e-02 & 1e-02 & 10.0 & 3.38 & 0.08 & 3.38 \\ \hline
$A_2$ & 3.47e-11 & 1e-10 & 1e-02 & 10.0 & 4.0 & 3.37 & -1.0 & 3.37 \\ \hline
$A_2$ & 7.53e-11 & 1e-10 & 1e-05 & 3.0 & 1e-01 & 3.38 & 2.17 & 3.38 \\ \hline
$B_1$ & 0.0 & 1e-10 & 1 & 100 & 102 & 3.50 & 4.04 & 4.04 \\ \hline
$B_1$ & 0.0 & 1e-10 & 2.0 & 5e-01 & 6.0 & 3.43 & 407.54 & 407.54 \\ \hline
$B_2$ & 0.0 & 1e-10 & 1.5 & 2.0 & 1.75 & 3.5 & 229.04 & 229.04 \\ \hline
$B_2$ & 0.0 & 1e-10 & 1.5 & 225 & 1e-02 & 4.46 & 15.50 & 15.50 \\ \hline\hline
\end{tabular}
\end{table}

\begin{table}[!ht]
    \centering
    \captionof{table}{Analysis of the cases $D_1$, $D_2$, $E_1$ and $E_2$ of Subsection \ref{S42}, \\ for the Example \ref{E26} and for the function \ref{F1}}
    \label{tab:Ex26Fct1}
    \begin{tabular}{|| c || c | c | c | c | c | c | c || c ||}
\hline\hline
\diagbox[width=7em]{$Cases$}{$Params.$} & $error$ & $\varepsilon$ & $\mu$ & $a$ & $b$ & $N_2$ & $N^\prime$ & $N$  \\ \hline
$D_1$ & 3.37e-11  & 1e-10 & 0.0 & 0.25 & 3.5 & 3.18 & 0.83 & 3.18 \\ \hline
$D_1$ & 9.79e-11  & 1e-10 & 1e-3 & 1.25 & 5.5 & 3.18 & -0.17 & 3.18 \\ \hline
$D_2$ & 3.28e-11  & 1e-10 & 1e-3 & 5.5 & 1.25 & 3.19 & -0.17 & 3.19 \\ \hline
$D_2$ & 1.11e-11  & 1e-10 & 1e-3 & 3.5 & 0.25 & 3.19 & 0.83 & 3.19 \\ \hline
$E_1$ & 0.0  & 1e-10 & 2.0 & 21.0 & 24.0 & 5.82 & -0.97 & 5.82 \\ \hline
$E_2$ & 0.0  & 1e-10 & 2.0 & 24.0 & 21.0 & 6.09 & -0.97 & 6.09 \\ \hline\hline
\end{tabular}
\end{table}

\begin{table}[!ht]
    \centering
    \captionof{table}{Analysis of the cases $D_1$, $D_2$, $E_1$ and $E_2$ of Subsection \ref{S42}, \\ for the Example \ref{E26} and for the function \ref{F2}}
    \label{tab:Ex26Fct2}
    \begin{tabular}{|| c || c | c | c | c | c | c | c || c ||}
\hline\hline
\diagbox[width=7em]{$Cases$}{$Params.$} & $error$ & $\varepsilon$ & $\mu$ & $a$ & $b$ & $N_2$ & $N^\prime$ & $N$  \\ \hline
$D_1$ & 4.21e-12 & 1e-10 & 0.0 & 0.25 & 3.5 & 3.23 & 0.76 & 3.23 \\ \hline
$D_1$ & 9.18e-11  & 1e-10 & 1e-03 & 1.25 & 5.5 & 3.23 & -0.24 & 3.23 \\ \hline
$D_2$ & 5.68e-11 & 1e-10 & 1e-03 & 5.5 & 1.25 & 3.23 & -0.24 & 3.23 \\ \hline
$D_2$ & 7.18e-11  & 1e-10 & 1e-03 & 3.5 & 0.25 & 3.23 & 0.76 & 3.23 \\ \hline
$E_1$ & 0.0  & 1e-10 & 2.0 & 21.0 & 24.0 & 4.14 & -0.97 & 4.14 \\ \hline
$E_2$ & 0.0  & 1e-10 & 2.0 & 24.0 & 21.0 & 4.24 & -0.97 & 4.24 \\ \hline\hline
\end{tabular}
\end{table}

\clearpage
\section*{Conclusions and Discussions}\label{S5}
The novelty of the present paper consists in the following both theoretical and practical aspects regarding Hessian driven damping dynamical systems and also their numerical methods. We conclude the link between the discrete case and the continuous counterpart of the optimization algorithms, by pointing out our major contributions and some discussions that arise naturally from our research: 
\begin{itemize}
\item[\textbf{1)}] In Subsection \ref{S21}, we have shown that Nesterov's algorithm \ref{AGM2} is in fact a natural discretization of a Hessian-type dynamical system. Even though other authors (see \cite{ChenLuo} and \cite{Oberman}) have pointed out a link between continuous Hessian dynamical systems and Nesterov algorithm, we have introduced a rigorous perspective such that Nesterov method can be interpreted as a Lie-Trotter discretization of the underlying second evolution equation. This approach is totally new with respect to inertial algorithms, in the following sense: only recently (see \cite{ALP} and \cite{Jordan}) it was shown that Nesterov's method is a natural discretization of the \ref{Extended-AVD} dynamical system. Before that, no other constructions were available such that Nesterov algorithm could be interpreted as an explicit numerical scheme of a continuous dynamical system without an extra-gradient step, so our method represents a totally new approach to optimization. Furthermore, the use of the splitting technique is somewhat new to optimization problems, so that this article strengthens the relationship between distinct mathematical fields. 
\item[\textbf{2)}] We have presented a geometrical-type perspective upon optimization methods in the setting of minimization problems that involve convex objective functions (see Subsection \ref{S23}). Likewise, we gave intuitive explanations over the choice of our algorithms in Subsection \ref{S24}. At the same time, in Subsection \ref{S32} we have constructed Lemma \ref{L15} that represents a generalization of Lemma \ref{L14} that was given in \cite{Attouch_Chbani}. Also, our proof of Theorem \ref{T20} is based upon the proof from the same paper, namely \cite{Attouch_Chbani}, but in a much more general setting.
\item[\textbf{3)}] We have introduced new algorithms that contain the discrete version of the Hessian system, namely \ref{LT-S-IGAHD}. Also, the particular case given in Example \ref{E25} resembles the newly introduced algorithm IGAHD in \cite{Attouch_Chbani}. Also, in Example \ref{E24} and Example \ref{E26} we have gave general algorithms that contain three arbitrary parameters (under some suitable conditions) that have numerical benefits in practical applications (see Subsections \ref{S42} and \ref{S43}). On the other hand, we emphasize that our algorithms have been developed in an elegant manner based upon the combination of splitting methods with symplectic algorithms, similar to the simpler case of Nesterov inertial algorithm (see \cite{Jordan}). 
\item[\textbf{4)}] We have developed some new optimization algorithms based upon two gradient evaluations at each iteration. For each of these numerical schemes we have presented the associated continuous versions, such that these numerical schemes are constructed through various geometric principles like splitting and symplectic methods. The associated continuous dynamical systems \ref{DynSys-ARDM} and the system from Theorem \ref{T7} are worth studying by the fact that their discretizations are obtained in a natural manner, as in the case with Nesterov's method (see \cite{ALP}and \cite{Jordan}).
\end{itemize}

\section*{Acknowledgments}
I am gratefully indebted to Szilard Csaba L\' aszl\' o for his suggestions concerning the organization of this article and for the helpful comments regarding the dynamical systems and the optimization algorithms of the present research paper.

\end{document}